\documentclass[11pt]{article}

\usepackage[margin=1in]{geometry}
\usepackage{amsmath}
\usepackage{amssymb}
\usepackage{paralist}
\usepackage{graphicx}
\graphicspath{{Figures/}}
\usepackage{epsfig}
\usepackage{epstopdf}
\usepackage{algpseudocode}
\usepackage{booktabs}
\usepackage{multirow}
\usepackage{subcaption}
\usepackage{xcolor}
\usepackage{float}
\usepackage{placeins}
\usepackage{algorithm}

\usepackage[colorlinks=true]{hyperref}
\hypersetup{
    urlcolor=blue,
    citecolor=red
}

\allowdisplaybreaks

\newcommand{\subjclasstext}[1]{\par\noindent\textbf{2020 Mathematics Subject Classification.} #1.}
\newcommand{\keywordstext}[1]{\par\noindent\textbf{Key words and phrases.} #1.}

\newcommand{\be}{\begin{equation}}
\newcommand{\bea}{\left[ \begin{array}}
\newcommand{\ee}{\end{equation}}
\newcommand{\eea}{\end{array} \right]}
\newcommand{\ba}{\begin{array}}
\newcommand{\ea}{\end{array}}

\algrenewcommand\algorithmicrequire{\textbf{Input:}}
\algrenewcommand\algorithmicensure{\textbf{Output:}}

\DeclareMathOperator*{\argmin}{arg\,min}

\def\I0{\bigl[ \, I_m \, , \, 0_{m\times p} \, \bigr]}
\def\0I{\bigl[ \, 0_{n\times p} \, I_p \, \bigr]}

\newcommand{\bvec}{\left[ \begin{array}{c} }
\newcommand{\evec}{\end{array} \right]}

\newcommand{\bA}{{\bf A}}
\newcommand{\bB}{{\bf B}}

\newcommand{\bD}{{\bf D}}

\newcommand{\bI}{{\bf I}}

\newcommand{\bL}{{\bf L}}

\newcommand{\bP}{{\bf P}}
\newcommand{\bQ}{{\bf Q}}
\newcommand{\bR}{{\bf R}}

\newcommand{\bU}{{\bf U}}
\newcommand{\bV}{{\bf V}}
\newcommand{\bW}{{\bf W}}
\newcommand{\bX}{{\bf X}}

\newcommand{\bb}{{\bf b}}
\newcommand{\bc}{{\bf c}}
\newcommand{\bd}{{\bf d}}
\newcommand{\bbe}{{\bf e}}

\newcommand{\bg}{{\bf g}}

\newcommand{\bq}{{\bf q}}

\newcommand{\bu}{{\bf u}}
\newcommand{\bw}{{\bf w}}
\newcommand{\bx}{{\bf x}}
\newcommand{\by}{{\bf y}}

\newcommand{\bbR}{\mathbb{R}}

\title{Wavelet-based multilevel framework for $\ell_1$-regularized image deblurring}

\author{
Danyh Tolah$^{1}$
\quad Malena I. Espa\~nol$^{2}$
\quad Misha E. Kilmer$^{3}$
\\[1ex]
{\footnotesize $^1$College of Sciences and Human Studies, Prince Mohammad Bin Fahd University, Al-Khobar, Saudi Arabia}\\
{\footnotesize \href{mailto:dtolah@pmu.edu.sa}{dtolah@pmu.edu.sa}}\\[1ex]
{\footnotesize $^2$School of Mathematical and Statistical Sciences, Arizona State University, Tempe, AZ, United States}\\
{\footnotesize \href{mailto:malena.espanol@asu.edu}{malena.espanol@asu.edu}}\\[1ex]
{\footnotesize $^3$Department of Mathematics, Tufts University, Medford, MA, United States}\\
{\footnotesize \href{mailto:misha.kilmer@tufts.edu}{misha.kilmer@tufts.edu}}\\
}

\date{}

\begin{document}

\maketitle


\begin{abstract}
Solving large-scale $\ell_1$-regularized image deblurring problems
efficiently while preserving sharp edges remains a significant computational
challenge. We propose a wavelet-based multilevel framework that embeds three
iterative solvers, Iteratively Reweighted Least Squares, Split Bregman, and
Majorization-Minimization, within a multilevel V-cycle. Discrete wavelet
transforms define the interlevel transfer operators, and regularization
parameters are selected automatically by Generalized Cross Validation. Two
information transfer strategies are introduced and compared: one transfers
only the coarse solution to the fine level, while the other transfers
solver-specific auxiliary quantities. Numerical experiments demonstrate substantial computational savings
for IRLS, with speedups exceeding an order of magnitude, while MM and SB exhibit more modest computational differences. The experiments generally show that transferring auxiliary iterates performs best with Haar wavelets, whereas transferring only the solution performs best with Daubechies wavelets.
\end{abstract}

\subjclasstext{Primary: 65F10, 65F22}
\keywordstext{Image deblurring, multilevel methods, $\ell_1$ regularization, wavelets, generalized cross-validation parameter selection}

\bigskip

\section{Introduction}

Many scientific and engineering problems can be modeled as discrete ill-posed inverse problems~\cite{Hansenbk, hansen2010discrete}. A prototypical example is image deblurring, where the goal is to recover a sharp image from blurred and noisy data~\cite{austin2022image, hansen2006deblurring}. This problem can be written as
\begin{equation}\label{eq:original}
    \bA\bx \approx \bb = \bb^{\mathrm{true}} + \bbe,
    \qquad \text{with} \qquad
    \bA\bx^{\mathrm{true}} = \bb^{\mathrm{true}},
\end{equation}
where $\bA\in \mathbb{R}^{m^2\times m^2}$ is the blurring operator, $\bx^{\mathrm{true}}\in \mathbb{R}^{m^2}$ represents the unknown image of size $m\times m$, and $\bbe$ represents additive noise. For simplicity, we will consider square images with $m=2^p$, where $p$ is a nonnegative integer. Since the singular values of $\bA$ decay to zero, least-squares solutions are unstable and require regularization. In particular, $\ell_1$ regularization~\cite{alotaibi2023krylov, buccini2020lp, RudinOsherFatemi1992} is widely used for edge-preserving reconstruction, leading to the problem
\begin{equation}\label{eq:l1b}
    \min_{\bx} \left\{\|\bA\bx-\bb\|_2^2+\mu\|\bL\bx\|_1\right\},
\end{equation}
where $\mu>0$ is a regularization parameter and $\bL$ is a regularization operator. We will focus on $\bL$ being a discrete gradient operator.

Solving~\eqref{eq:l1b} efficiently remains computationally challenging for large-scale imaging problems. A common strategy is to reformulate~\eqref{eq:l1b} as a sequence of $\ell_2$ regularized least squares subproblems: Iteratively Reweighted Least Squares (IRLS) \cite{wohlberg2007iteratively} solves a sequence of weighted least squares problems with iteratively updated weights, while Split Bregman (SB)~\cite{GoldsteinOsher2009SplitBregman} and Majorization-Minimization (MM)~\cite{huang2017majorization, HunterLange2004, sweeney2025thesis} instead introduce auxiliary variables to obtain related quadratic subproblems. In all three cases, the original nonsmooth optimization problem is converted into a sequence of large-scale least squares subproblems, each of which remains computationally expensive. Krylov subspace methods, such as LSQR~\cite{paige1982lsqr}, are well-suited for these systems, but their computational cost grows rapidly with image resolution.

Multigrid and multilevel strategies provide an effective alternative by shifting much of the computational effort to smaller coarse-level problems while maintaining accuracy. Such approaches have been successfully applied to image deblurring and related inverse problems~\cite{Amdouni2024, austin2022image, bolten2024multigrid, buccini2020multigrid, Donatelli, espanol2010multilevel, espanol2014wavelet, Lauga2024Thesis, Lauga2024IMLFISTA, MorigiReichelSgallariShyshkov}. Wavelet-based multilevel methods are particularly attractive for imaging problems due to their ability to represent images sparsely across scales. Early work~\cite{espanol2020multilevel, espanol2010multilevel, espanol2014wavelet, espanol2009multilevel, klann2011wavelet} introduced multilevel methods in which discrete wavelet transforms define interlevel transfer operators, as distinct from work using wavelets as a regularization basis to enforce sparsity in the transform domain~\cite{belge2000wavelet,vonesch2009fast}.

More recently, a variety of multilevel approaches for $\ell_1$-regularized inverse problems have been developed. In~\cite{buccini2020multigrid}, framelet soft-thresholding is combined with an iterated Tikhonov smoother within a multigrid framework, promoting sparsity indirectly in the transform domain; low-pass framelet filters are used as transfer operators, and the noise level is assumed known for parameter selection. In~\cite{Lauga2024IMLFISTA}, a multilevel inertial forward-backward method (IML-FISTA) is proposed, where Moreau envelope smoothing is employed to ensure first-order coherence between levels. Iteratively reweighted edge-preserving regularization has also been investigated~\cite{gazzola2020inner}, with algebraic multigrid preconditioners proposed for the resulting sequences of $\ell_2$-regularized least squares problems~\cite{bolten2024multigrid}. Despite this progress, to the best of our knowledge, no previous wavelet-based multilevel framework unifies IRLS, SB, and MM within a single V-cycle with automatic parameter selection, and the interaction between wavelet family and information transfer strategy has not been systematically studied.

\paragraph{\textbf{Main contributions}}
The contributions of this paper are threefold. First, we develop a unified wavelet-based multilevel V-cycle for IRLS, SB, and MM, using low-pass wavelet-based transfer operators. Edge preservation is enforced explicitly through an $\ell_1$ penalty at every level, and regularization parameters are selected automatically via Generalized Cross Validation (GCV)~\cite{golub1979gcv} within each $\ell_1$ solver, without requiring prior knowledge of the noise level. Second, two multilevel information transfer strategies are introduced and compared: Approach~X, which transfers only the coarse solution to initialize the fine-level solve, and Approach~D, which instead transfers solver-specific auxiliary quantities. All subproblems within the $\ell_1$ solvers are solved with LSQR. Third, systematic numerical experiments across two- and three-level cycles, two wavelet families, and both information transfer strategies reveal a general interaction: Approach~D benefits from Haar wavelets, while Approach~X benefits from Daubechies wavelets, a trade-off that constitutes the central qualitative finding of the study.

\paragraph{\textbf{Paper organization}} The remainder of the paper is organized as follows. Section~\ref{sec:solvers} describes three iterative methods for the $\ell_1$-regularized problem. Krylov subspace projection and LSQR are covered in Section~\ref{sec:LSQR}. The wavelet-based multilevel method, including its computational implementation, is described in Section~\ref{sec:multilevel}. Numerical experiments are presented in Section~\ref{sec:num}. Conclusions are given in Section~\ref{sec:conclusion}.

\section{Iterative Methods}
\label{sec:solvers}

In this section, we review three iterative methods for solving~\eqref{eq:l1b}: IRLS, MM, and SB. All three reformulate the $\ell_1$-regularized problem as a sequence of $\ell_2$-regularized least squares subproblems of the form
\begin{equation}\label{eq:generic_sub}
    \min_{\bx} \left\{ \|\bA\bx - \bb\|_2^2 + \gamma^{(k)} \|\bD^{(k)}\bx - \bc^{(k)}\|_2^2 \right\},
\end{equation}
where $\gamma^{(k)} > 0$ is a regularization parameter and  $\bD^{(k)}$ and $\bc^{(k)}$ are a matrix operator and a shift vector, respectively, that differ by method. Each subproblem is solved via LSQR (Section~\ref{sec:LSQR}), and the parameter $\gamma^{(k)}$ is estimated by GCV, with the estimation frequency differing by method (see Sections~\ref{sub:IRLS}--\ref{sub:SB}).

\subsection{Iteratively Reweighted Least Squares}\label{sub:IRLS}

IRLS~\cite{Daubechies2010, wohlberg2007iteratively} solves a sequence of
weighted $\ell_2$ subproblems whose weights depend on the current iterate.
Here, we apply the iteratively reweighted norm construction to the
regularization vector $\bL\bx$. At iteration $k$, the method computes
$\bx^{(k+1)}$ by solving
\begin{equation}\label{eq:IRLS1}
\min_{\bx}\left\{\|\bA\bx-\bb\|_2^2+\mu^{(k)}\|\operatorname{diag}(\bw^{(k)})\bL\bx\|_2^2\right\}.
\end{equation}
This has the form~\eqref{eq:generic_sub} with $\gamma^{(k)}=\mu^{(k)}$,$\bD^{(k)}=\operatorname{diag}(\bw^{(k)})\bL$, and $\bc^{(k)}=\mathbf{0}$. For the $\ell_1$ regularization term, the weight vector $\bw^{(k)}$ is defined component-wise by
$$w_i^{(k)}=\left(\max\left\{\left|(\bL\bx^{(k)})_i\right|,\epsilon\right\}\right)^{-1/2}.$$
Here, $\epsilon>0$ is a small stabilization constant that prevents the weights from becoming unbounded when a component of $\bL\bx^{(k)}$ is close to zero. With this choice,
$$\|\operatorname{diag}(\bw^{(k)})\bL\bx\|_2^2 =\sum_i\frac{\left|(\bL\bx)_i\right|^2}{\max\left\{\left|(\bL\bx^{(k)})_i\right|, \epsilon\right\}}.$$
Thus, the weighted quadratic term provides an approximation to the $\ell_1$ regularization term. Large weights at small-gradient locations penalize those components more strongly, promoting sparsity in $\bL\bx$, while smaller weights at large-gradient locations relax the penalty and help preserve edges. The regularization parameter $\mu^{(k)}$ is estimated via GCV as described in Section~\ref{sec:LSQR}. Iterations terminate when
\begin{equation}\label{eq:stop}
    \mathrm{RC}(\bx^{(k+1)})=\frac{\|\bx^{(k+1)} - \bx^{(k)}\|_2}{\|\bx^{(k)}\|_2} < tol,
\end{equation}
or when $k = K_{max}$. For convergence analysis see~\cite{Daubechies2010}. Algorithm \ref{alg:IRLS} summarizes IRLS.

\begin{algorithm}
\caption{The IRLS Method}
\label{alg:IRLS}
\begin{algorithmic}[1]
\Require $\bA, \bb, \bL, \epsilon, \bx^{(0)}, tol, K_{max}$
\Ensure $\bx$
    \For{$k=0,1,\dots$ until $\textrm{RC}(\bx^{(k+1)}) < tol$ or $k = K_{max}$}
    \State $w_i^{(k)} = \left(\max\left\{|(\bL\bx^{(k)})_i|,\epsilon\right\}\right)^{-1/2}$
    \State Estimate $\mu^{(k)}$ via GCV applied to~\eqref{eq:IRLS1} (see Section~\ref{sec:LSQR})
    \State Solve~\eqref{eq:IRLS1} using LSQR to obtain $\bx^{(k+1)}$
    \EndFor
\end{algorithmic}
\end{algorithm}

\subsection{Majorization-Minimization}\label{sub:MM}

MM~\cite{HunterLange2004} replaces the $\ell_1$ term in \eqref{eq:l1b} at each iteration by a smooth quadratic surrogate that upper-bounds it and is tight at the current iterate. Using the fixed quadratic majorant of~\cite{huang2017majorization} with stabilization parameter $\varepsilon > 0$, the $k$-th iteration solves
\begin{equation} \label{eq:MMsolve}
    \min_{\bx}\left\{\|\bA\bx - \bb\|_2^2 + \mu^{(k)} \|\bL\bx - \bw^{(k)}\|_2^2\right\},
\end{equation}
which has the form~\eqref{eq:generic_sub} with $\gamma^{(k)} = \mu^{(k)}$, $\bD^{(k)}=\bL$, and the shift vector $\bc^{(k)} = \bw^{(k)}$, where
\begin{equation} \label{eq:wreg}
    \bw_i^{(k)} = \bu_i^{(k)} \left( 1 - \left(\frac{\varepsilon^2}{(\bu_i^{(k)})^2 + \varepsilon^2}\right)^{1/2}\right),
    \quad \bu^{(k)}= \bL\bx^{(k)},
\end{equation}
with all operations being component-wise. The shift $\bw^{(k)}$ acts as a spatially adaptive target for $\bL\bx$: near edges where $|\bu_i^{(k)}|$ is large, $\bw_i^{(k)} \approx \bu_i^{(k)}$ and the gradient penalty is relaxed, preserving the edge; in smooth regions where $|\bu_i^{(k)}|$ is small, $\bw_i^{(k)} \approx 0$ and a standard $\ell_2$ gradient penalty is recovered. Following~\cite{sweeney2025parameter}, $\mu^{(k)}$ is estimated via GCV (using the pre-computed generalized singular value decomposition (GSVD) of $\{\bA,\bL\}$) in the first three outer iterations and held fixed thereafter; fixing it avoids instability caused by the rapidly changing surrogate in early iterations~\cite{sweeney2025parameter}. The stopping criteria are the same as for IRLS.  Algorithm \ref{alg:MML2} summarizes MM for solving \eqref{eq:l1b}.

\begin{algorithm}[H]
\caption{The MM Method}
\label{alg:MML2}
\begin{algorithmic}[1]
\Require $\bA, \bb, \bL, \bx^{(0)}, \varepsilon, tol, K_{max}$
\Ensure $\bx$
    \For{$k=0,1,\dots$ until $\textrm{RC}(\bx^{(k+1)}) < tol$ or $k = K_{max}$}
    \State $\bu^{(k)} = \bL\bx^{(k)}$
    \State $\bw_i^{(k)} = \bu_i^{(k)}\!\left(1-\left(\frac{\varepsilon^2}{(\bu_i^{(k)})^2+\varepsilon^2}\right)^{\!1/2}\right)$
    \State Estimate $\mu^{(k)}$ via GCV for $k \le 3$; set $\mu^{(k)} = \mu^{(3)}$ for $k > 3$
    \State Solve~\eqref{eq:MMsolve} using LSQR to obtain $\bx^{(k+1)}$
   \EndFor
\end{algorithmic}
\end{algorithm}

\subsection{Split Bregman}\label{sub:SB}

SB~\cite{GoldsteinOsher2009SplitBregman} introduces a splitting
variable $\bd$ to decouple the $\ell_1$ and $\ell_2$ terms. The
problem~\eqref{eq:l1b} can be written in constrained form as
\begin{equation*}
    \min_{\bx,\bd}
    \|\bA\bx-\bb\|_2^2+\mu\|\bd\|_1
    \quad \textrm{s.t.}\quad
    \bL\bx=\bd.
\end{equation*}
At iteration $k$, the $\bx$-update used in our implementation is
\begin{equation}\label{eq:Breg3}
    \bx^{(k+1)}
    =
    \argmin_{\bx}
    \left\{
    \|\bA\bx-\bb\|_2^2
    +
    \lambda^{(k)}
    \|\bL\bx-(\bd^{(k)}-\bg^{(k)})\|_2^2
    \right\}.
\end{equation}
This has the form~\eqref{eq:generic_sub} with
$\gamma^{(k)}=\lambda^{(k)}$,
$\bD^{(k)}=\bL$, and
$\bc^{(k)}=\bd^{(k)}-\bg^{(k)}$.
The split variable is then updated by
\begin{equation*}
    \bd^{(k+1)}
    =
    \mathrm{shrink}\!\left(
    \bL\bx^{(k+1)}+\bg^{(k)},\,\tau
    \right),
\end{equation*}
where
$\mathrm{shrink}(x,\tau)
=\mathrm{sign}(x)\max(|x|-\tau,0)$
component-wise. The multiplier is subsequently updated as
$$\bg^{(k+1)}=\bg^{(k)} + \left(\bL\bx^{(k+1)}-\bd^{(k+1)}\right).$$
Following~\cite{sweeney2025parameter}, we keep the shrinkage threshold
$\tau$ fixed across iterations while estimating the regularization parameter $\lambda^{(k)}$ by GCV during the first three outer iterations and fixing it thereafter. In the notation of~\cite{sweeney2025parameter}, our parameter $\lambda^{(k)}$ corresponds to the square of their parameter, since their inner problem is written with $(\lambda^{(k)})^2/2$ multiplying the quadratic regularization term. Thus, fixing $\tau$ is equivalent to keeping the ratio between the $\ell_1$ parameter and the quadratic penalty parameter fixed. In our experiments, we use $\tau=0.005$ as in~\cite{sweeney2025parameter}. The same stopping criteria as in IRLS and MM are used. Algorithm~\ref{alg:SB2} summarizes the SB iteration.
\begin{algorithm}[H]
\caption{The SB Method}
\label{alg:SB2}
\begin{algorithmic}[1]
\Require $\bA, \bb, \bL, \tau, \bd^{(0)} = \bg^{(0)} = {\mathbf{0}}, tol, K_{max}$
\Ensure $\bx$
    \For{$k=0,1,\dots$ until $\textrm{RC}(\bx^{(k+1)}) < tol$ or $k = K_{max}$}
    \State Estimate $\lambda^{(k)}$ via GCV (using pre-computed GSVD) for $k \le 3$; set $\lambda^{(k)} = \lambda^{(3)}$ for $k > 3$
    \State Solve~\eqref{eq:Breg3} using LSQR to obtain $\bx^{(k+1)}$
    \State $\bd^{(k+1)} = \mathrm{shrink}\!\left(\bL\bx^{(k+1)}+\bg^{(k)},\,\tau\right)$ \quad (component-wise)
    \State $\bg^{(k+1)} = \bg^{(k)} + \left(\bL\bx^{(k+1)} - \bd^{(k+1)}\right)$
    \EndFor
\end{algorithmic}
\end{algorithm}

\section{GCV Parameter Selection and LSQR}\label{sec:LSQR}
The regularization parameter $\gamma^{(k)}$ in each subproblem~\eqref{eq:generic_sub} is chosen by minimizing the GCV functional. Following~\cite{sweeney2025parameter}, which extends GCV to problems with a nonzero shift $\bc^{(k)}$, the solution of~\eqref{eq:generic_sub} for a trial parameter $\gamma$, with $\bD^{(k)}$ and $\bc^{(k)}$ held fixed, is
\begin{equation}\label{eq:xgamma}
    \bx_\gamma = \bA_\gamma^{\sharp}\bb + \bD_\gamma^{\sharp}\bc^{(k)},
\end{equation}
where
$$
\bA_\gamma^{\sharp} := \bA_\gamma^{-1}\bA^\top, \qquad
\bD_\gamma^{\sharp} := \gamma\,\bA_\gamma^{-1}(\bD^{(k)})^\top, \qquad
\bA_\gamma := \bA^\top\bA + \gamma(\bD^{(k)})^\top\bD^{(k)}.
$$
The GCV functional is then
\begin{equation}\label{eq:GCV}
    \mathrm{GCV}(\gamma) = \frac{\|\bA\bx_\gamma - \bb\|_2^2}{\bigl(m^2 - \mathrm{tr}(\bA\bA_\gamma^{\sharp})\bigr)^2},
\end{equation}
where $\bA\bA_\gamma^{\sharp}$ is the resolution matrix associated with $\bb$; note that the trace depends only on $\bA_\gamma^\sharp$, since GCV's leave-one-out argument treats $\bb$ as the noisy data and $\bc^{(k)}$ as fixed~\cite{sweeney2025parameter}. The parameter $\gamma^{(k)} := \arg\min_\gamma \mathrm{GCV}(\gamma)$ is then used in~\eqref{eq:generic_sub}. In our implementation, we evaluate~\eqref{eq:GCV} efficiently using the GSVD of the pair $\{\bA, \bD^{(k)}\}$, avoiding the need to form $\bA_\gamma^\sharp$ explicitly. At each outer iteration where GCV is needed, the stored generalized singular values are used to evaluate the GCV functional over a grid of 50 log-spaced candidate parameters to identify a coarse minimizer; \texttt{fminbnd} is then used to refine this estimate within a neighborhood of the grid minimizer, since the true minimizer generally falls between grid points. This makes parameter selection computationally negligible compared to the LSQR solve.

The cost of this approach depends on how often the GSVD itself must be recomputed, which differs by method. For MM and SB, the operator pair $\{\bA, \bD^{(k)}\}$ is the same for all $k$ since $\bD^{(k)}=\bL$, so the GSVD is computed once, before the outer iteration begins, and the stored values are reused in the first three outer iterations, after which the regularization parameter is held fixed; the parameter values stabilize quickly, so this introduces negligible bias. For IRLS, the weight matrix $\mathrm{diag}(\bw^{(k)})$ changes at every step, so the parameter must instead be re-estimated from the current weighted pair $\{\bA, \mathrm{diag}(\bw^{(k)})\bL\}$: for two-level IRLS, the GSVD is recomputed and the GCV re-evaluated at every outer iteration; for three-level IRLS, a two-pass strategy is used instead, in which an initial GCV estimate is obtained from the unweighted pair, three warm-up iterations are performed to update the weights, and a second GCV estimate is obtained from the resulting weighted pair, with this fixed parameter then used throughout the fine-level solve.

Once $\gamma^{(k)}$ has been determined, the subproblem~\eqref{eq:generic_sub} is solved using LSQR. To this end, we rewrite~\eqref{eq:generic_sub} as the unweighted augmented least-squares problem
\begin{equation}\label{eq:augmented}
\min_{\bx}\|\hat{\bA}\bx - \hat{\bb}\|_2^2,
\qquad
\hat{\bA} = \begin{bmatrix}\bA \\ \sqrt{\gamma^{(k)}}\,\bD^{(k)}\end{bmatrix},
\qquad
\hat{\bb} = \begin{bmatrix}\bb \\ \sqrt{\gamma^{(k)}}\,\bc^{(k)}\end{bmatrix},
\end{equation}
which specializes in the three solvers as:
\begin{align*}
&\text{IRLS:}\quad
&\hat{\bA} & = \begin{bmatrix}\bA \\ \sqrt{\mu^{(k)}}\,\mathrm{diag}(\bw^{(k)})\bL\end{bmatrix},\quad
&\hat{\bb} &= \begin{bmatrix}\bb \\ \mathbf{0}\end{bmatrix};\\
&\text{MM:}\quad
&\hat{\bA} & = \begin{bmatrix}\bA \\ \sqrt{\mu^{(k)}}\bL\end{bmatrix},\quad
&\hat{\bb}& = \begin{bmatrix}\bb \\ \sqrt{\mu^{(k)}}\bw^{(k)}\end{bmatrix};\\
& \text{SB:} \quad
&\hat{\bA} &=  \begin{bmatrix}\bA \\ \sqrt{\lambda^{(k)}}\bL\end{bmatrix},\quad
&\hat{\bb}& = \begin{bmatrix}\bb \\ \sqrt{\lambda^{(k)}}(\bd^{(k)}-\bg^{(k)})\end{bmatrix}.
\end{align*}

We use LSQR~\cite{paige1982lsqr} to solve~\eqref{eq:augmented}, which applies the Golub-Kahan bidiagonalization~\cite{golub1965calculating} to build orthonormal bases $\bU_{j+1}$ and $\bV_j$ that satisfy $\hat{\bA}\bV_j = \bU_{j+1}\bB_j$ for a lower bidiagonal matrix $\bB_j \in \bbR^{(j+1)\times j}$ using only matrix-vector products with $\hat{\bA}$ and $\hat{\bA}^\top$. The iterate at step $j$ is $\bx_j = \bV_j\by_j$, where
$\by_j = \argmin_{\by} \|\bB_j\by - \beta_1 \bbe_1\|_2^2$. Because $\hat{\bA}$ is never formed explicitly, LSQR requires only products with $\bA$, $\bA^\top$, $\bL$, and $\bL^\top$, and it scales well to large problems. Other iterative solvers could be used in place of LSQR~\cite{chung2024computational, KilmerHansenEspanol, lampe2012large}, but are not considered here.

\section{Wavelet-Based Multilevel Methods}
\label{sec:multilevel}

The multilevel approach defines a hierarchy of compressed operators and right-hand sides of decreasing size~\cite{MGM},
$$
\bA^{(i)} \in \bbR^{(m/2^i)^2 \times (m/2^i)^2}, \qquad \bb^{(i)} \in \bbR^{(m/2^i)^2},
\qquad 0 \le i \le n,
$$
where $i=0$ is the finest level and $i=n$ is the coarsest. At each level the $\ell_1$ problem~\eqref{eq:l1b} is posed on the compressed operator $\bA^{(i)}$ and right-hand side $\bb^{(i)}$. Restriction and prolongation operators $\bR^{(i)}$ and $\bP^{(i)}$
transfer information between levels according to
$$\bA^{(i+1)}=\bR^{(i)}\bA^{(i)}\bP^{(i)},\qquad\bb^{(i+1)}=\bR^{(i)}\bb^{(i)}.$$
We construct $\bR^{(i)}$ and $\bP^{(i)}$ from the low-pass approximation component of the Discrete Wavelet Transform (DWT) \cite{austin2022image, espanol2014wavelet}. Setting $\bR^{(i)} = \bW^\top$ and $\bP^{(i)} = \bW$, where $\bW = \bW_1 \otimes \bW_1$ is the 2D low-pass operator and $\bW_1^\top$ is the 1D low-pass analysis matrix:
\begin{equation*}
\bW_1^\top=\left[ \begin{smallmatrix}
h_0     &h_1    &\dots  &\dots  &h_{\kappa-2}&h_{\kappa-1}&       &       &       &     &       & \\
        &       &h_0    &h_1    &\dots  &\dots  &h_{\kappa-2}&h_{\kappa-1}&       &     &       & \\
        &       &       &       &\dots  &\dots  &\dots  &\dots  &\dots  &\dots&       &\\
        &       &       &       &       &       &h_0    & h_1   & \dots &\dots&h_{\kappa-2}&h_{\kappa-1}\\
h_{\kappa-2} &h_{\kappa-1}&       &       &       &       &       &       &     h_0    &h_1      &\dots & \dots\\
\dots   &\dots  &\dots  &\dots  &       &       &       &       &       &     &  \dots&\dots\\
\dots   &\dots  &h_{\kappa-2}&h_{\kappa-1}&       &       &       &       &       &     &h_0    &h_1
\end{smallmatrix}\right],
\end{equation*}
where $h_0,\ldots,h_{\kappa-1}$ are the low-pass filter coefficients. Each row applies the filter with a stride of two, halving the signal length at each level. For the Haar wavelet, the coefficients are
$$(h_0,h_1)=\left(\tfrac{\sqrt{2}}{2},\tfrac{\sqrt{2}}{2}\right)$$
and for the Daubechies four-tap (Db4) filter, they are
$$ (h_0,h_1,h_2,h_3)=\tfrac{1}{4\sqrt{2}} (1+\sqrt{3},\,3+\sqrt{3},\,3-\sqrt{3},\,1-\sqrt{3}).$$
The Haar wavelet is a two-tap filter whose approximation coefficients are local averages of adjacent pixel pairs; the resulting coarse image is piecewise-constant, preserving edge locations sharply at the cost of a blocky representation of smooth regions. The Db4 filter has two vanishing moments, so it reproduces linear polynomials exactly and produces a smoother, more faithful coarse approximation of gradually varying regions. The trade-off is that Db4's wider support spreads edge transitions across more coefficients, so edges are represented less sharply at the coarse level. This structural difference influences which information can be transferred effectively between levels, as we discuss next.

\subsection{Information Transfer Strategies}\label{sec:init}

We consider two multilevel information transfer strategies together with a single-level baseline.

\noindent \textbf{Approach~I} (single-level baseline): the solver is applied at the finest level only, initialized with identity weights (IRLS) or zero auxiliary variables (SB and MM).

\noindent \textbf{Approach~X} (solution transfer): only the coarse solution $\bx^{(i+1)}$ is interpolated to initialize the fine-level solve. The interpolation is performed as
$$\bX^{(i)} = \bW_1 \bX^{(i+1)} \bW_1^\top,$$
where $\bX^{(i+1)}$ is the coarse solution reshaped as a $(m/2^{i+1}) \times (m/2^{i+1})$ matrix, followed by column-wise vectorization to obtain $\bx^{(i)}$. This interpolated solution serves as the starting point for the fine-level outer iteration. For IRLS and MM, the initial weights and surrogate variables are then generated from $\bx^{(i)}$ via their respective update formulas, so the auxiliary variables are implicitly initialized from the interpolated solution rather than from zero. For SB, the initial split variable and multiplier are similarly derived from one shrinkage step applied to $\bx^{(i)}$.

\noindent \textbf{Approach~D} (auxiliary transfer): coarse-level auxiliary quantities associated with the regularization operator are transferred to the fine level. These solver-specific quantities provide an informed initialization for the fine-level iteration. Because $\bL$ separates horizontal and vertical finite differences, the transferred auxiliary quantities are not interpolated exactly as in Approach~X; the blockwise construction is described in Section~\ref{sec:computational}.

Combining the two wavelet families with the two information transfer strategies yields four multilevel variants: HWT-D, HWT-X, DWT-D, and DWT-X. Table~\ref {tab:method_summary} summarizes all combinations, and Algorithm~\ref{alg:MGM} summarizes these multilevel approaches.

\begin{table}[htbp]
\caption{Summary of method combinations tested. Approach I is a single-level
baseline (no wavelet transform); Approaches X and D are multilevel and differ
in what is transferred from the coarse to the fine level. All combinations
are tested with IRLS, MM, and SB.}
\label{tab:method_summary}
\centering
\begin{tabular}{lll}
\toprule
Approach & Wavelet & Transfer \\
\midrule
I & --- & None (single level) \\
\midrule
X & HWT & Coarse solution only \\
X & DWT & Coarse solution only \\
\midrule
D & HWT & Auxiliary iterates \\
D & DWT & Auxiliary iterates \\
\bottomrule
\end{tabular}
\end{table}

\begin{algorithm}[ht]
\caption{Multilevel V-Cycle (MGM)}
\label{alg:MGM}
\begin{algorithmic}[1]
\Require Current level $i$, system $\bA^{(i)}\bx \approx \bb^{(i)}$, and $s \in \{\text{Approach X}, \text{Approach D}\}$
\Ensure $\bx^{(i)}$; and $\bq^{(i)}$ if $s=\text{Approach D}$
\If{$i = n$} \Comment{coarsest level}
\State Solve $\min_{\bx}\|\bA^{(i)}\bx-\bb^{(i)}\|_2^2+\mu\|\bL\bx\|_1$ using Alg.~\ref{alg:IRLS}, \ref{alg:MML2}, or~\ref{alg:SB2}
\State Set $\bx^{(i)}$ from the coarse-level solve
\If{$s=\text{Approach D}$}
\State Set the solver-dependent auxiliary quantity $\bq^{(i)}$ from the completed level-$i$ solve
\EndIf
\Else
\State Form the coarse system:
$\bA^{(i+1)}=\bW^\top \bA^{(i)}\bW,
\qquad
\bb^{(i+1)}=\bW^\top \bb^{(i)}$
\If{$s=\text{Approach X}$}
\State $\bx^{(i+1)} \gets \text{MGM}\bigl(i+1, \bA^{(i+1)},\bb^{(i+1)},s\bigr)$
\State Interpolate the coarse solution: $\bx^{(i)}=\bW\bx^{(i+1)}$
\State Initialize the level-$i$ solver with $\bx^{(i)}$
\Else \Comment{Approach D}
\State $\bigl[\bx^{(i+1)},\bq^{(i+1)}\bigr] \gets \text{MGM}\bigl(i+1, \bA^{(i+1)},\bb^{(i+1)},s\bigr)$
\State Transfer the auxiliary quantity: $\bq^{(i)}=\mathcal{T}\bigl(\bq^{(i+1)}\bigr)$
\State Initialize the level-$i$ solver with the transferred auxiliary quantity $\bq^{(i)}$
\EndIf
\State Solve $\min_x \|\bA^{(i)}\bx - \bb^{(i)}\|_2^2 + \mu\|\bL\bx\|_1$ at level $i$ using the selected initialization, obtaining $\bx^{(i)}$
\If{s = Approach D}
    \State Set the auxiliary quantity $\bq^{(i)}$ from the completed level-$i$ solve \Comment{needed if this call returns to level $i-1$}
\EndIf
\EndIf
\end{algorithmic}
\end{algorithm}

\subsection{Computational Implementation}\label{sec:computational}

The 2D operator $\bW = \bW_1 \otimes \bW_1$ is applied by filtering first along rows and then along columns, so neither the full matrix nor its inverse is formed explicitly. The coarse operator $\bA^{(i+1)} = \bW^\top \bA^{(i)} \bW$ preserves the Toeplitz or block-Toeplitz structure of the blurring matrix under the Haar wavelet transform, enabling fast matrix-vector products at all levels via FFT-based convolution. Under Db4, this structure does not hold exactly, but the separable Kronecker structure of $\bW$ is sufficient for efficient computation since LSQR requires only products with $\bA$ and $\bA^\top$. We emphasize that this Kronecker/Toeplitz structure is not a requirement of the multilevel framework itself: Algorithm~\ref{alg:MGM} and the solvers of Section~\ref{sec:solvers} only require the ability to form matrix-vector products with $\bA^{(i)}$, $(\bA^{(i)})^\top$, $\bL$, and $\bL^\top$ at each level, so the method's formulation is not restricted to Toeplitz or Kronecker-structured operators. Thus, in principle, the framework can also be applied to more general
blurring operators, including spatially varying or non-separable blurs,
provided the required matrix-vector products can be computed efficiently. The structural assumption is used only to accelerate the matrix-vector products in our specific experiments; for a general operator lacking this structure, $\bA^{(i)}$ and $(\bA^{(i)})^\top$ can still be applied by forming the action of $\bW^\top$, $\bA^{(i)}$, and $\bW$ in sequence, at correspondingly higher per-iteration cost.

The regularization operator $\bL$ is defined for 2D images as
$$\bL =\begin{bmatrix} \bI \otimes \bD \\ \bD \otimes \bI \end{bmatrix} \in \mathbb{R}^{2m(m-1)\times m^2},$$
where $\bD$ is the 1D finite difference matrix. Importantly, $\bL$ is \emph{not} transported across levels; instead, the regularization operator is defined directly at the resolution of each level. This ensures that edge-preserving regularization acts on the gradient at each level's resolution, rather than on a compressed approximation of the fine-level gradient.

For Approach~D, the solver-specific auxiliary information computed at the
coarse level is transferred directly to initialize the corresponding
quantities at the next finer level. Specifically, the IRLS weight vector,
the MM surrogate vector, and the SB auxiliary quantity $\bd-\bg$ are
transferred for IRLS, MM, and SB, respectively. For SB, since the transfer
operator is linear, transferring $\bd$ and $\bg$ separately and then forming
their difference is equivalent to transferring $\bd-\bg$ directly. Although
all these quantities have different interpretations, their components are indexed
according to the horizontal and vertical finite-difference blocks associated
with the regularization operator $\bL$. Consequently, each
auxiliary vector can be partitioned as
$$
\bq^{(i+1)}
=
\begin{bmatrix}
\bq_x^{(i+1)}\\
\bq_y^{(i+1)}
\end{bmatrix}.
$$
Here, $\bq_x^{(i+1)}$ and $\bq_y^{(i+1)}$ denote the components associated
with the horizontal and vertical finite-difference blocks, respectively.
They are reshaped as
$(m/2^{i+1}-1)\times(m/2^{i+1})$ and
$(m/2^{i+1})\times(m/2^{i+1}-1)$ matrices.
Because the two blocks have different shapes, each is interpolated
separately using a rectangular restriction $\widetilde{\bW}_1$ of
$\bW_1$ matched to the reduced dimension introduced by the
finite-difference operator $\bD$:
$$
\bQ_x^{(i)}
=
\widetilde{\bW}_1\,\bQ_x^{(i+1)}\,\bW_1^\top,
\qquad
\bQ_y^{(i)}
=
\bW_1\,\bQ_y^{(i+1)}\,\widetilde{\bW}_1^\top.
$$
Here, $\widetilde{\bW}_1$ consists of the first $m/2^i-1$ rows and
$m/2^{i+1}-1$ columns of $\bW_1$. This blockwise transfer is distinct
from the interpolation of the solution $\bx^{(i+1)}$ under Approach~X,
which is defined on a single square grid and is interpolated using
$\bW_1$ alone.

\section{Numerical Results}
\label{sec:num}
In this section, we assess the proposed multilevel framework through a series of image deblurring experiments. In each experiment, an unknown image is reconstructed from a blurred, noisy observation.

\begin{figure}[tbp]
\centering
\begin{subfigure}[t]{0.25\textwidth}
\centering
\includegraphics[width=\linewidth]{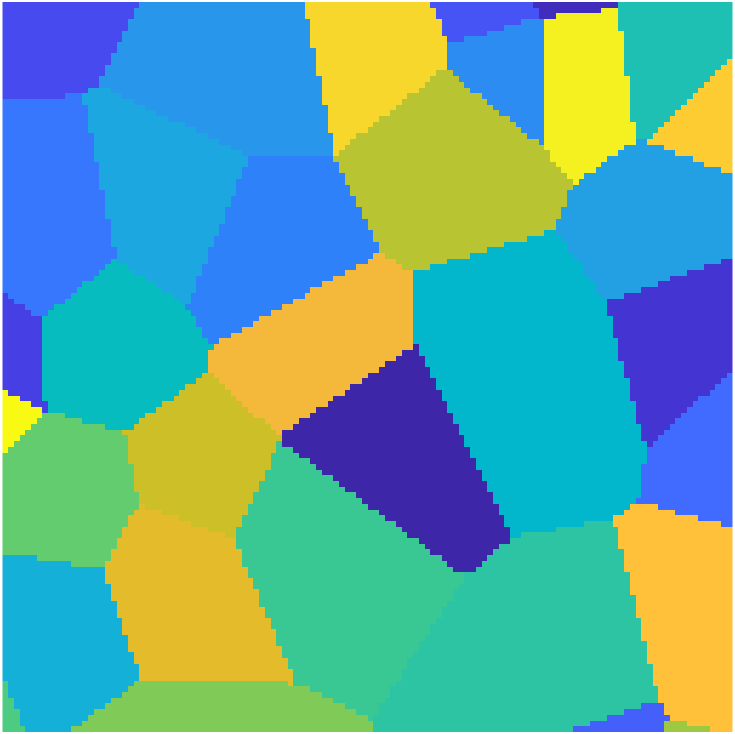}
\caption{True image}
\label{fig:ANa}
\end{subfigure}\hfill
\begin{subfigure}[t]{0.25\textwidth}
\centering
\includegraphics[width=\linewidth]{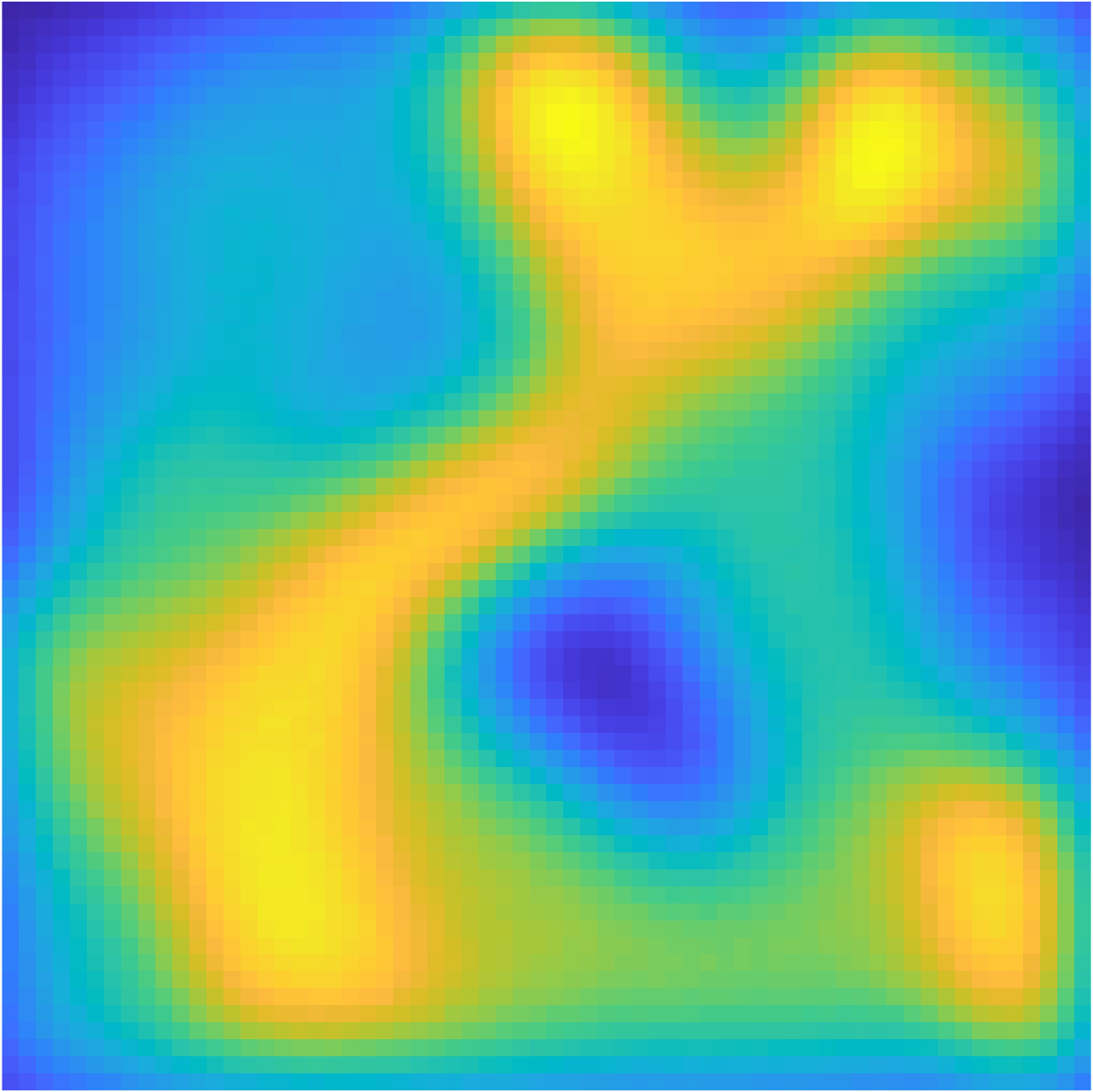}
\caption{$\sigma=3.5$}
\label{fig:ANb3}
\end{subfigure}\hfill
\begin{subfigure}[t]{0.25\textwidth}
\centering
\includegraphics[width=\linewidth]{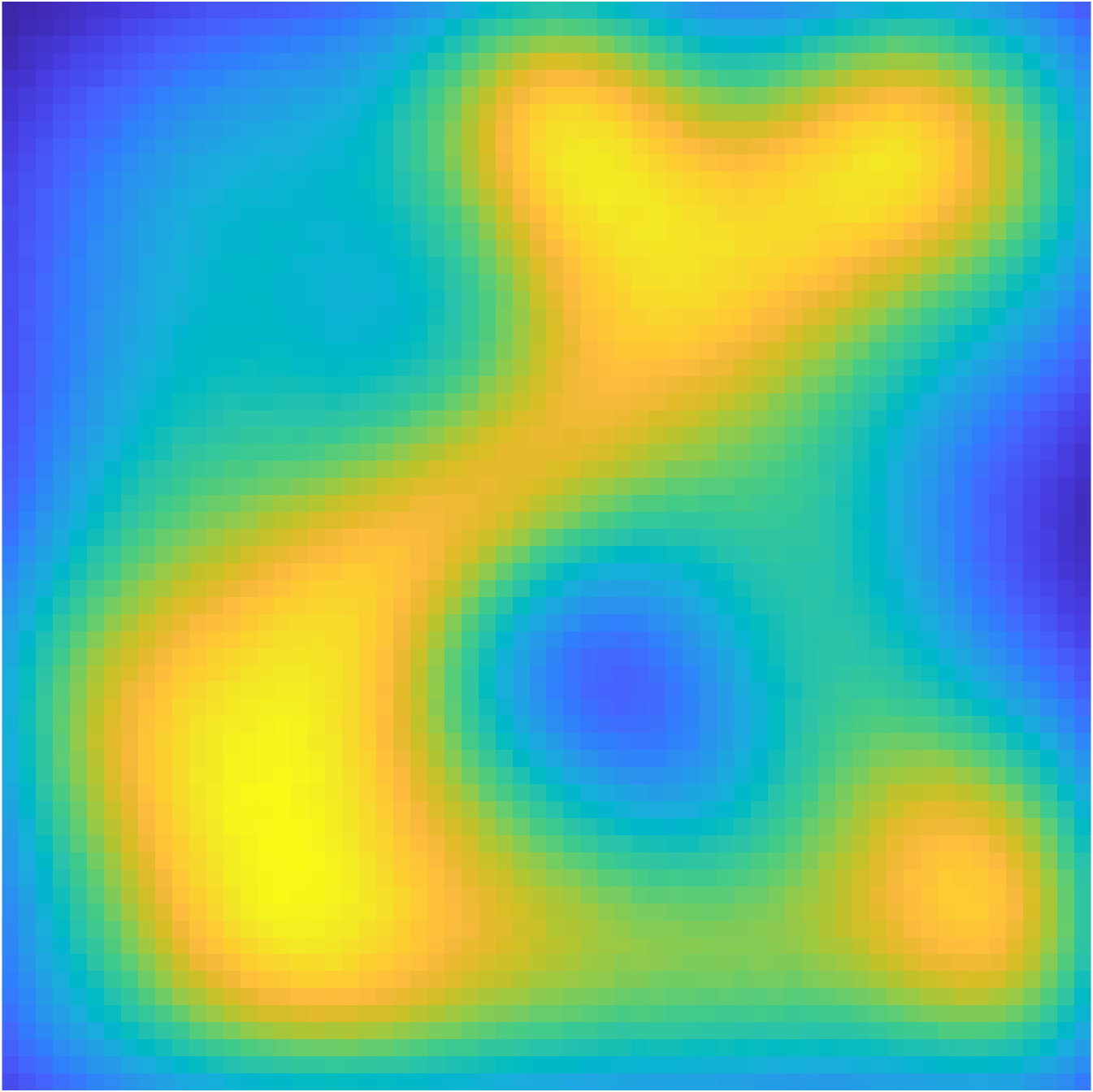}
\caption{$\sigma=4.5$}
\label{fig:ANb4}
\end{subfigure}
\caption{True Voronoi image and corresponding blurred, noisy observations.}
\label{fig:ANimages}
\end{figure}

\begin{figure}[tbp]
\centering
\begin{subfigure}[t]{0.25\textwidth}
\centering
\includegraphics[width=\linewidth]{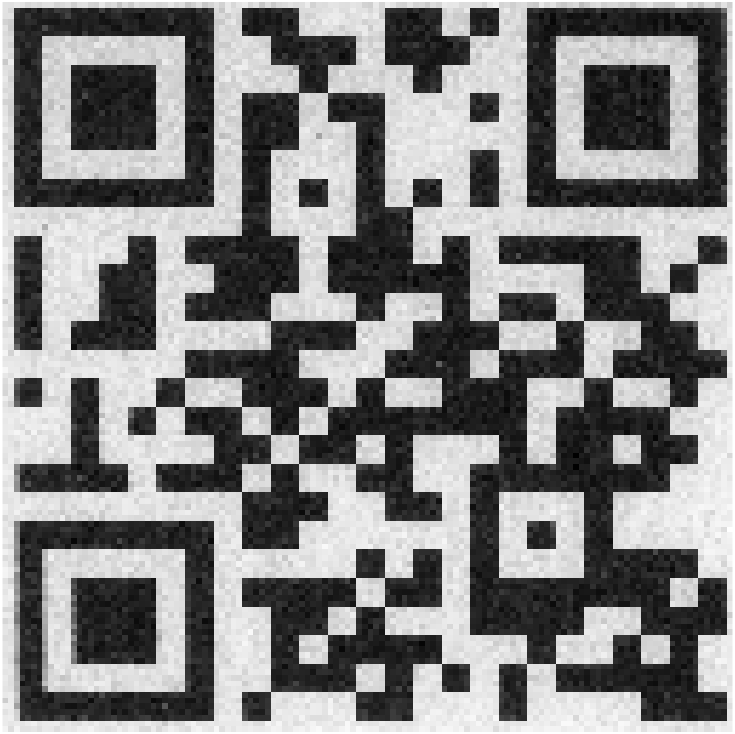}
\caption{True image}
\label{fig:QRa}
\end{subfigure}\hfill
\begin{subfigure}[t]{0.25\textwidth}
\centering
\includegraphics[width=\linewidth]{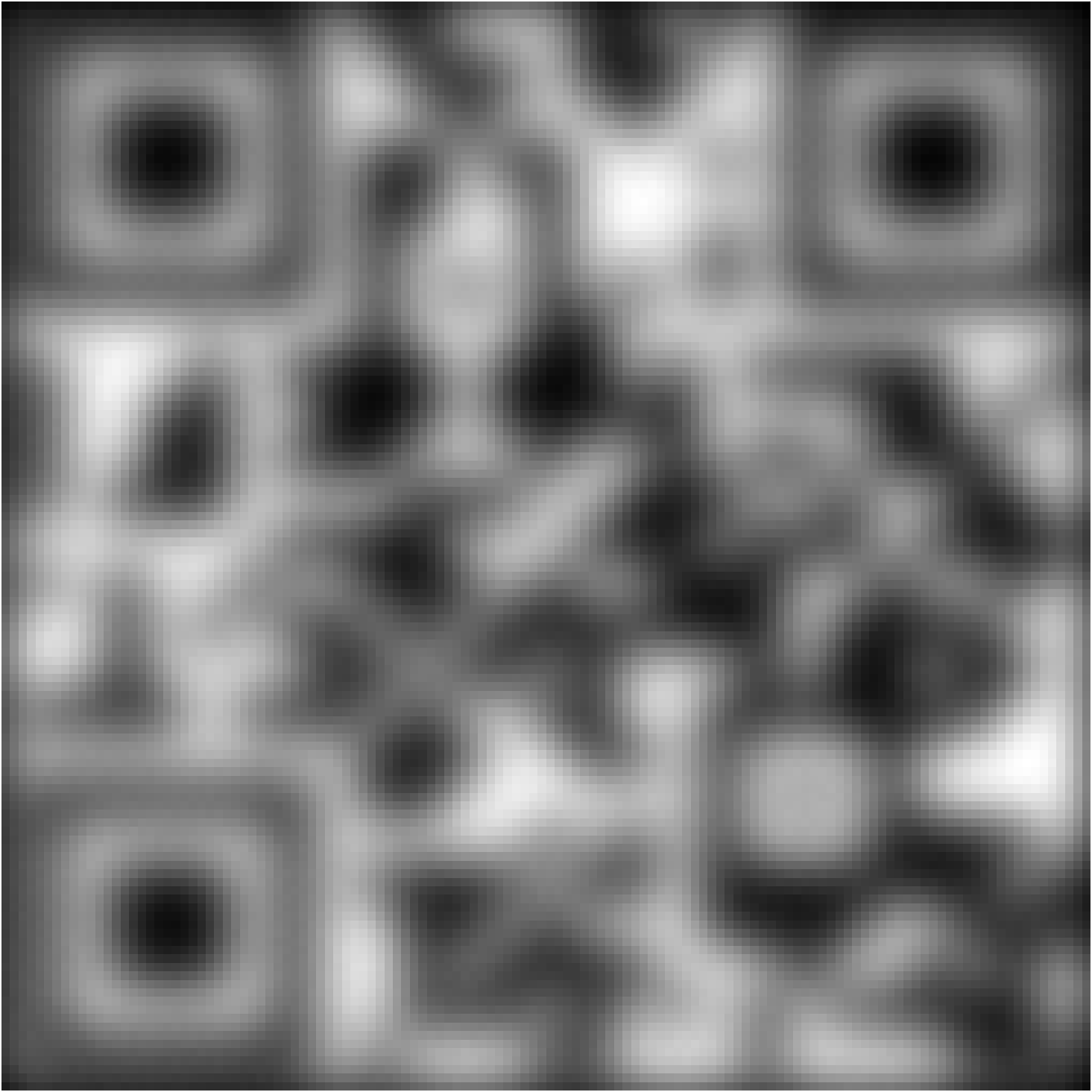}
\caption{$\sigma=3.5$}
\label{fig:QRb3}
\end{subfigure}\hfill
\begin{subfigure}[t]{0.25\textwidth}
\centering
\includegraphics[width=\linewidth]{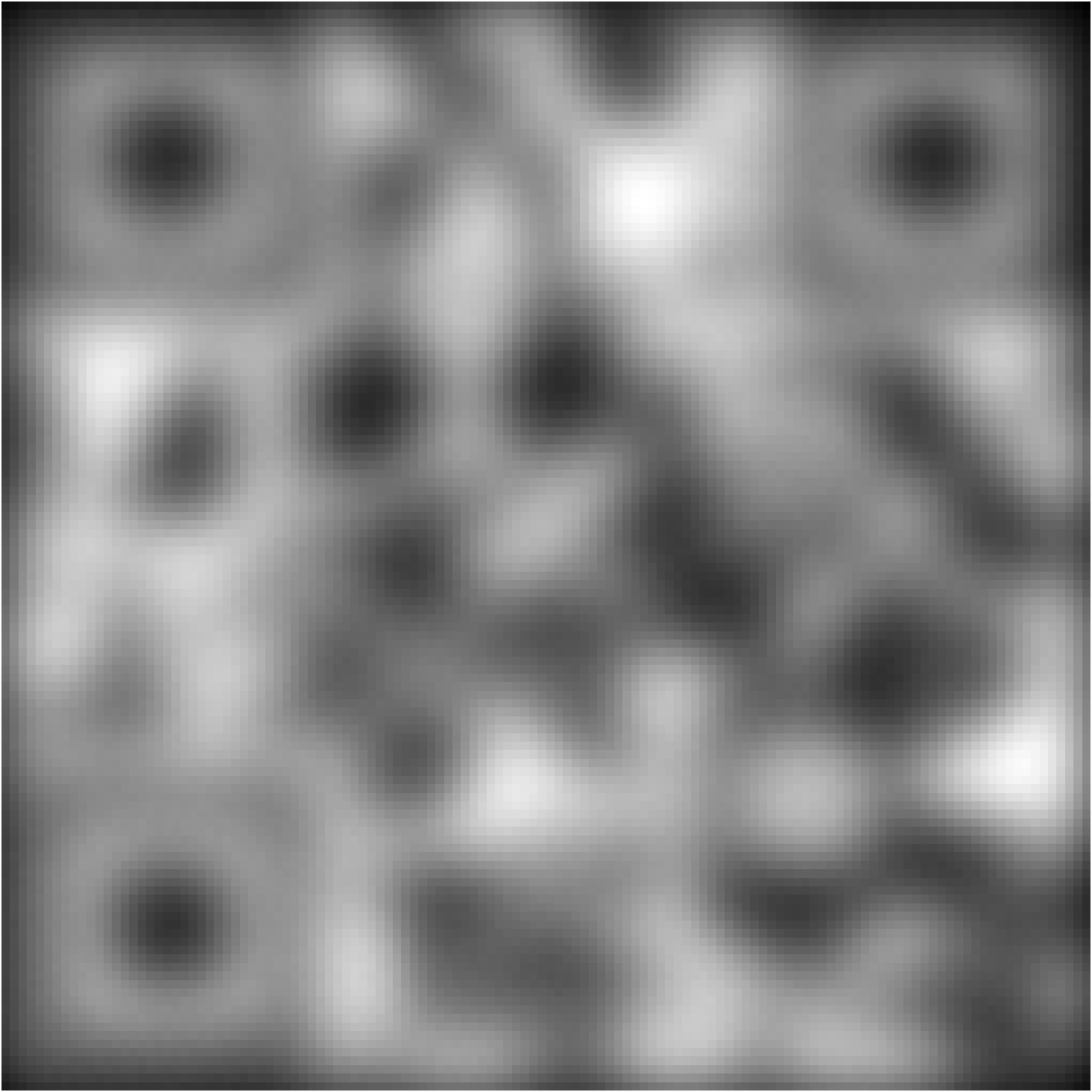}
\caption{$\sigma=4.5$}
\label{fig:QRb4}
\end{subfigure}
\caption{True QR image and corresponding blurred, noisy observations.}
\label{fig:QRimages}
\end{figure}

We use two test images: a Voronoi diagram (Figure~\ref{fig:ANimages}) and a QR code image (Figure~\ref{fig:QRimages}), the latter obtained from \texttt{fips.fi}~\cite{fips}. For the two-level experiments, images of size $64\times 64$ are used, vectorized column-wise into vectors of length $4{,}096$; for the three-level experiments, images of size $128\times 128$ are used, vectorized into vectors of length $16{,}384$. The Voronoi image is tested at both resolutions, while the QR image is tested only at $128\times 128$. Three-level experiments report results only for Approaches~D and~X, since Approach~I is prohibitively expensive at $128\times 128$.

Each true image $\bx^{\mathrm{true}}$ is blurred by an operator $\bA \in \bbR^{m^2 \times m^2}$, formed by the Kronecker product of two symmetric banded Toeplitz matrices $\tilde{\bA} \in \bbR^{m \times m}$, with the first row defined by the Gaussian kernel
$$
z=\left[\exp\left(-\frac{[0:\mathrm{band}-1]^{\cdot 2}}{2\sigma^2}\right),\; \mathbf{0} \right],$$
where $\sigma>0$ controls the blur width. Blur parameters $\sigma = 3.5$ and $\sigma = 4.5$ are tested throughout, producing  $\bb^{\mathrm{true}} = \bA\bx^{\mathrm{true}}$. The two blur levels considered here correspond to different rates of singular value
decay for $\bA$: larger $\sigma$ produces faster decay. As discussed in~\cite{Amdouni2024},
faster singular value decay allows the wavelet-compressed coarse operator
$\bA^{(i+1)} = \bR^{(i)}\bA^{(i)}\bP^{(i)}$ to retain a larger share of the dominant
singular subspace of $\bA^{(i)}$, so the coarse-level problem remains a faithful
reduced model of the fine-level problem. When the blur is mild and the singular
values decay slowly, coarsening instead discards a larger share of significant
singular value information, so the coarse-level solve is a less reliable surrogate
for the fine-level problem and provides a less informative starting point. Zero-mean Gaussian noise is then added, scaled to a $1\%$ relative noise level, $\|\bbe\|_2/\|\bb^{\mathrm{true}}\|_2 = 10^{-2}$, to produce the observed data $\bb = \bb^{\mathrm{true}} + \bbe$.
Stopping tolerances ($tol$) are $10^{-6}$ at coarse levels and $10^{-3}$ at the fine or single level. For IRLS and MM, we define $\epsilon = \varepsilon=0.03$. The value of $\tau$ is fixed at $0.005$. All computations are carried out in MATLAB.

Reconstruction quality is measured primarily by the relative reconstruction error (RRE),
$$
\mathrm{RRE}(\bx) = \frac{\|\bx-\bx^{\mathrm{true}}\|_2}{\|\bx^{\mathrm{true}}\|_2}.
$$

For three-level experiments, we also report SSIM~\cite{Wang2004SSIM}, computed using MATLAB's \texttt{ssim} function, which evaluates
\begin{equation*}
\mathrm{SSIM}(\bx^{\mathrm{true}}, \bx) =
\frac{(2\mu_{\bx^{\mathrm{true}}}\mu_{\bx} + c_1)(2\sigma_{\bx^{\mathrm{true}}\bx} + c_2)}
{(\mu_{\bx^{\mathrm{true}}}^2 + \mu_{\bx}^2 + c_1)(\sigma_{\bx^{\mathrm{true}}}^2 + \sigma_{\bx}^2 + c_2)}
\end{equation*}
over local Gaussian-weighted windows and averages the result over all windows to obtain a single scalar value, where $\mu_{\bx^{\mathrm{true}}}, \mu_{\bx}$, $\sigma_{\bx^{\mathrm{true}}}^2, \sigma_{\bx}^2$, and $\sigma_{\bx^{\mathrm{true}}\bx}$ denote the local means, variances, and covariance within each window, $c_1 = (0.01\, L)^2$, $c_2 = (0.03\, L)^2$, and $L$ is the dynamic range. A higher SSIM (closer to 1) indicates greater structural similarity. SSIM conclusions are consistent with RRE throughout the two-level experiments and the three-level Voronoi results; one exception, involving the three-level QR results, is discussed in Section~\ref{sec:threelevel}.

For the single-level baseline (Approach~I), iteration $k$ in Figures~\ref{fig:RRENirntv}-\ref{fig:RRENirnSB}
refers to the $k$-th outer iteration of the solver (Algorithm~\ref{alg:IRLS},
\ref{alg:MML2}, or~\ref{alg:SB2}) applied directly at the fine level from a default
initialization. For the multilevel approaches (D and X), iteration $k$ refers instead
to the $k$-th outer iteration of the \emph{fine-level} solve only, counted from the
coarse-to-fine transfer described in Section~\ref{sec:init}; the outer iterations
performed at the coarse level to produce this transferred initialization are not
included in the count. Consequently, iteration~1 for D and X already reflects a
warm start informed by the (unshown) coarse-level solve, so the iteration axis in these figures is not a like-for-like measure of total computational work across approaches. However, the wall-clock comparison in Table~\ref{tab:time_all_solvers} is the fair basis for that comparison.

\subsection{Two-Level Results}

\begin{figure}[tbp]
\centering
\begin{subfigure}[t]{0.19\textwidth}
\centering
\includegraphics[width=\linewidth]{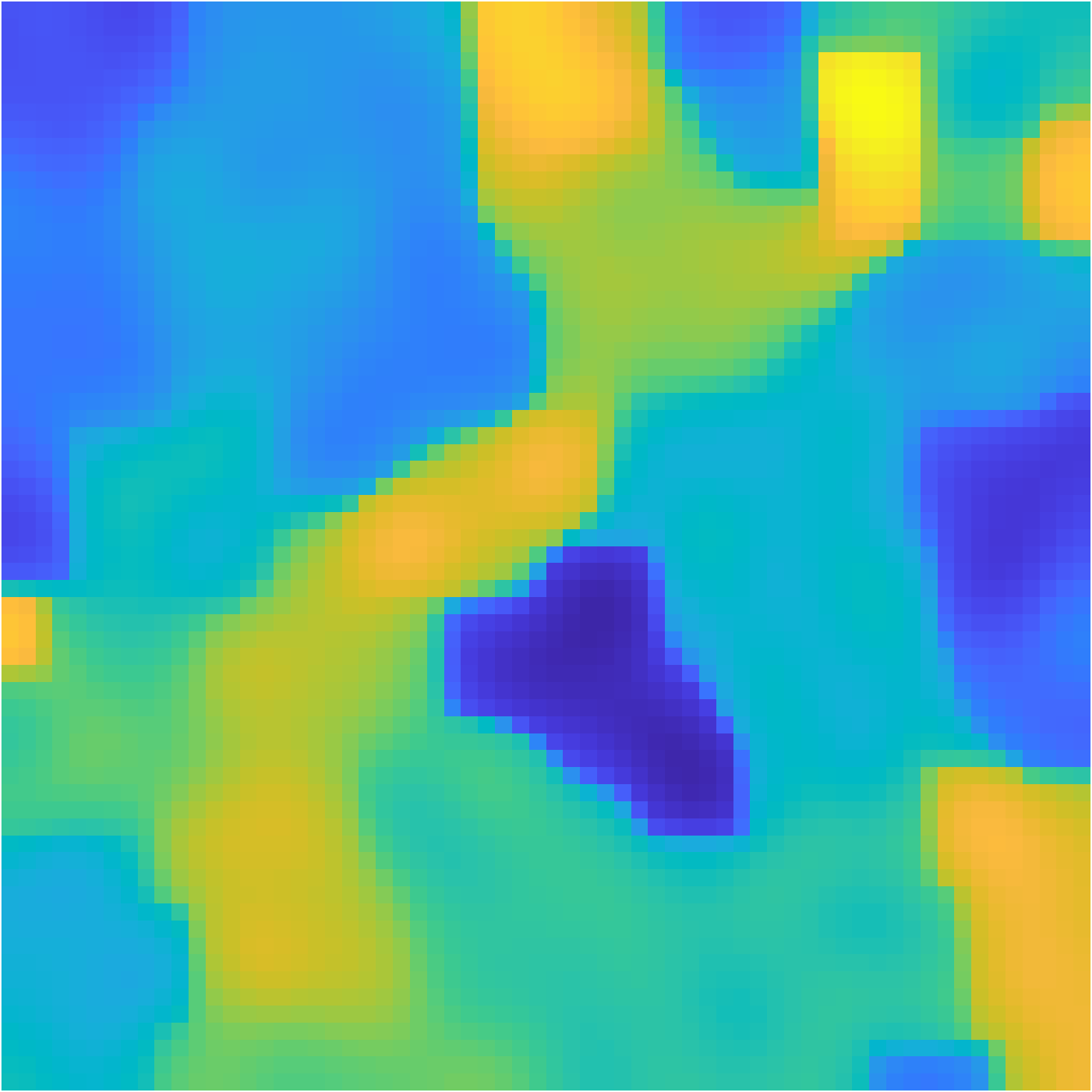}
\caption{I}
\end{subfigure}\hfill
\begin{subfigure}[t]{0.19\textwidth}
\centering
\includegraphics[width=\linewidth]{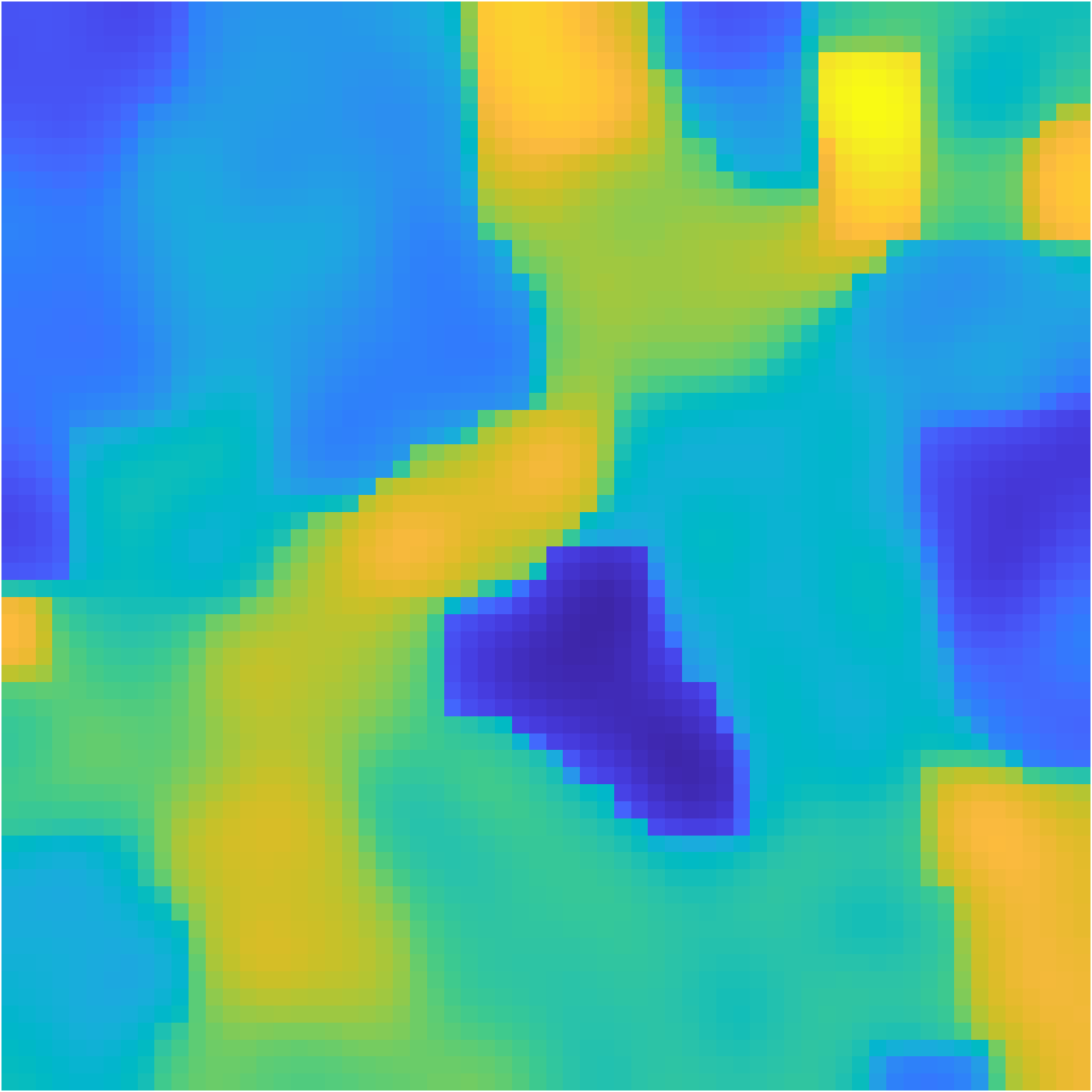}
\caption{HWT-D}
\end{subfigure}\hfill
\begin{subfigure}[t]{0.19\textwidth}
\centering
\includegraphics[width=\linewidth]{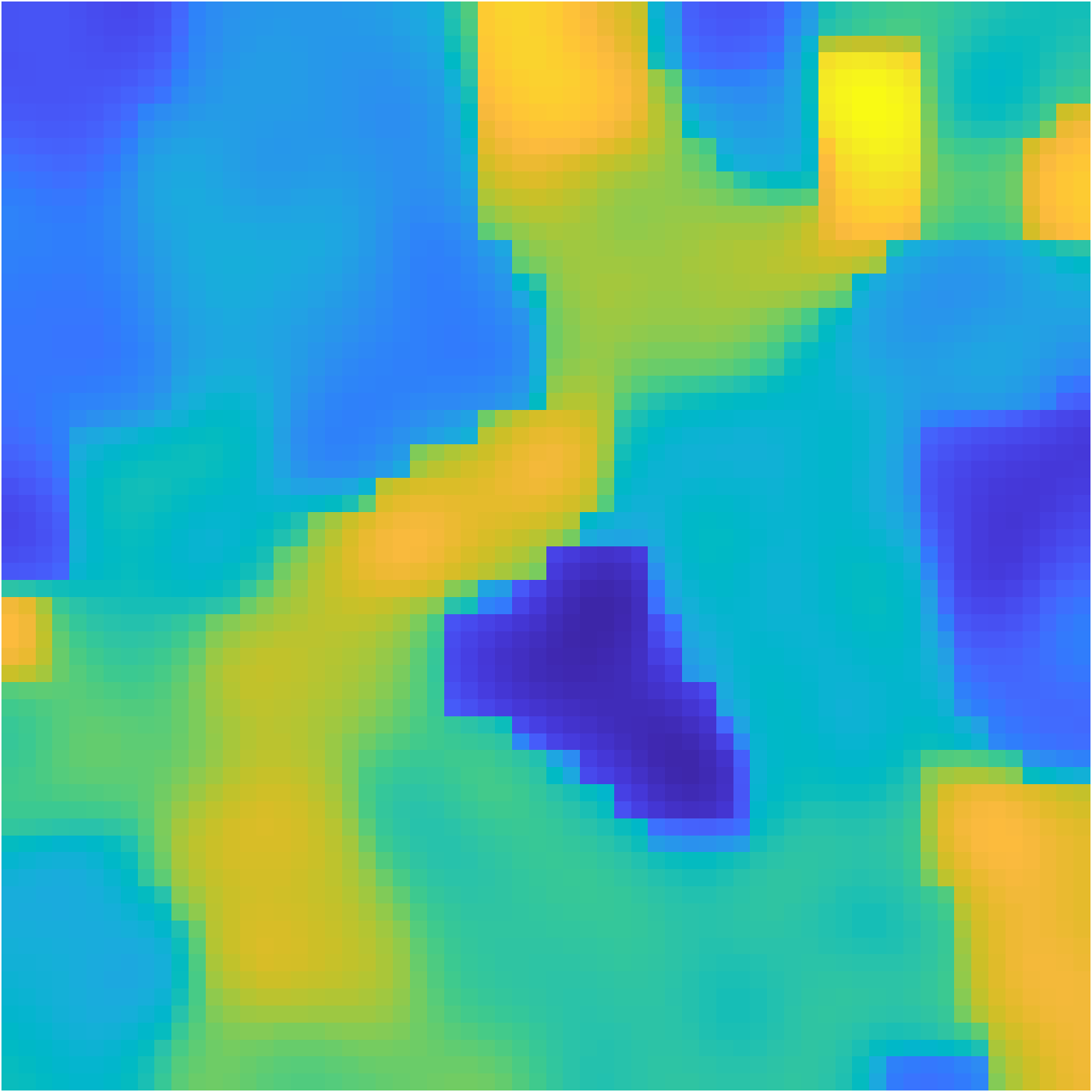}
\caption{HWT-X}
\end{subfigure}\hfill
\begin{subfigure}[t]{0.19\textwidth}
\centering
\includegraphics[width=\linewidth]{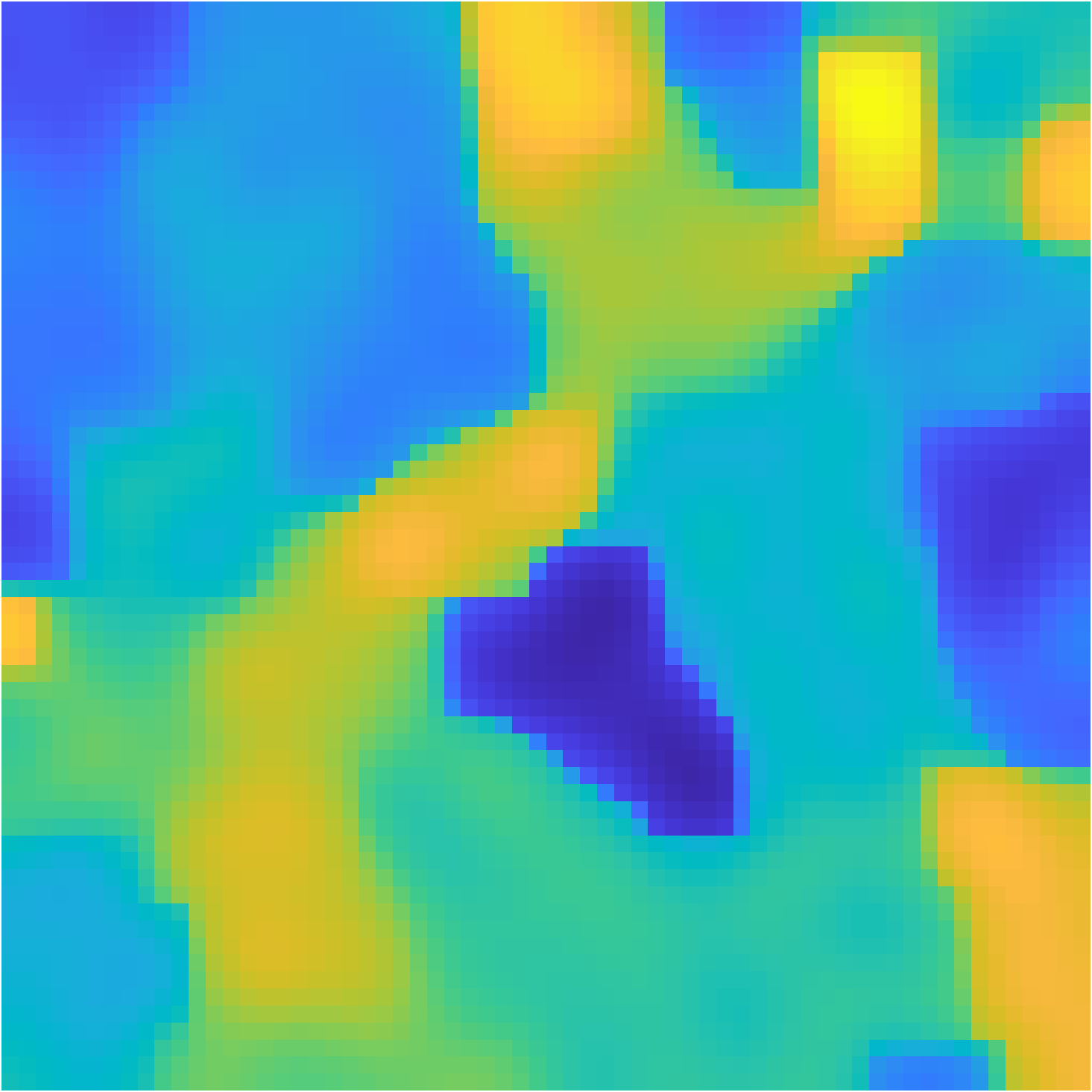}
\caption{DWT-D}
\end{subfigure}\hfill
\begin{subfigure}[t]{0.19\textwidth}
\centering
\includegraphics[width=\linewidth]{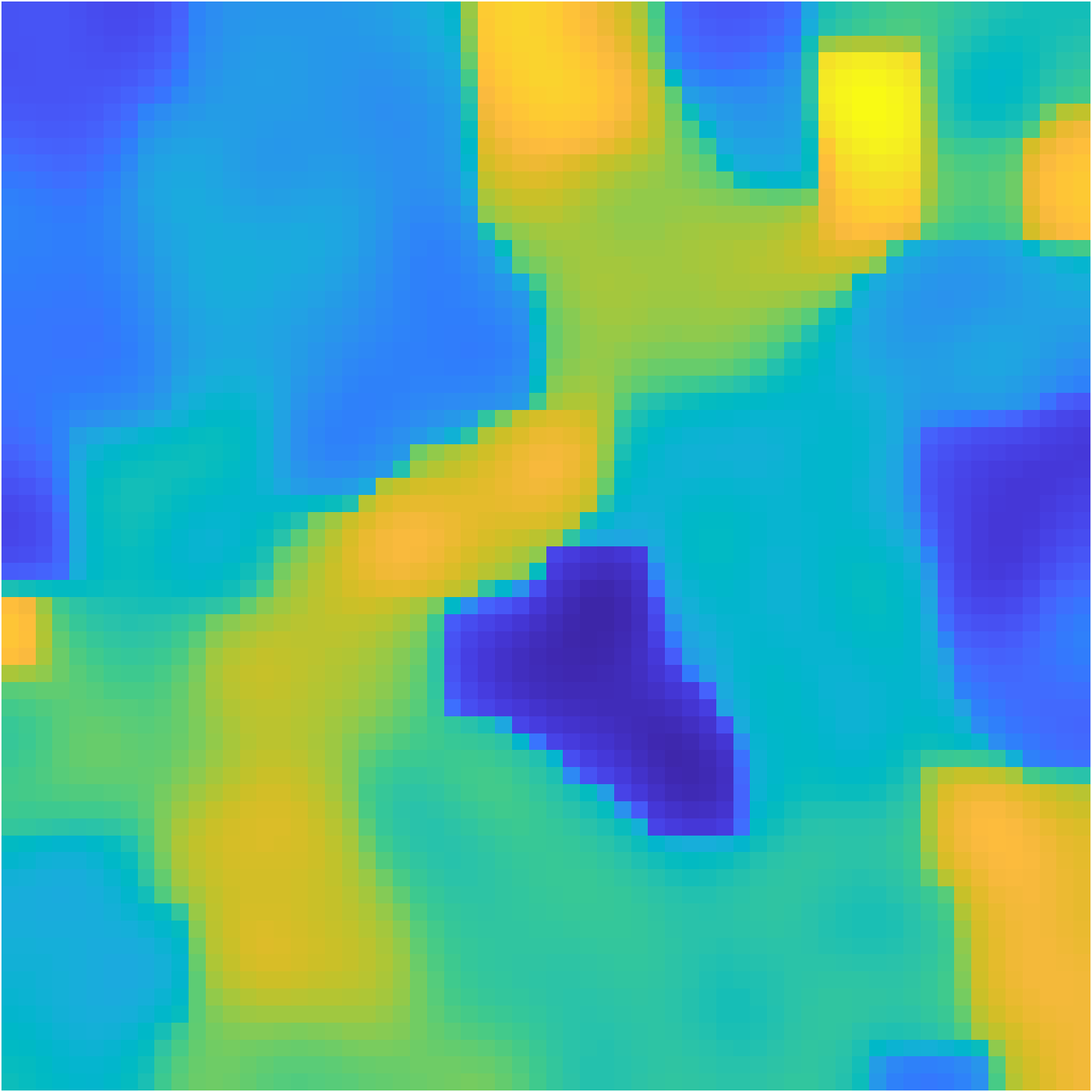}
\caption{DWT-X}
\end{subfigure}
\vspace{0.8em}
\begin{subfigure}[t]{0.19\textwidth}
\centering
\includegraphics[width=\linewidth]{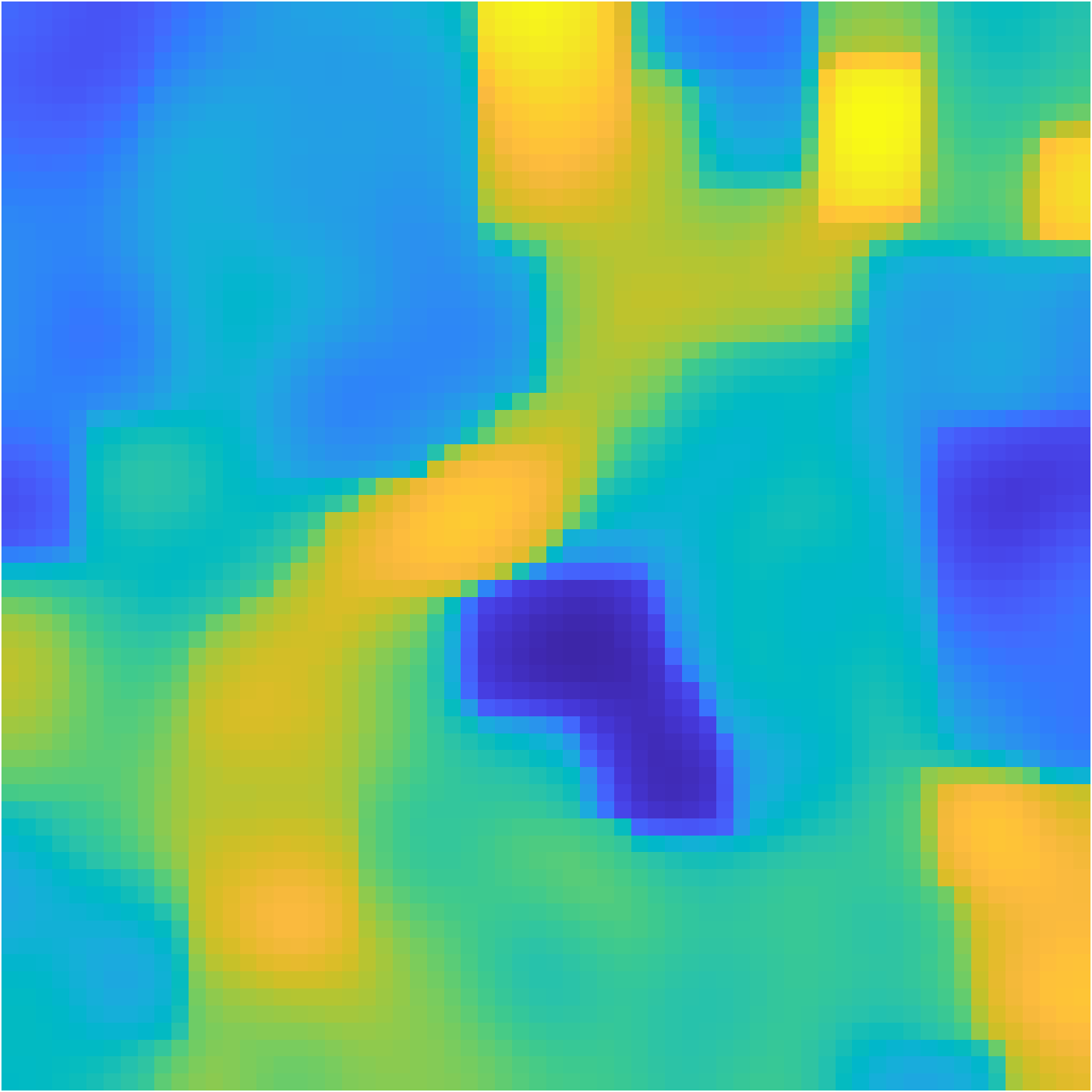}
\caption{I}
\end{subfigure}\hfill
\begin{subfigure}[t]{0.19\textwidth}
\centering
\includegraphics[width=\linewidth]{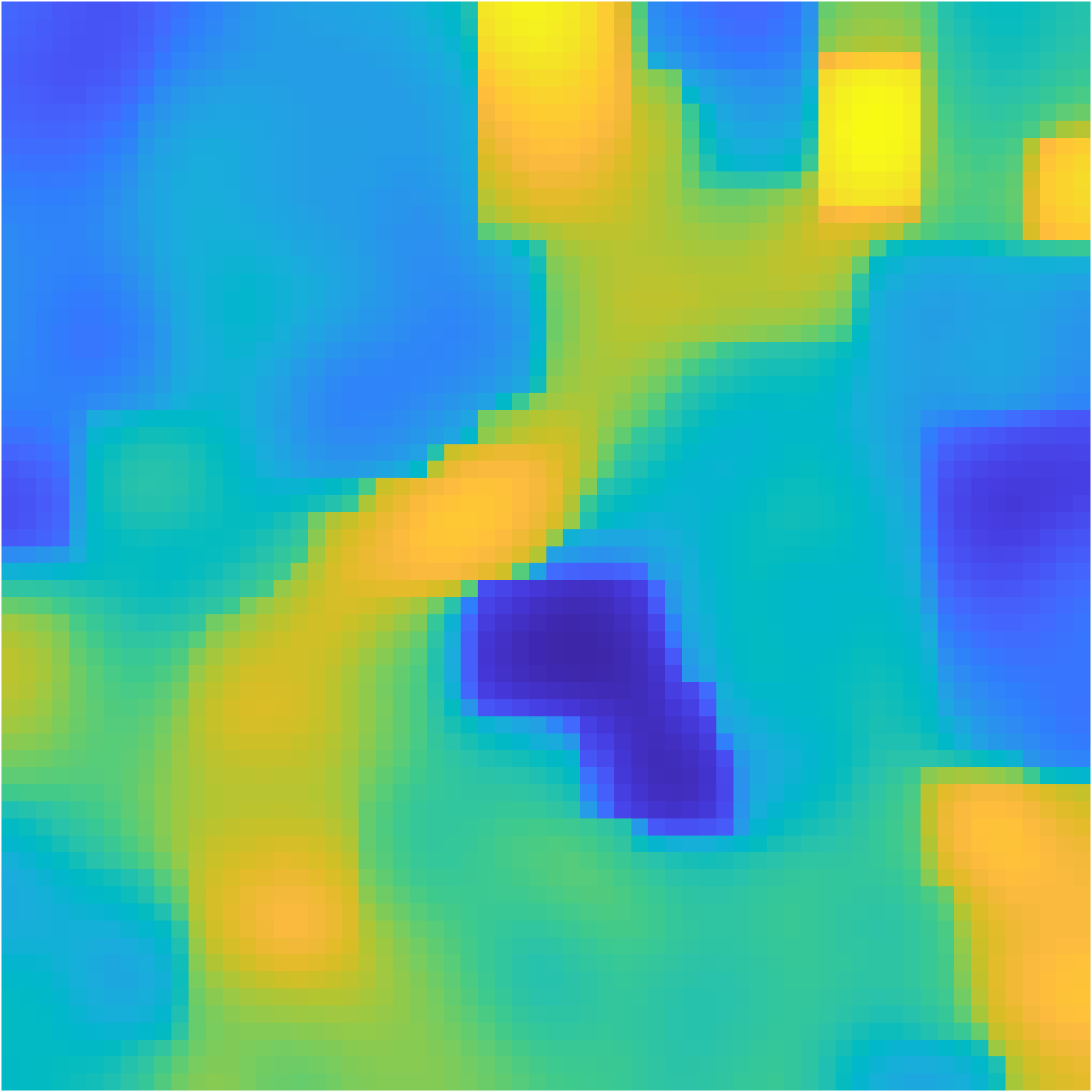}
\caption{HWT-D}
\end{subfigure}\hfill
\begin{subfigure}[t]{0.19\textwidth}
\centering
\includegraphics[width=\linewidth]{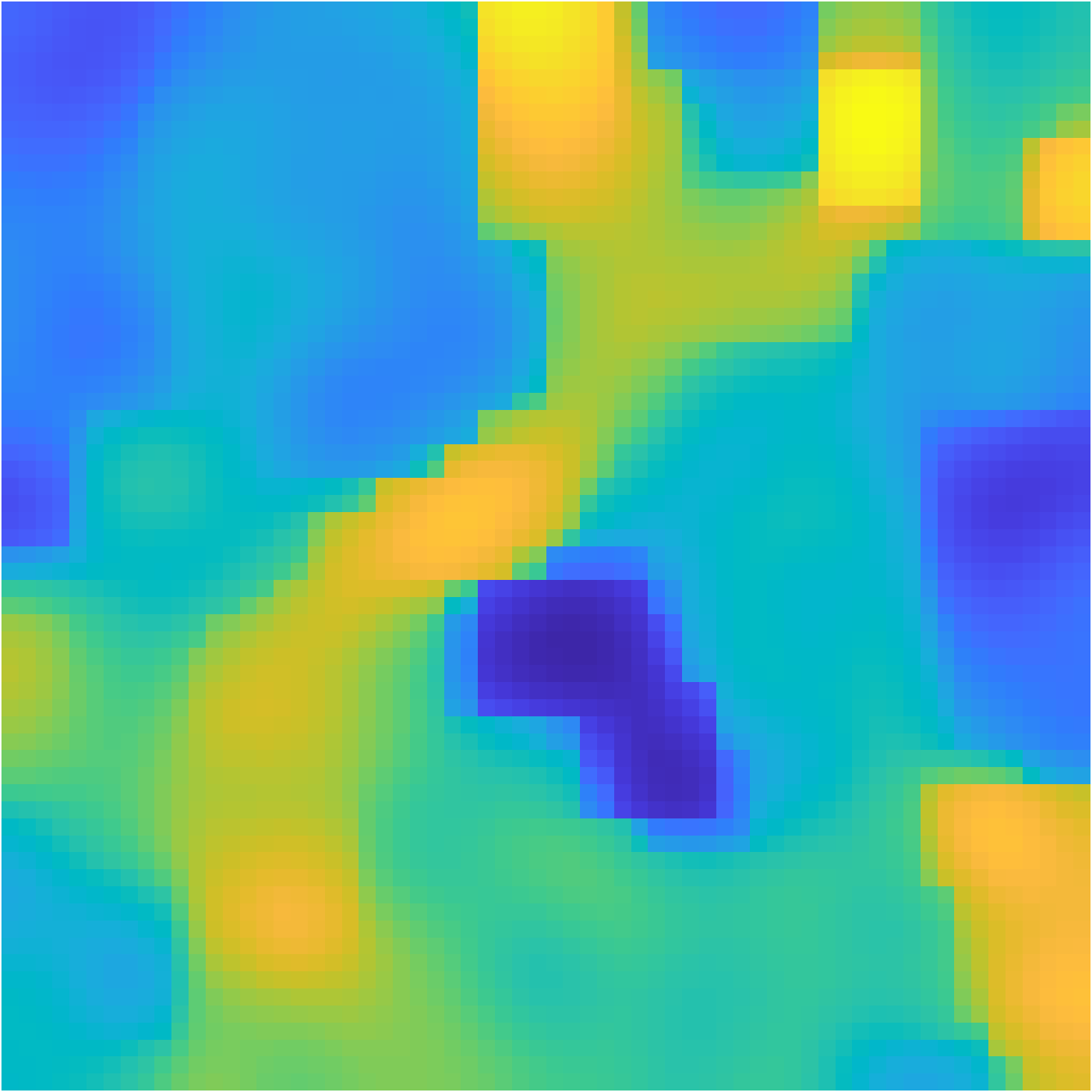}
\caption{HWT-X}
\end{subfigure}\hfill
\begin{subfigure}[t]{0.19\textwidth}
\centering
\includegraphics[width=\linewidth]{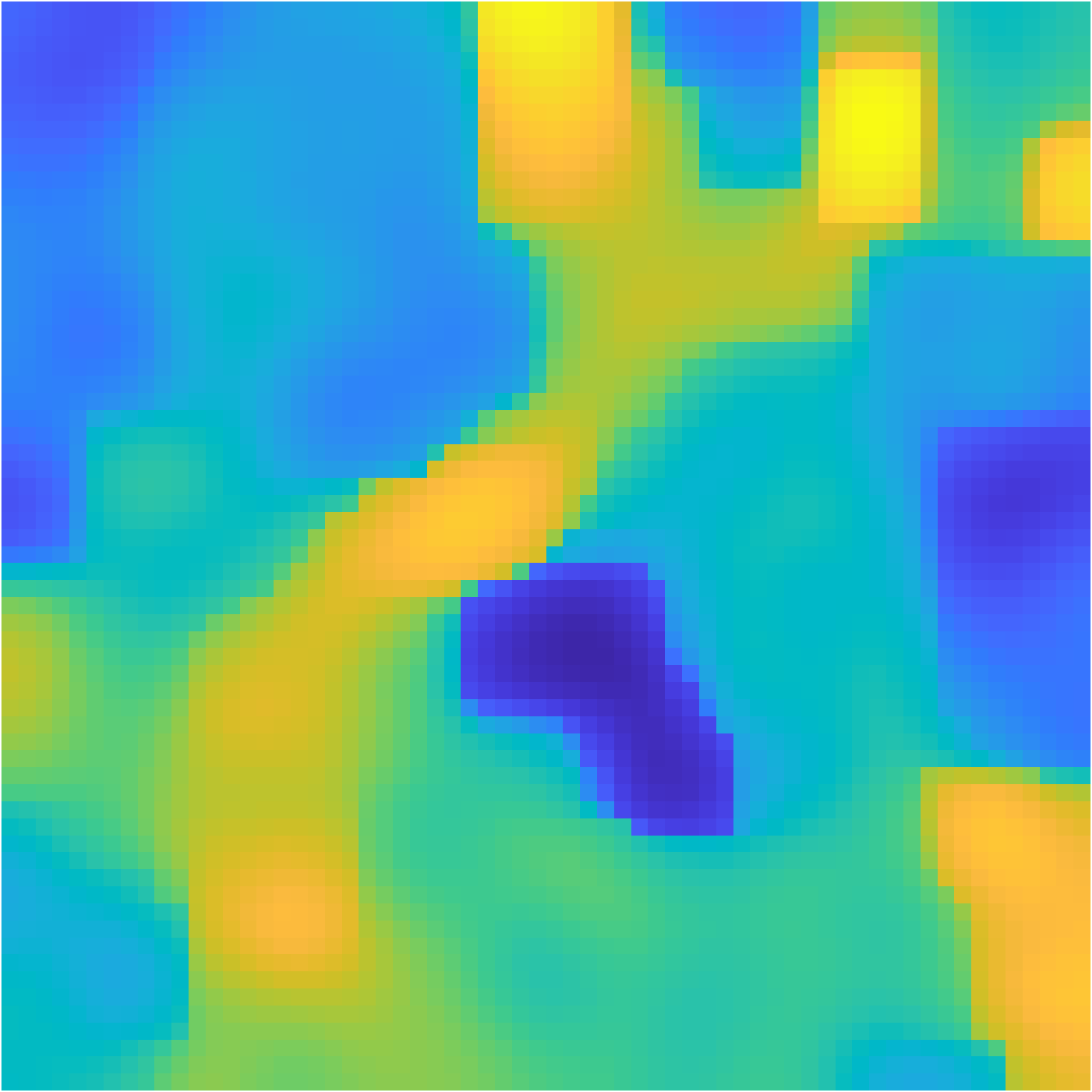}
\caption{DWT-D}
\end{subfigure}\hfill
\begin{subfigure}[t]{0.19\textwidth}
\centering
\includegraphics[width=\linewidth]{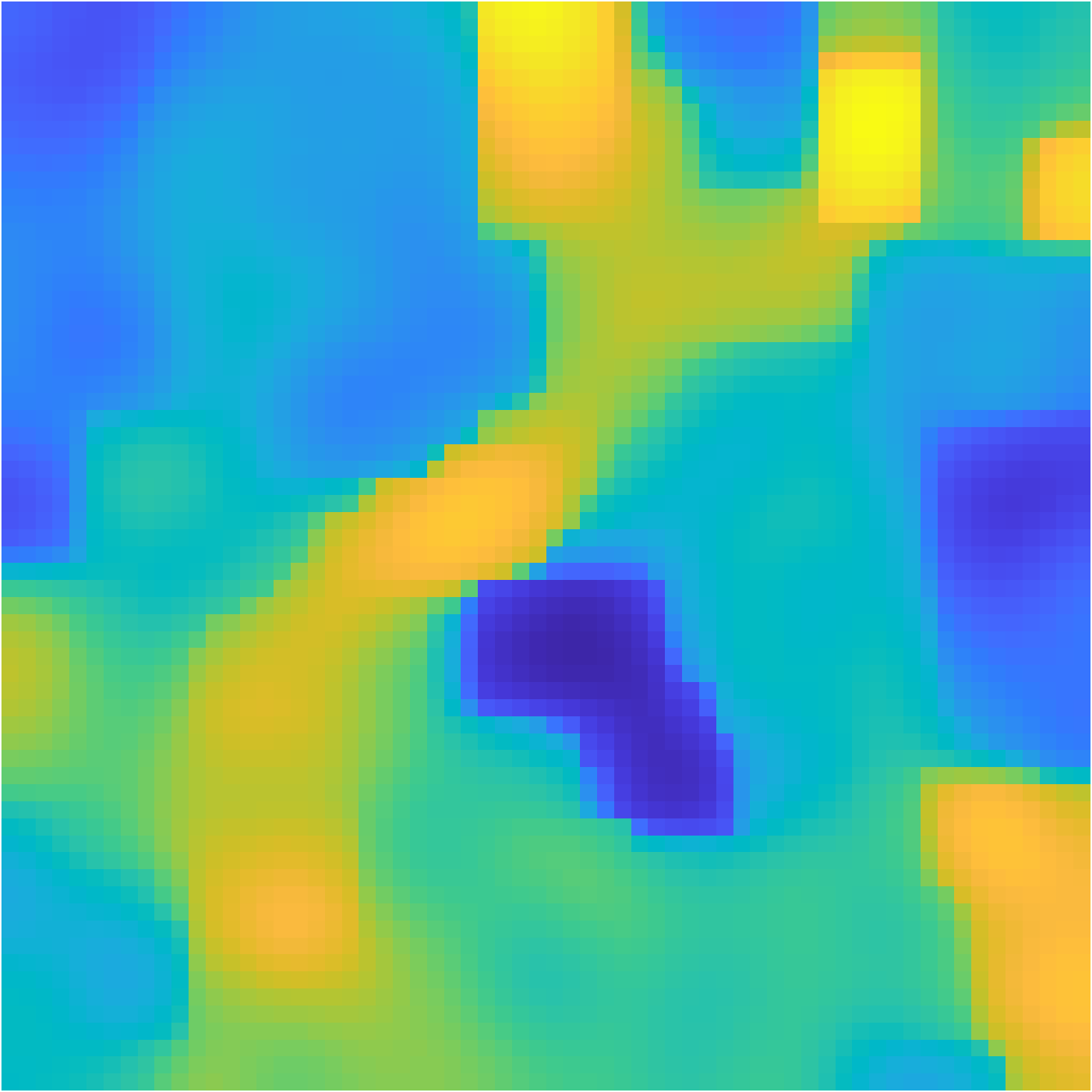}
\caption{DWT-X}
\end{subfigure}
\caption{Deblurred images using IRLS. Top row: $\sigma=3.5$; bottom row: $\sigma=4.5$.}
\label{fig:SolutionTV}
\end{figure}

\begin{figure}[tbp]
\centering
\begin{subfigure}[t]{0.19\textwidth}
\centering
\includegraphics[width=\linewidth]{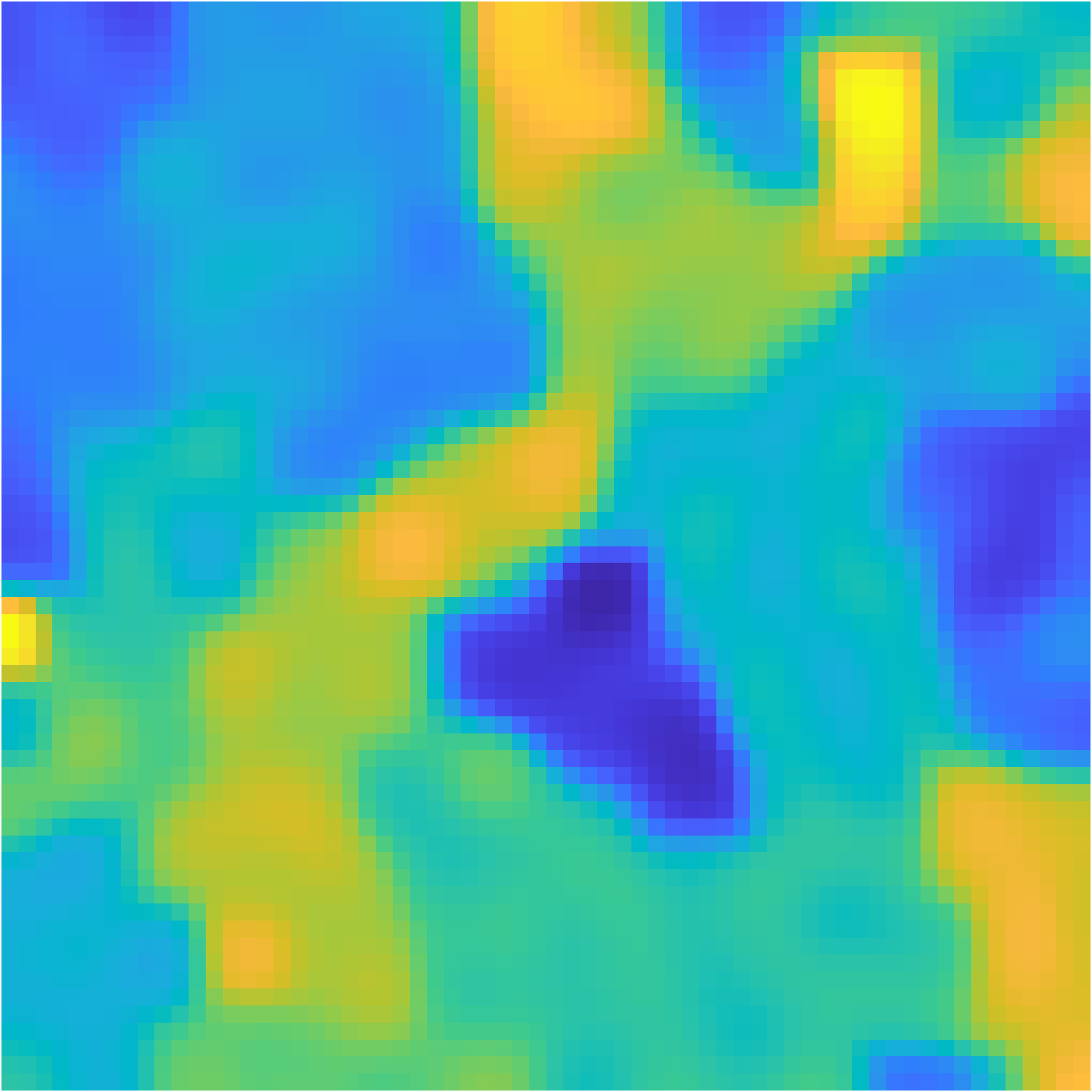}
\caption{I}
\end{subfigure}\hfill
\begin{subfigure}[t]{0.19\textwidth}
\centering
\includegraphics[width=\linewidth]{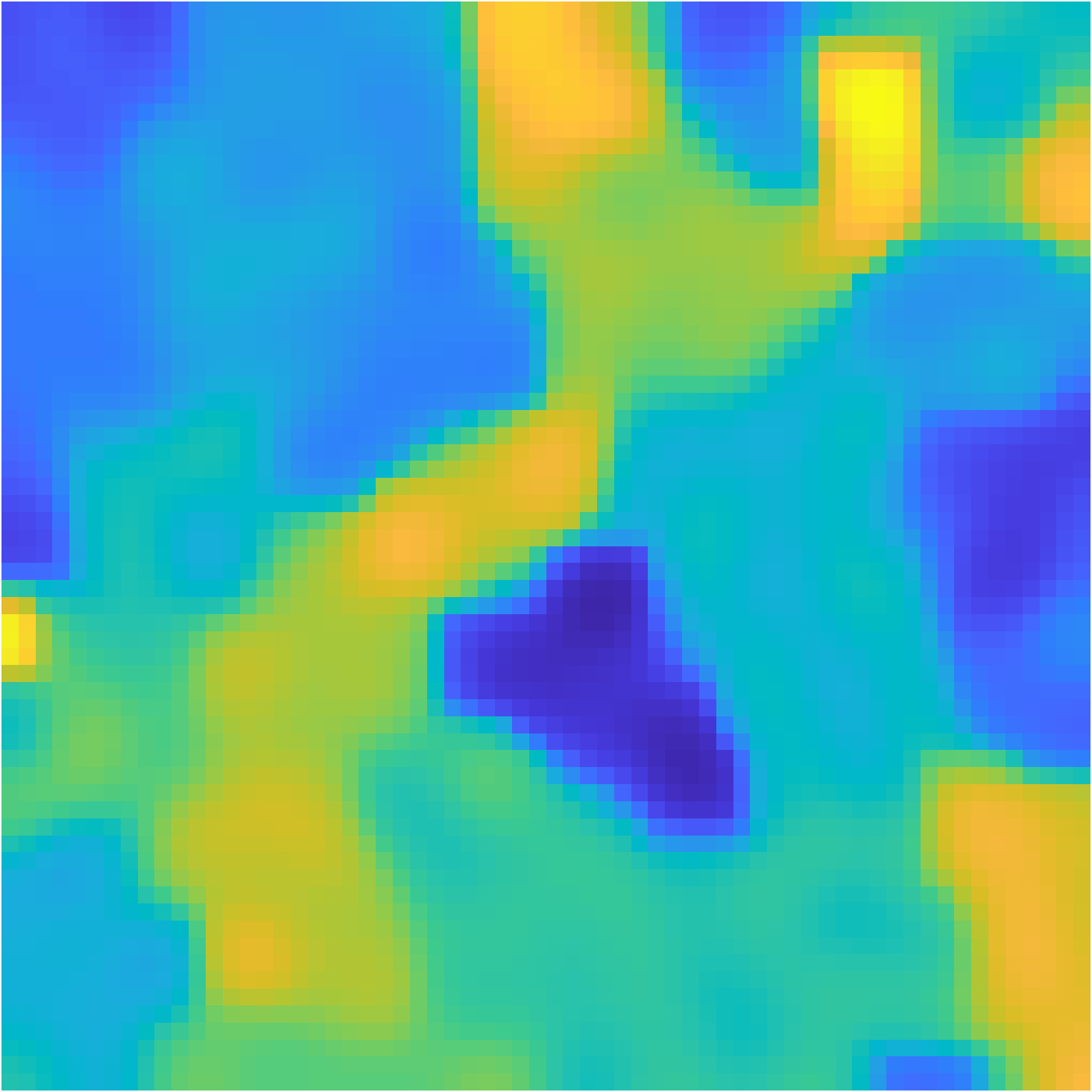}
\caption{HWT-D}
\end{subfigure}\hfill
\begin{subfigure}[t]{0.19\textwidth}
\centering
\includegraphics[width=\linewidth]{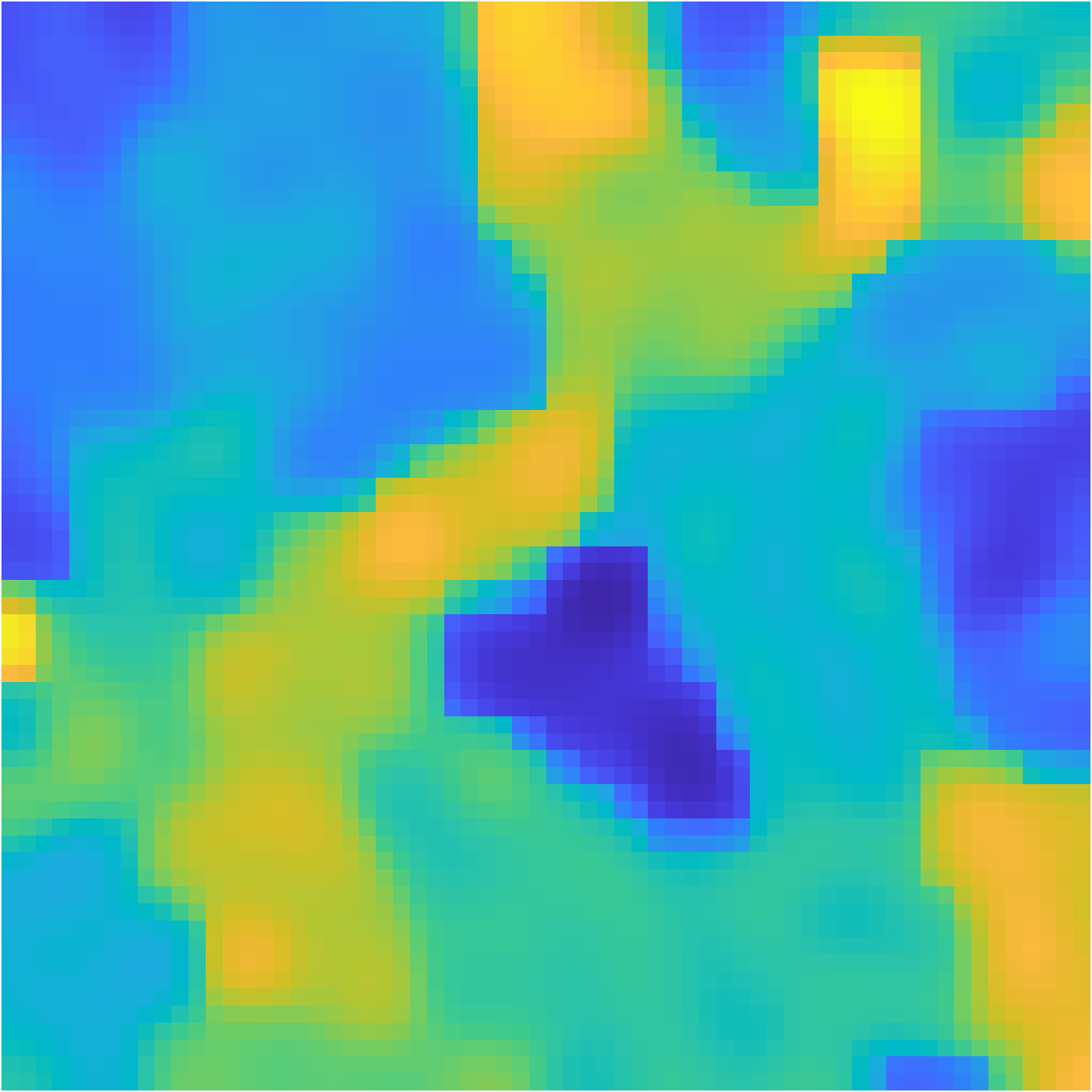}
\caption{HWT-X}
\end{subfigure}\hfill
\begin{subfigure}[t]{0.19\textwidth}
\centering
\includegraphics[width=\linewidth]{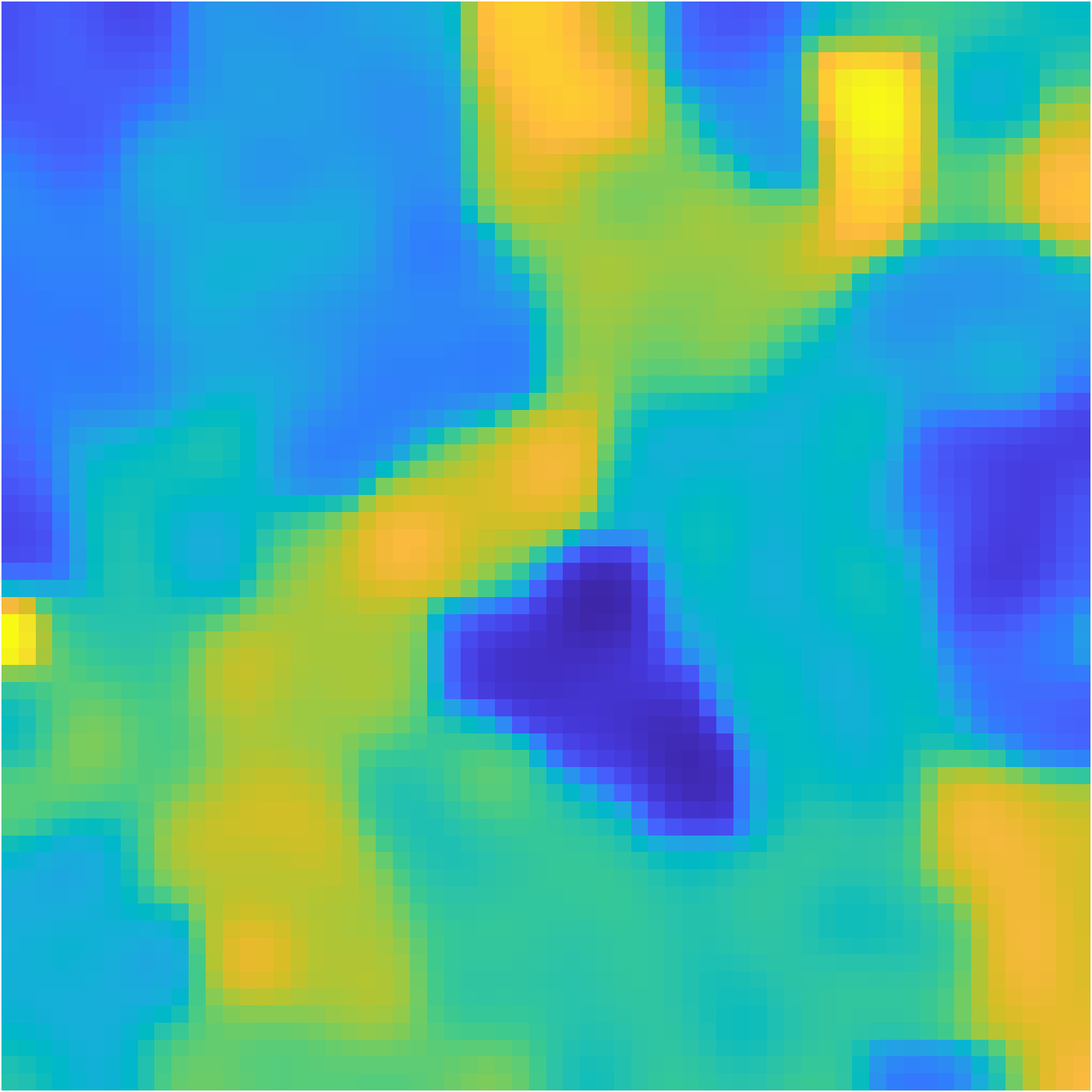}
\caption{DWT-D}
\end{subfigure}\hfill
\begin{subfigure}[t]{0.19\textwidth}
\centering
\includegraphics[width=\linewidth]{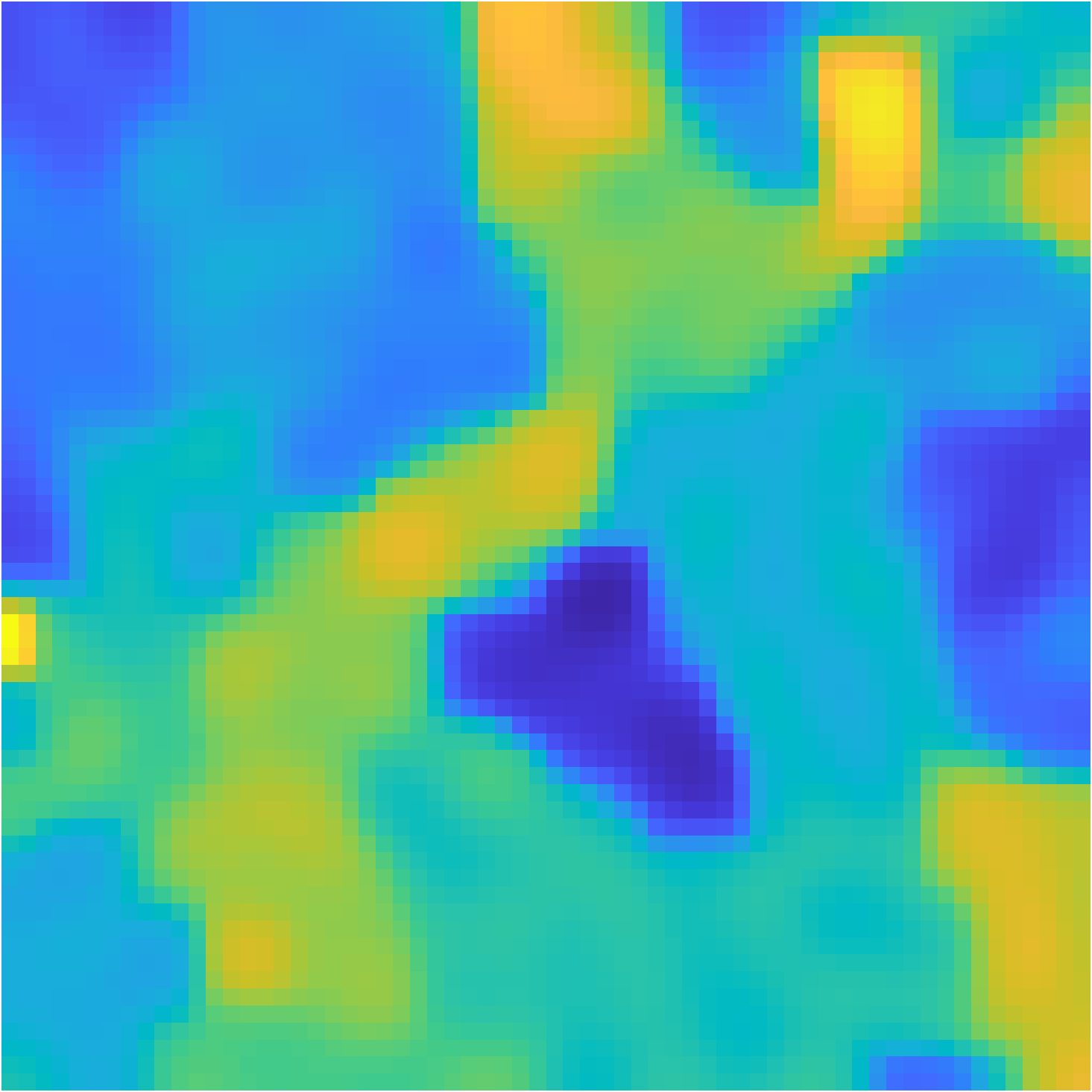}
\caption{DWT-X}
\end{subfigure}
\vspace{0.8em}
\begin{subfigure}[t]{0.19\textwidth}
\centering
\includegraphics[width=\linewidth]{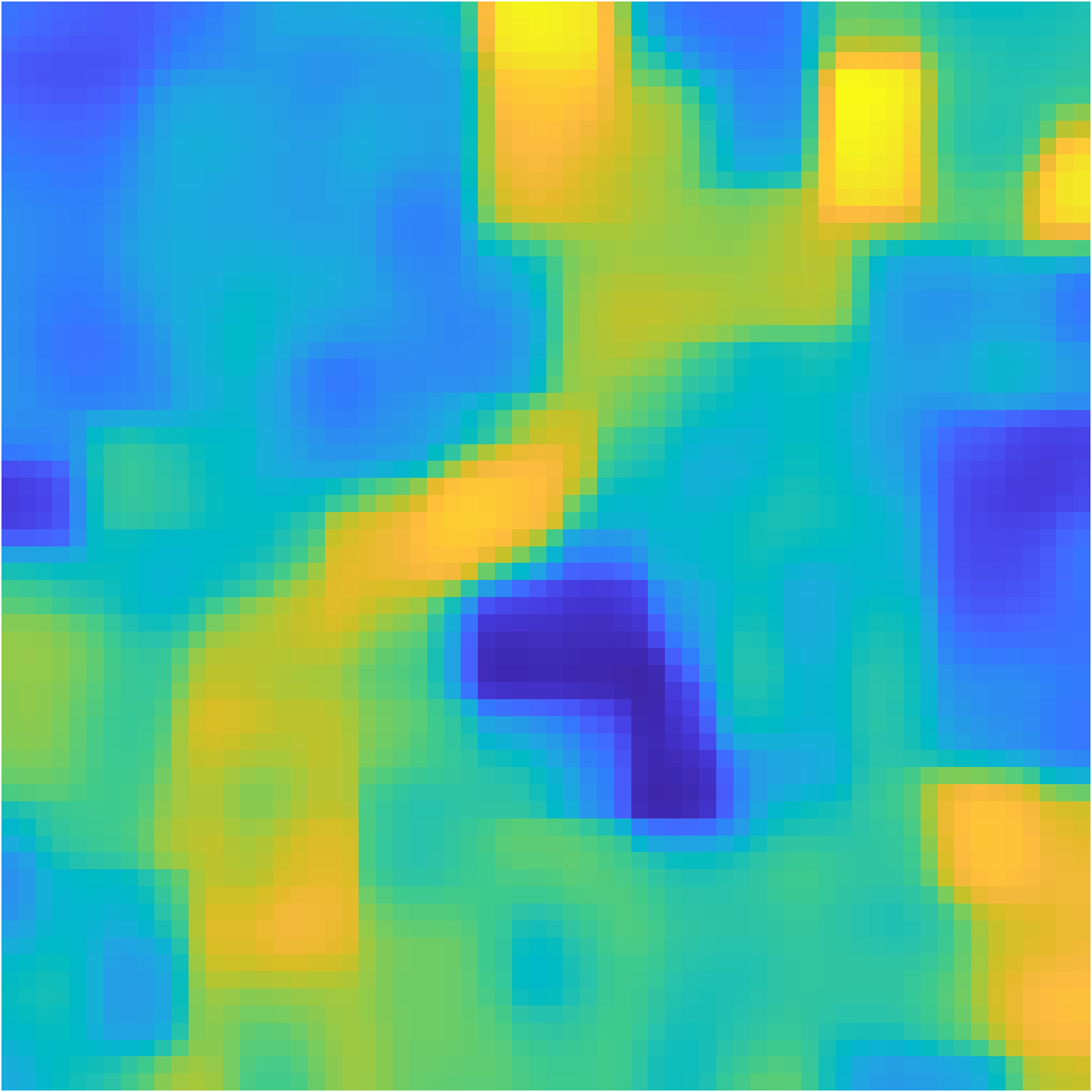}
\caption{I}
\end{subfigure}\hfill
\begin{subfigure}[t]{0.19\textwidth}
\centering
\includegraphics[width=\linewidth]{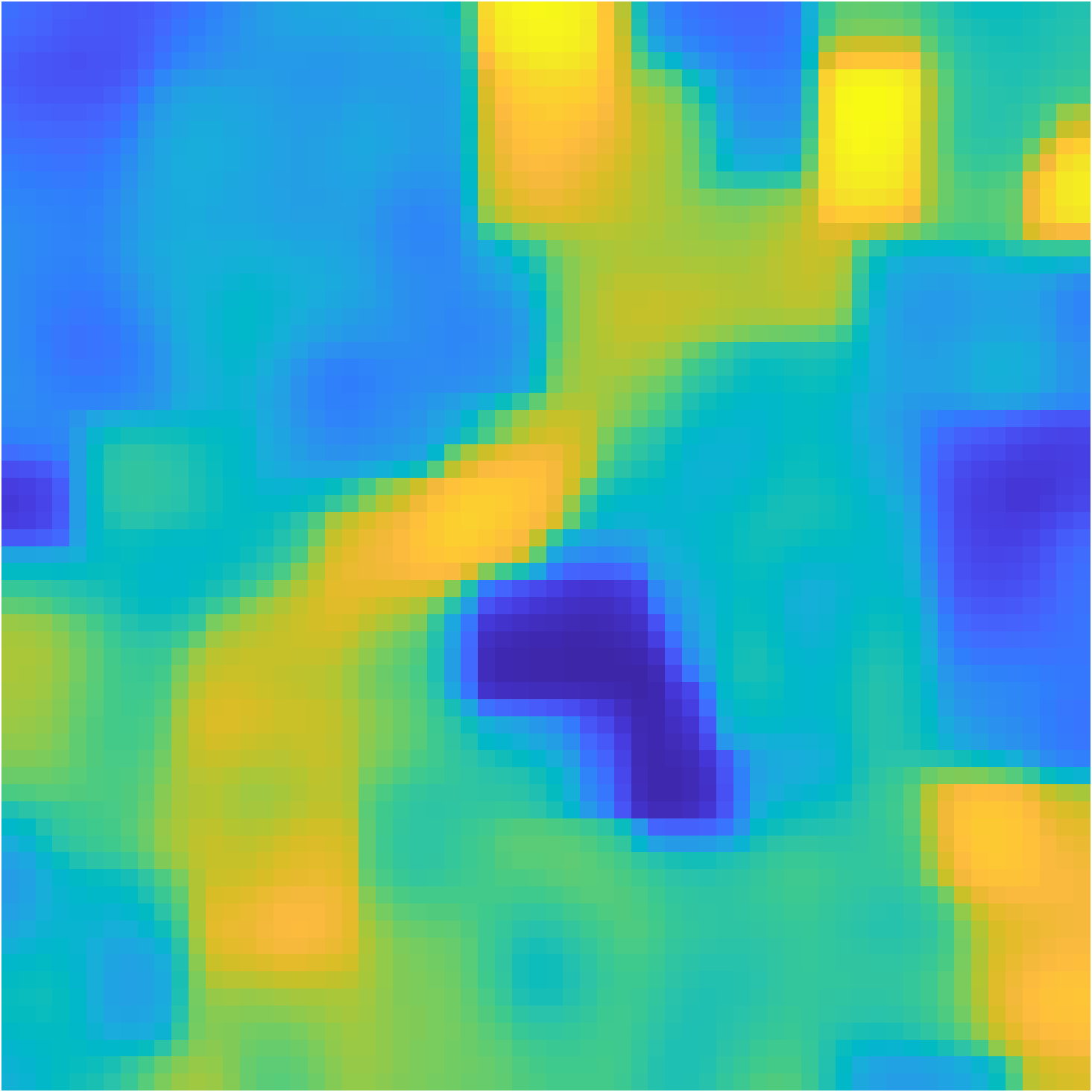}
\caption{HWT-D}
\end{subfigure}\hfill
\begin{subfigure}[t]{0.19\textwidth}
\centering
\includegraphics[width=\linewidth]{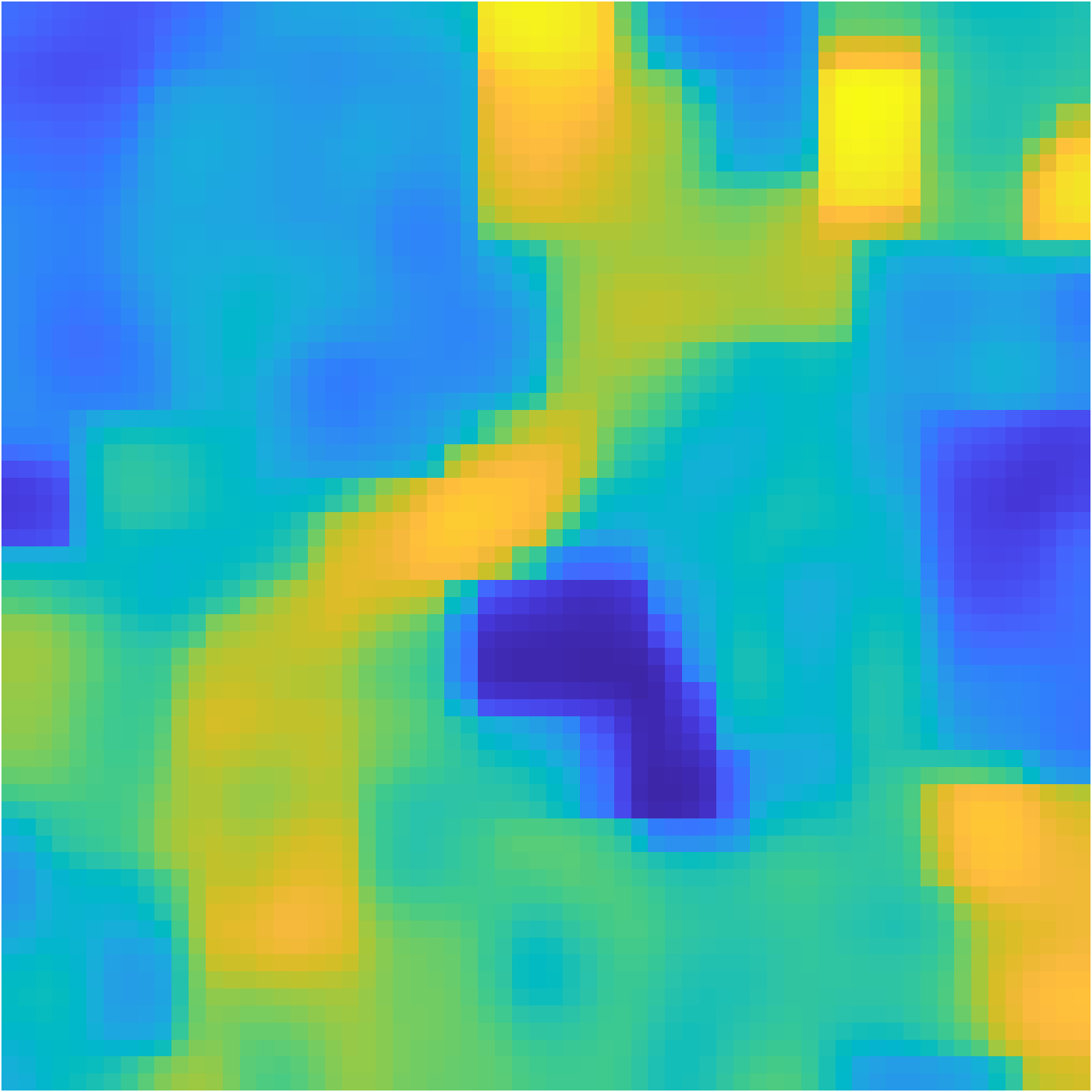}
\caption{HWT-X}
\end{subfigure}\hfill
\begin{subfigure}[t]{0.19\textwidth}
\centering
\includegraphics[width=\linewidth]{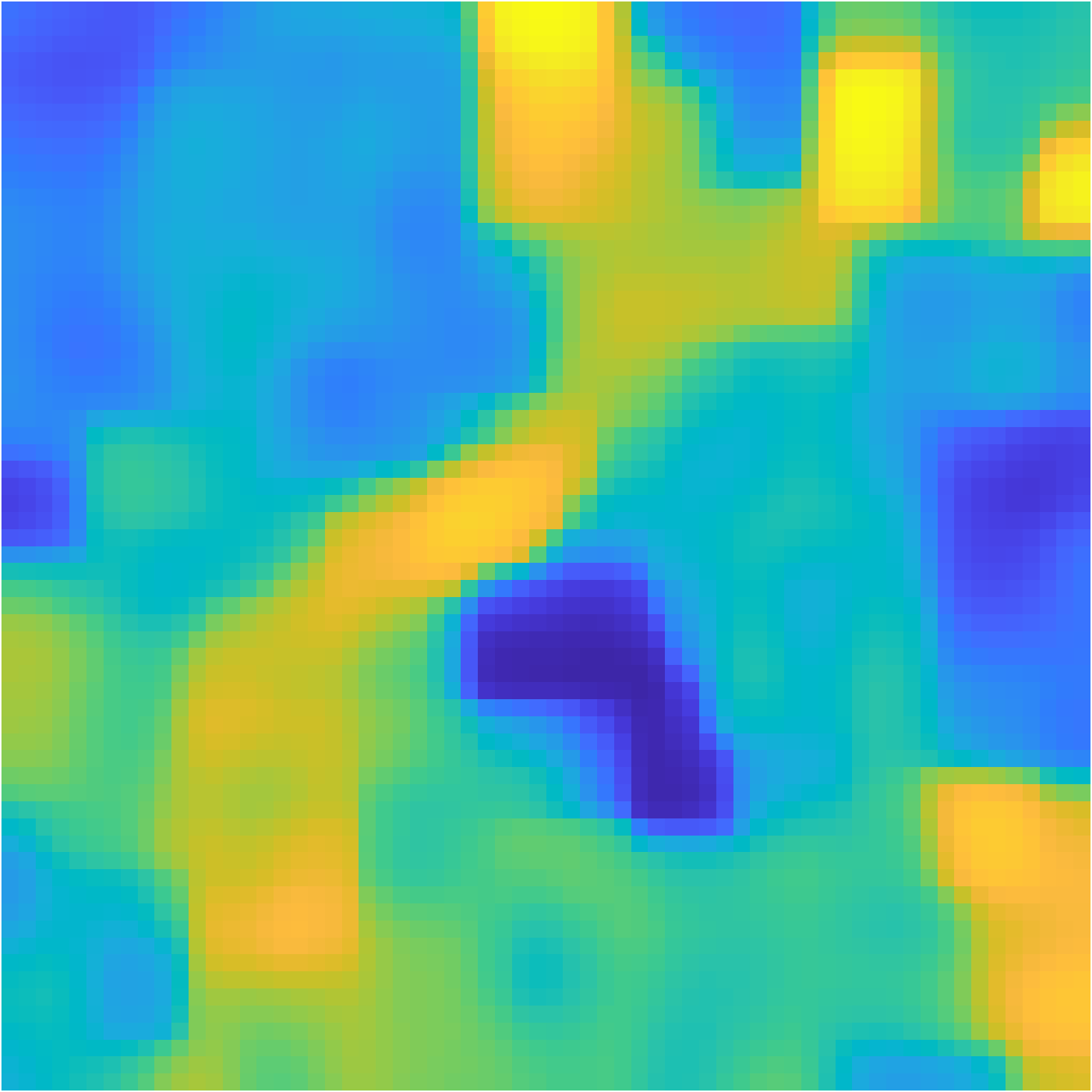}
\caption{DWT-D}
\end{subfigure}\hfill
\begin{subfigure}[t]{0.19\textwidth}
\centering
\includegraphics[width=\linewidth]{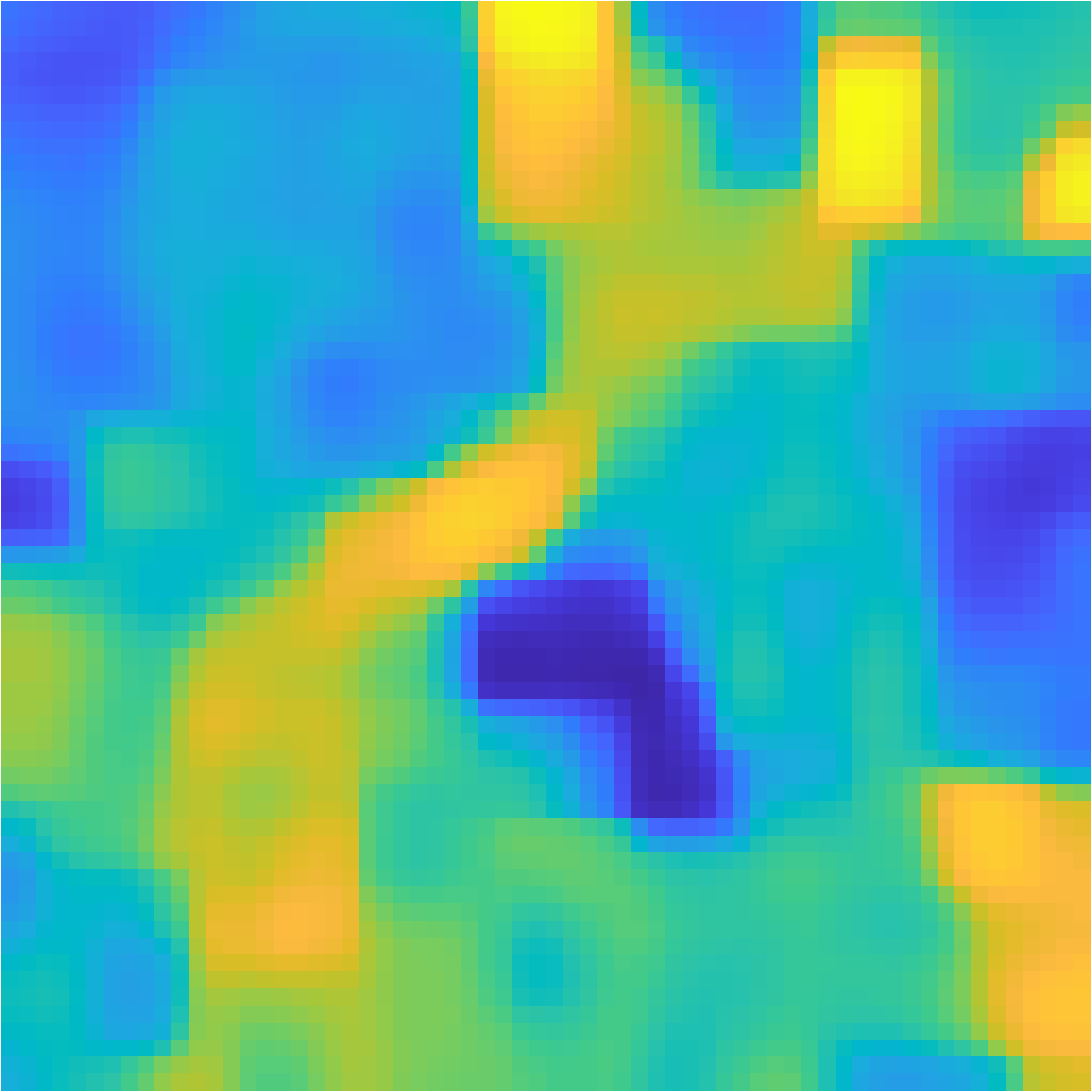}
\caption{DWT-X}
\end{subfigure}
\caption{Deblurred images using MM. Top row: $\sigma=3.5$; bottom row: $\sigma=4.5$.}
\label{fig:SolutionMM}
\end{figure}

\begin{figure}[tbp]
\centering
\begin{subfigure}[t]{0.19\textwidth}
\centering
\includegraphics[width=\linewidth]{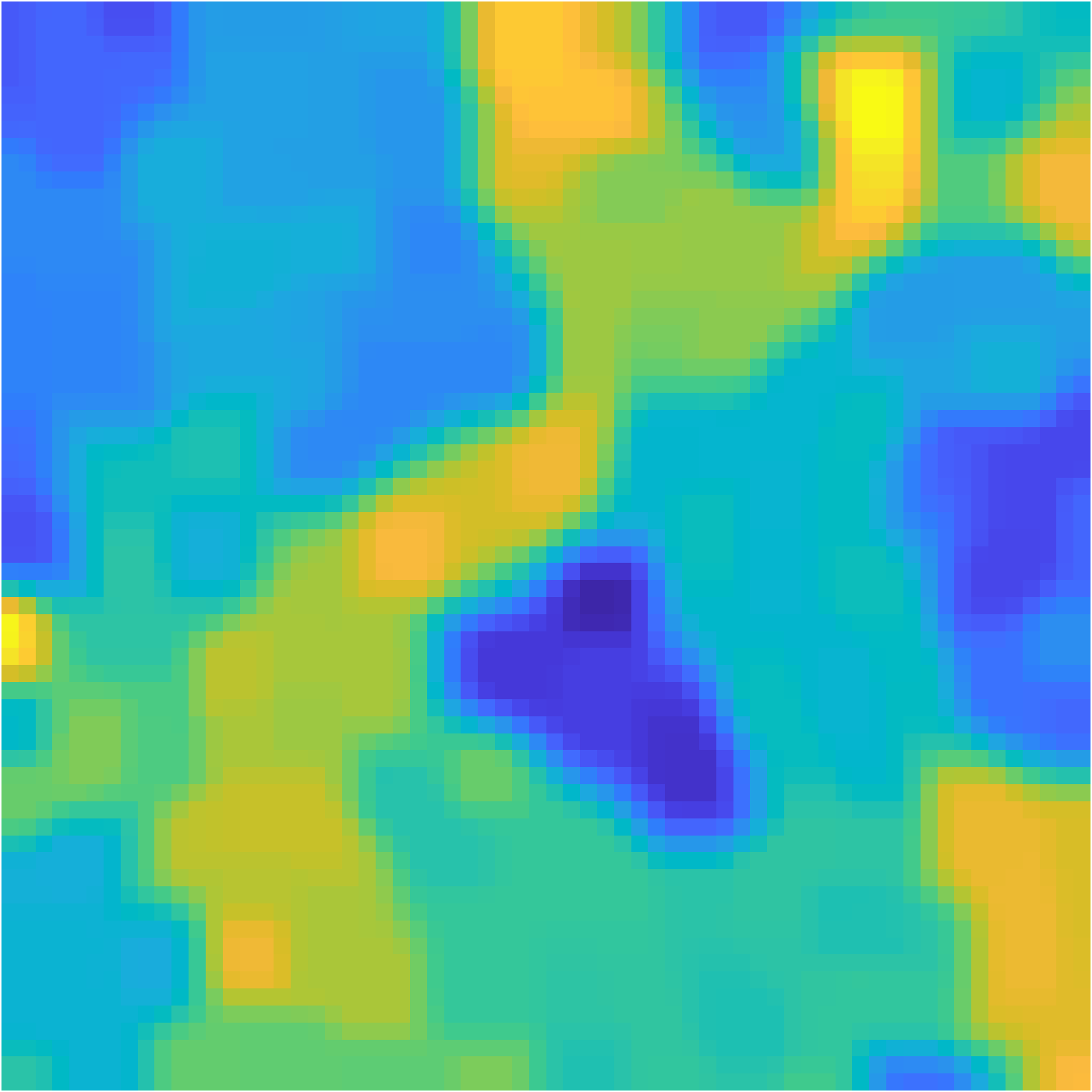}
\caption{I}
\end{subfigure}\hfill
\begin{subfigure}[t]{0.19\textwidth}
\centering
\includegraphics[width=\linewidth]{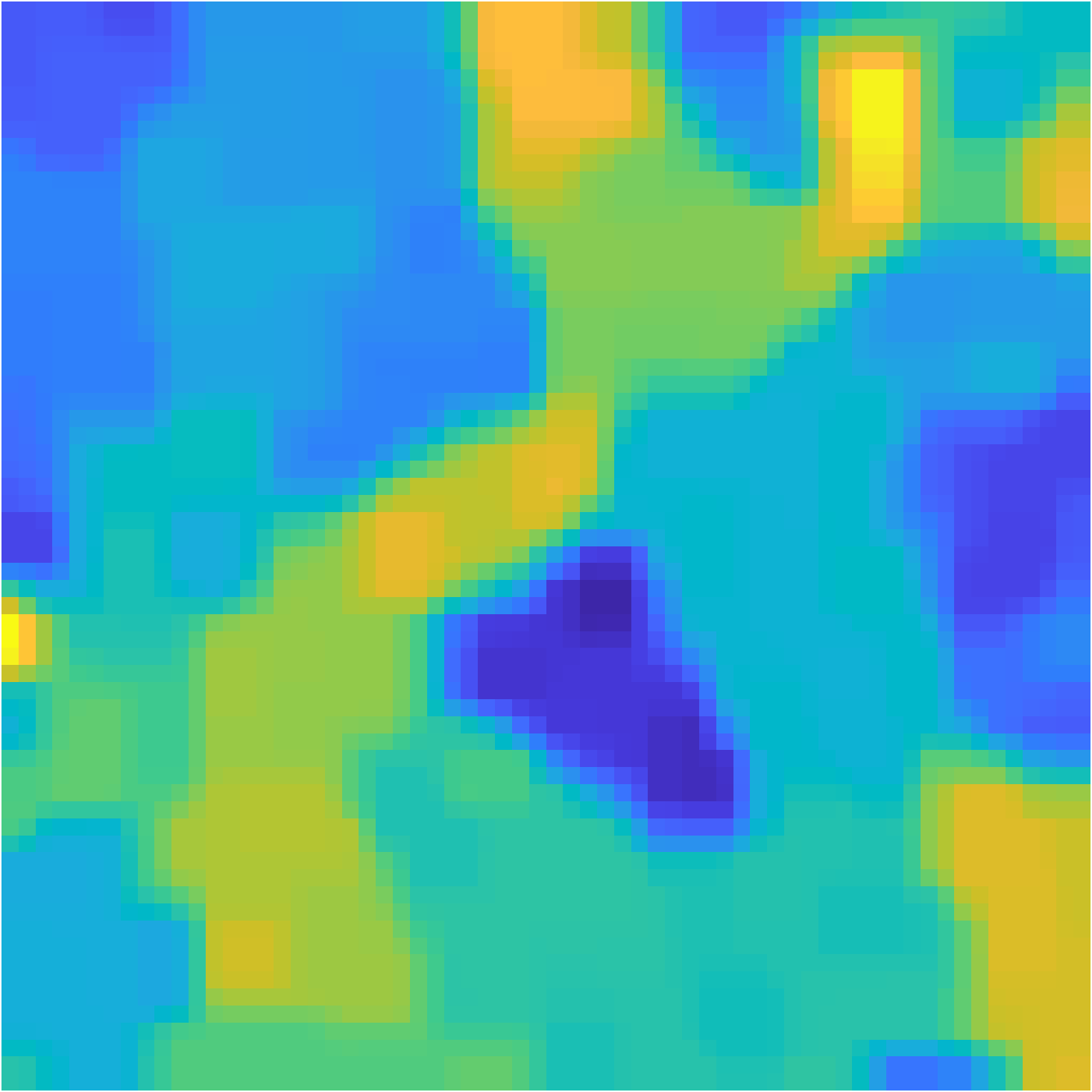}
\caption{HWT-D}
\end{subfigure}\hfill
\begin{subfigure}[t]{0.19\textwidth}
\centering
\includegraphics[width=\linewidth]{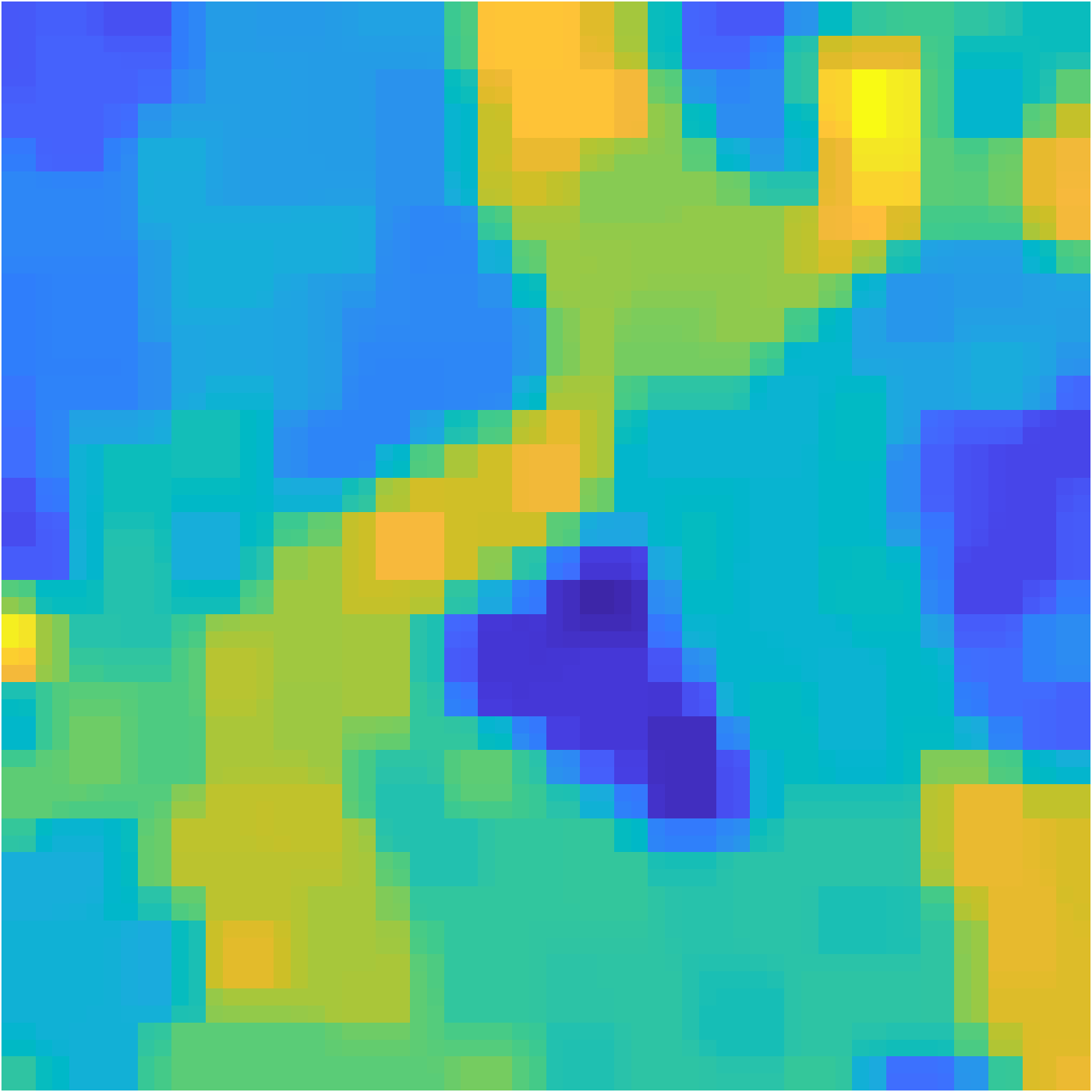}
\caption{HWT-X}
\end{subfigure}\hfill
\begin{subfigure}[t]{0.19\textwidth}
\centering
\includegraphics[width=\linewidth]{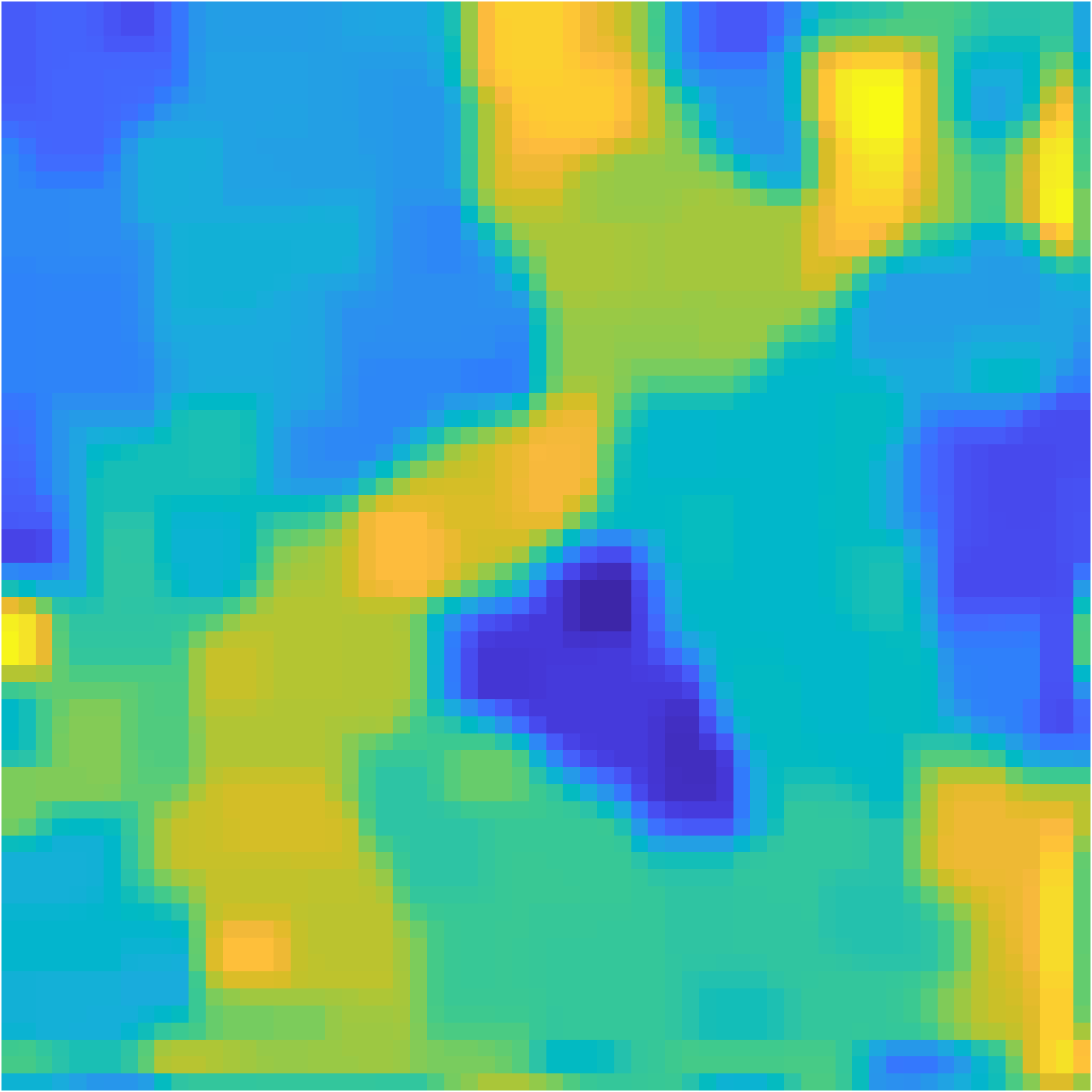}
\caption{DWT-D}
\end{subfigure}\hfill
\begin{subfigure}[t]{0.19\textwidth}
\centering
\includegraphics[width=\linewidth]{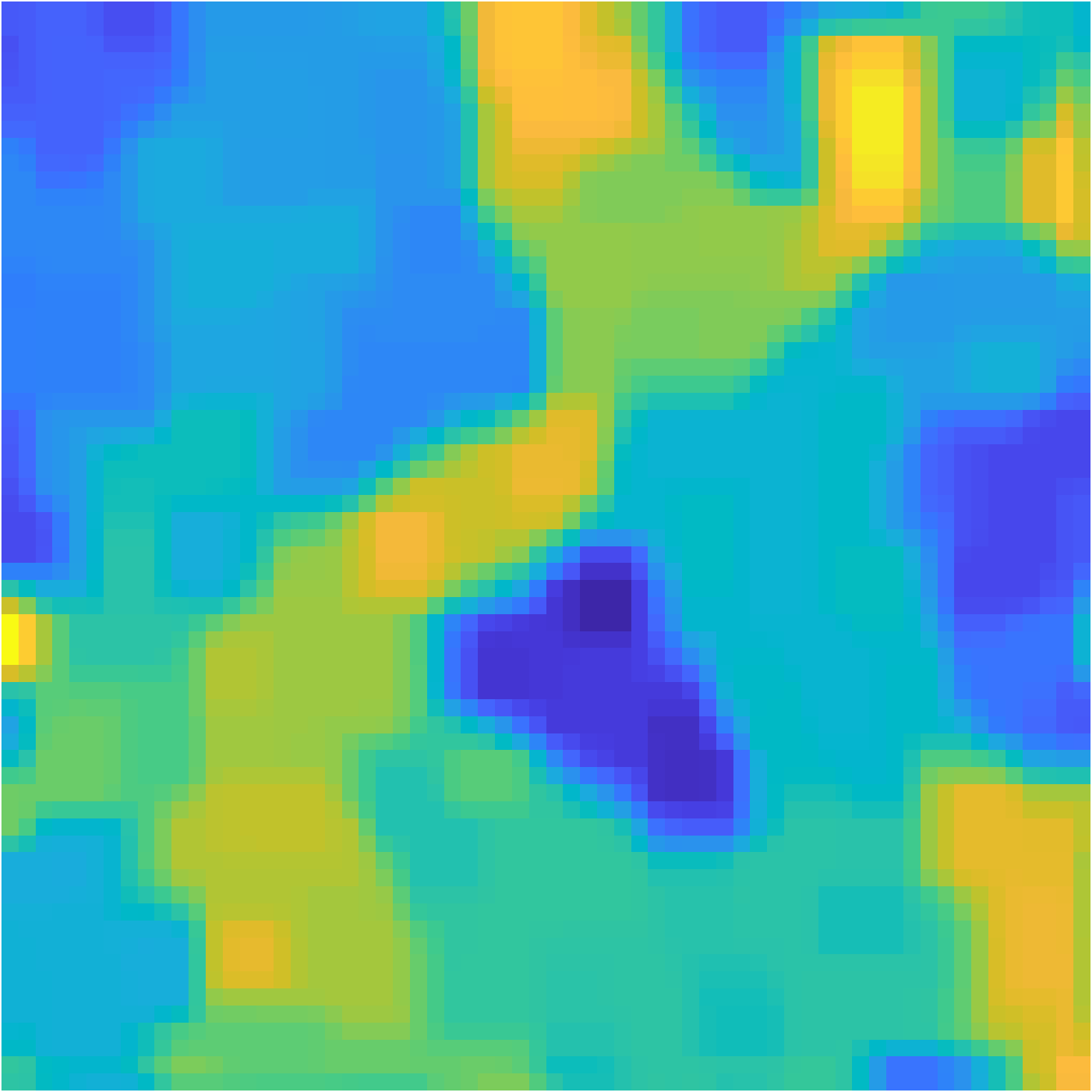}
\caption{DWT-X}
\end{subfigure}
\vspace{0.8em}
\begin{subfigure}[t]{0.19\textwidth}
\centering
\includegraphics[width=\linewidth]{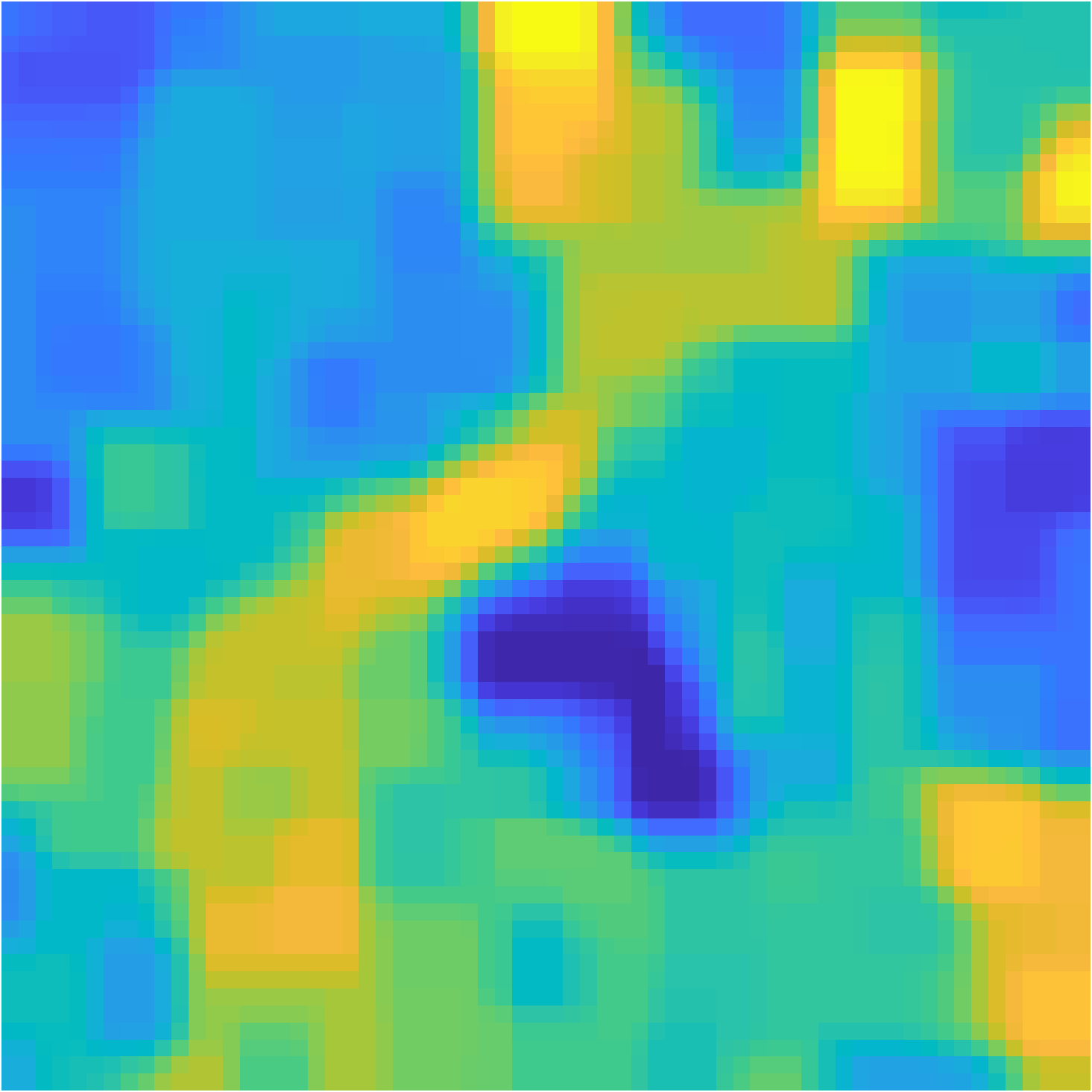}
\caption{I}
\end{subfigure}\hfill
\begin{subfigure}[t]{0.19\textwidth}
\centering
\includegraphics[width=\linewidth]{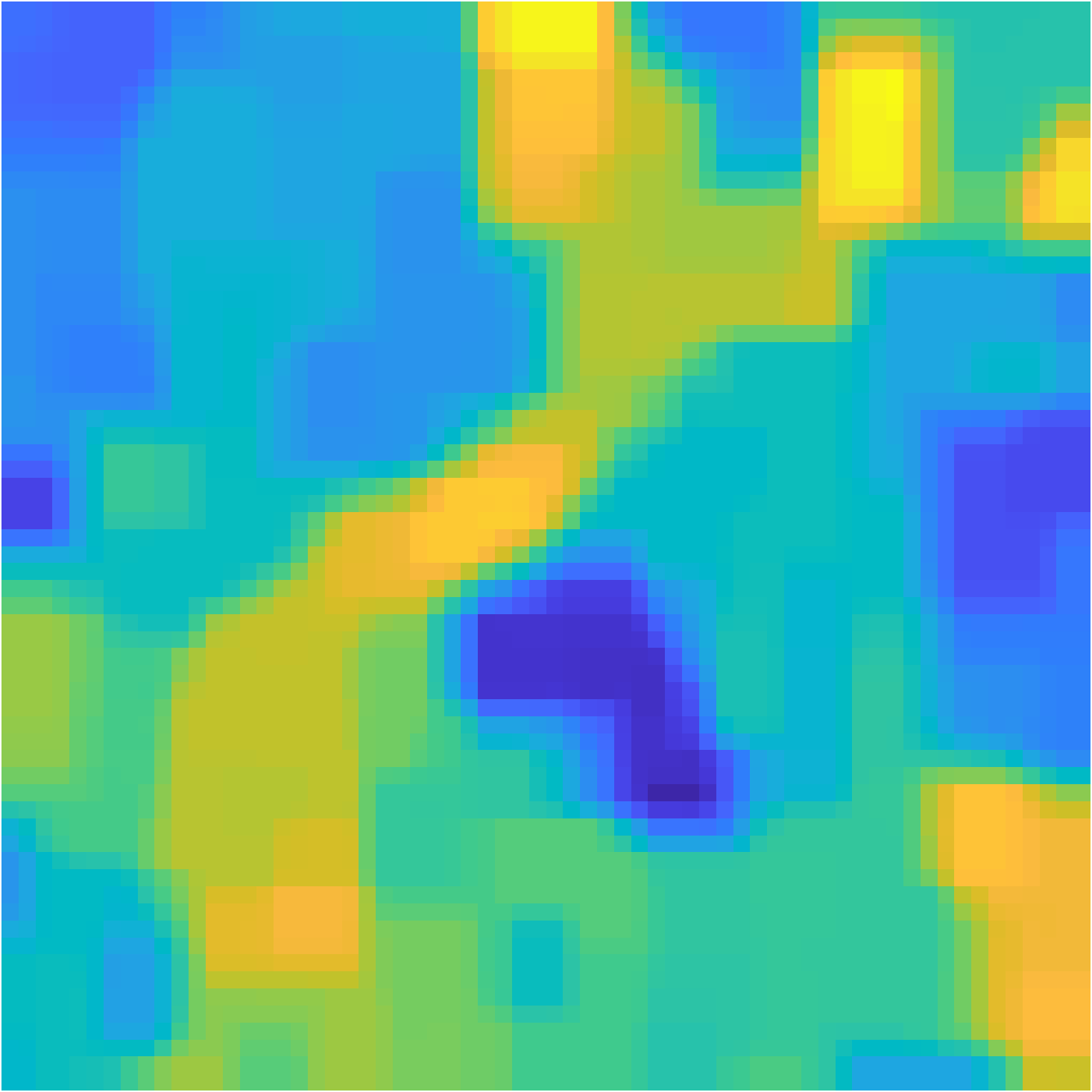}
\caption{HWT-D}
\end{subfigure}\hfill
\begin{subfigure}[t]{0.19\textwidth}
\centering
\includegraphics[width=\linewidth]{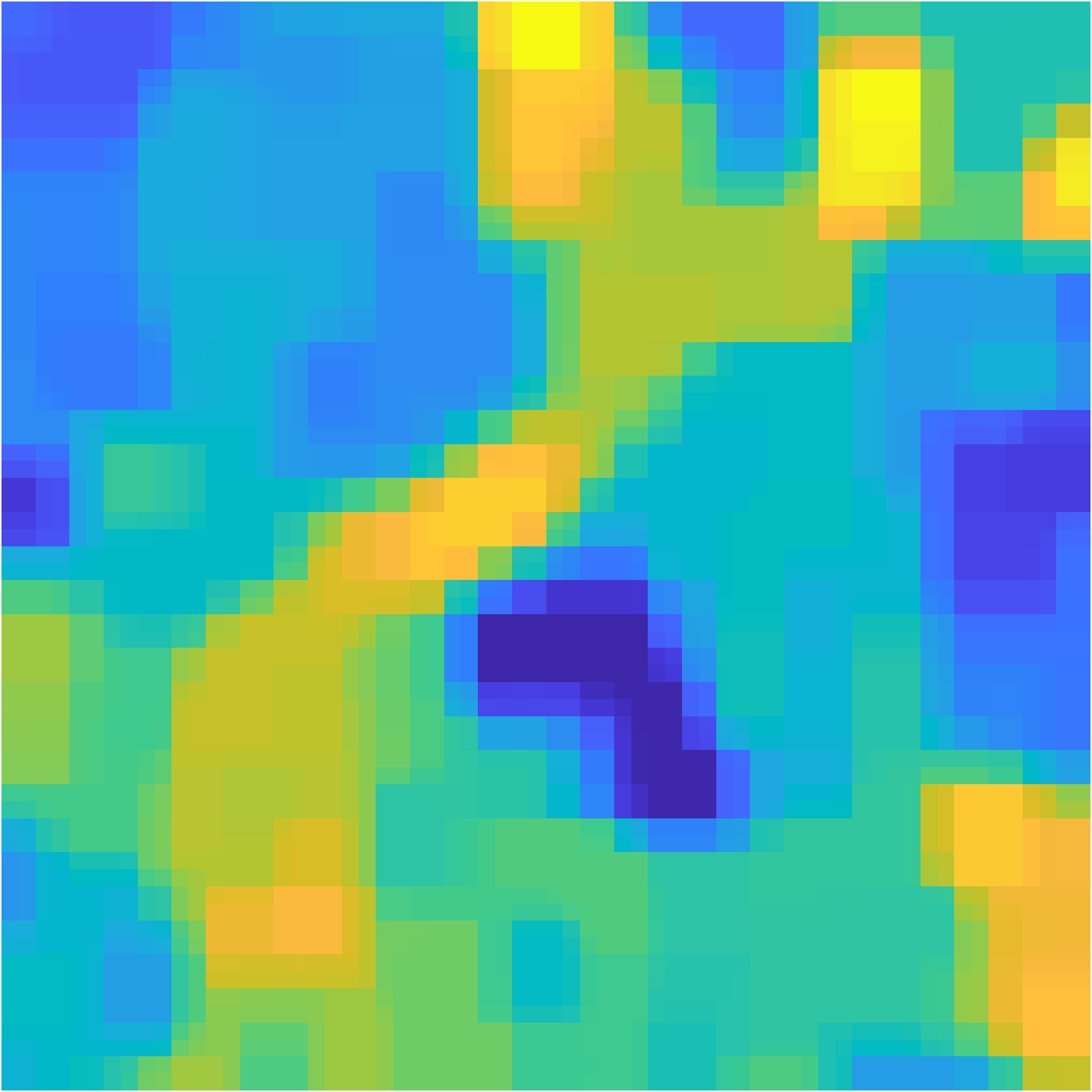}
\caption{HWT-X}
\end{subfigure}\hfill
\begin{subfigure}[t]{0.19\textwidth}
\centering
\includegraphics[width=\linewidth]{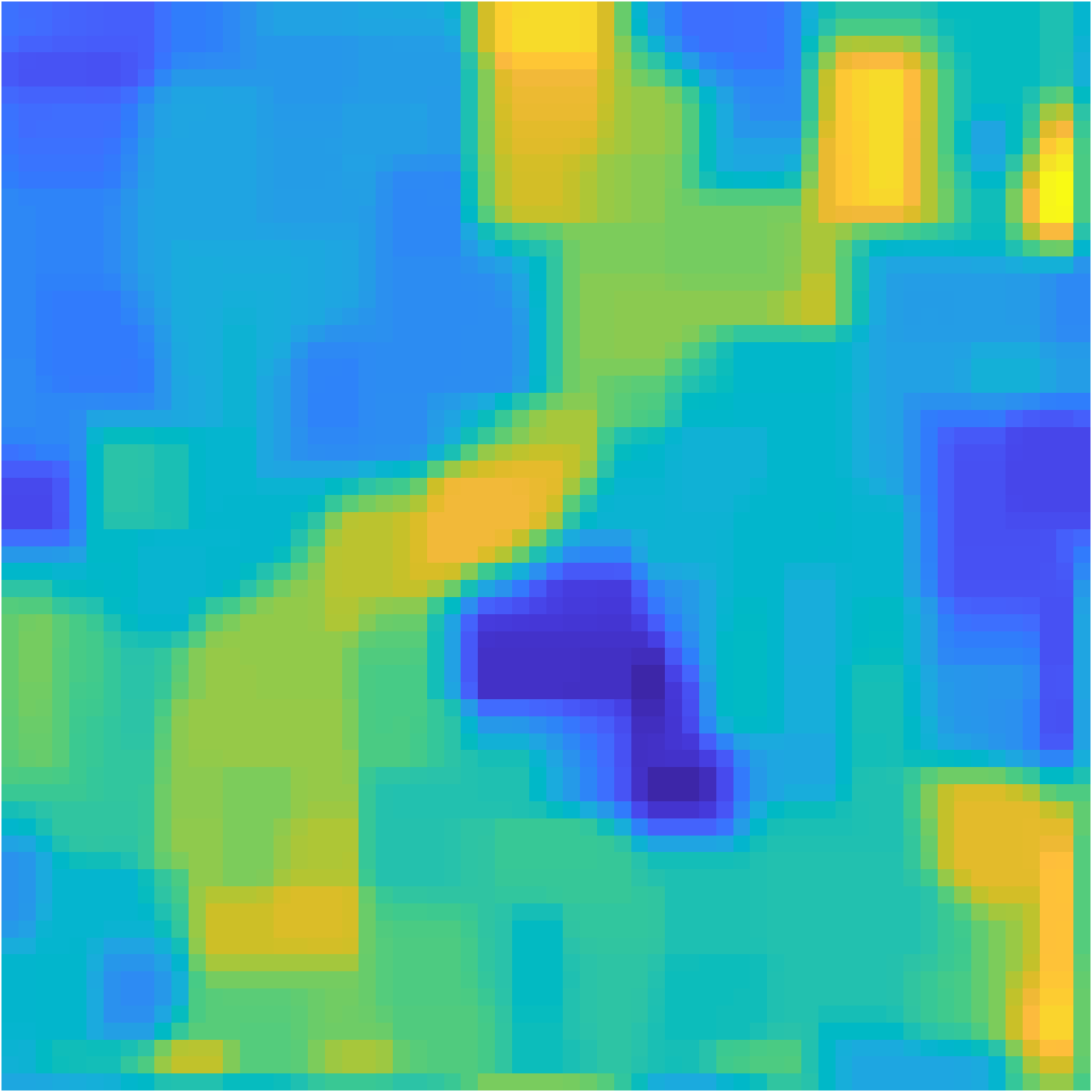}
\caption{DWT-D}
\end{subfigure}\hfill
\begin{subfigure}[t]{0.19\textwidth}
\centering
\includegraphics[width=\linewidth]{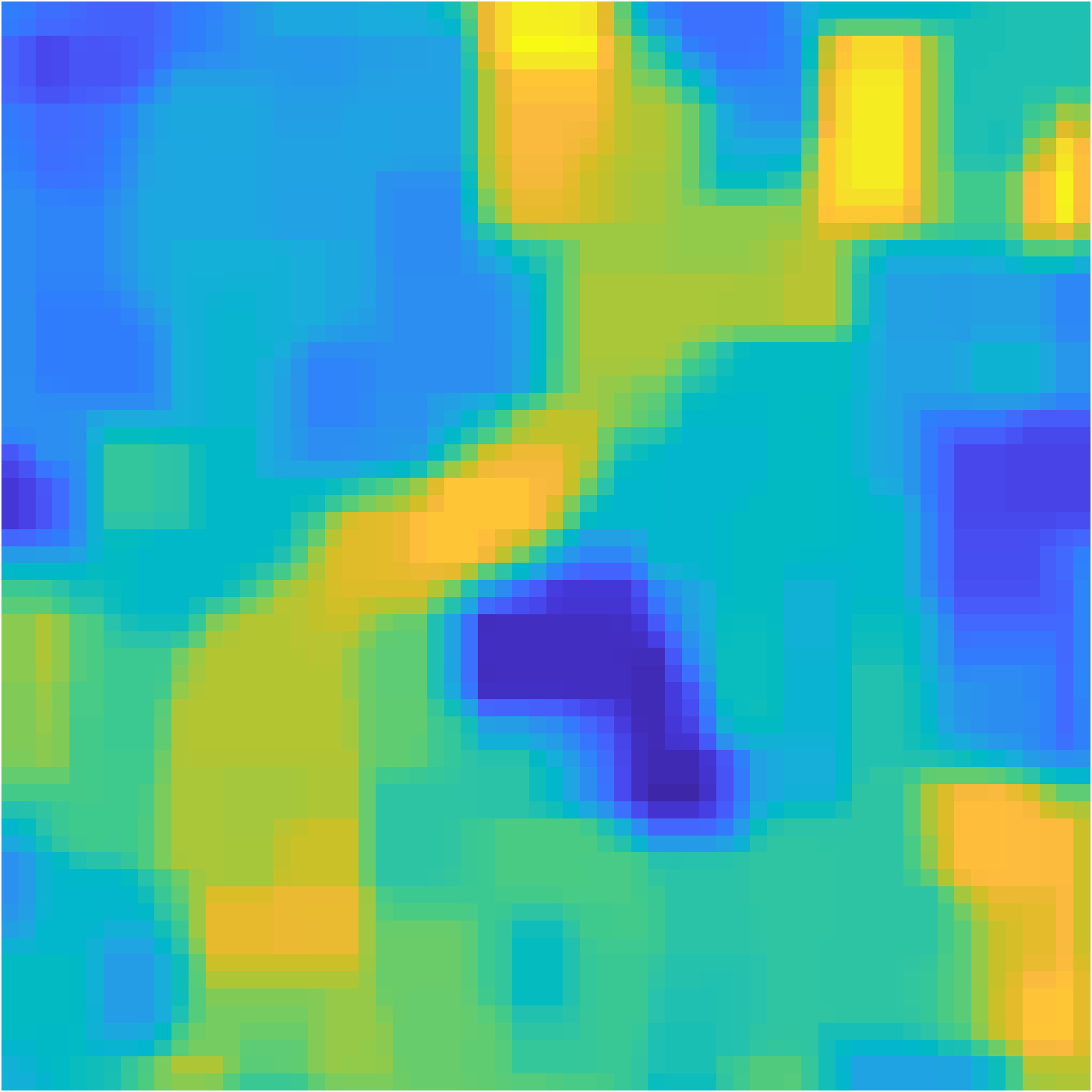}
\caption{DWT-X}
\end{subfigure}
\caption{Deblurred images using SB. Top row: $\sigma=3.5$; bottom row: $\sigma=4.5$.}
\label{fig:SolutionSB}
\end{figure}

\begin{figure}[tbp]
\centering
\begin{subfigure}[t]{0.45\textwidth}
\centering
\includegraphics[width=\linewidth]{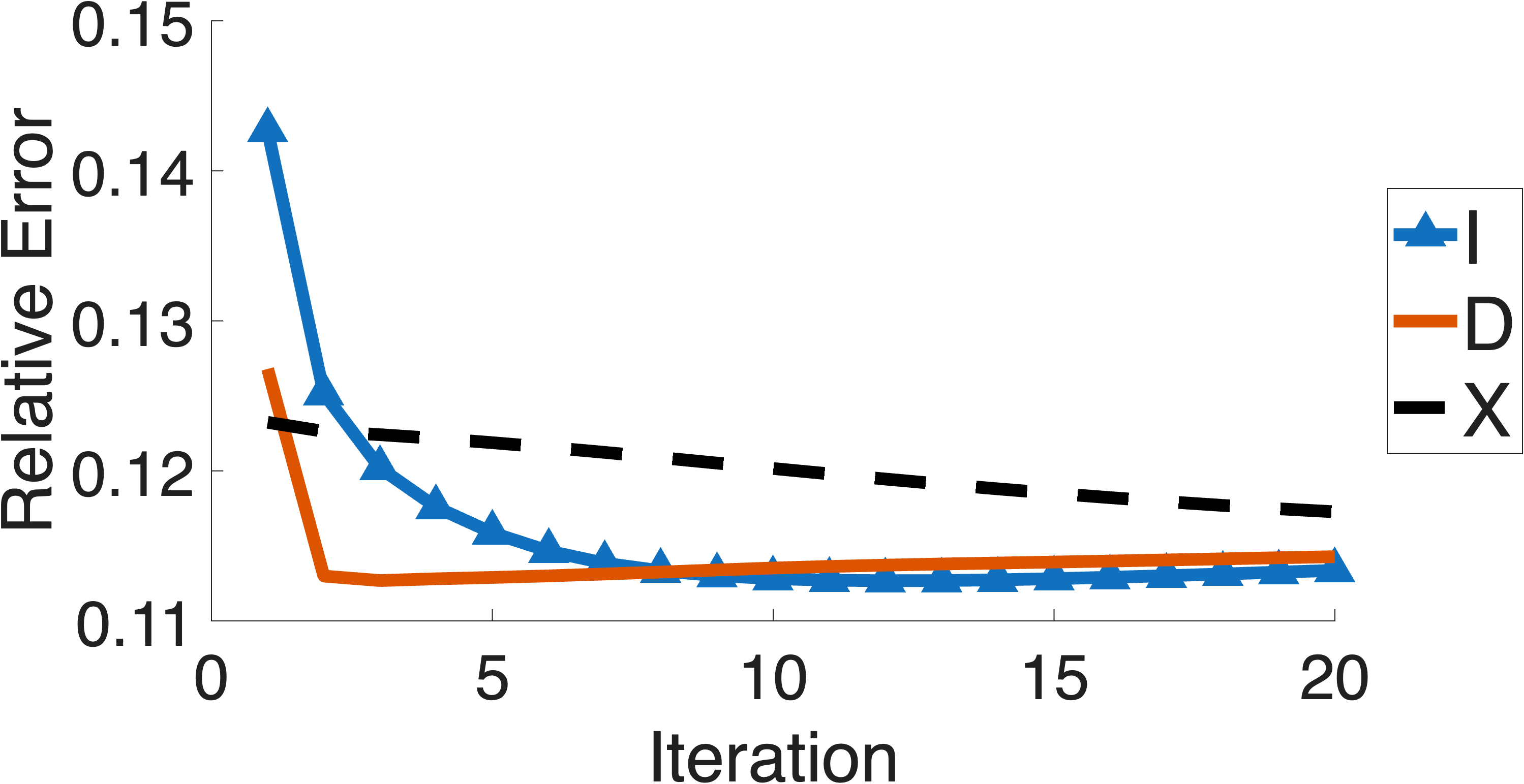}
\caption{HWT, $\sigma=3.5$}
\end{subfigure}\hfill
\begin{subfigure}[t]{0.45\textwidth}
\centering
\includegraphics[width=\linewidth]{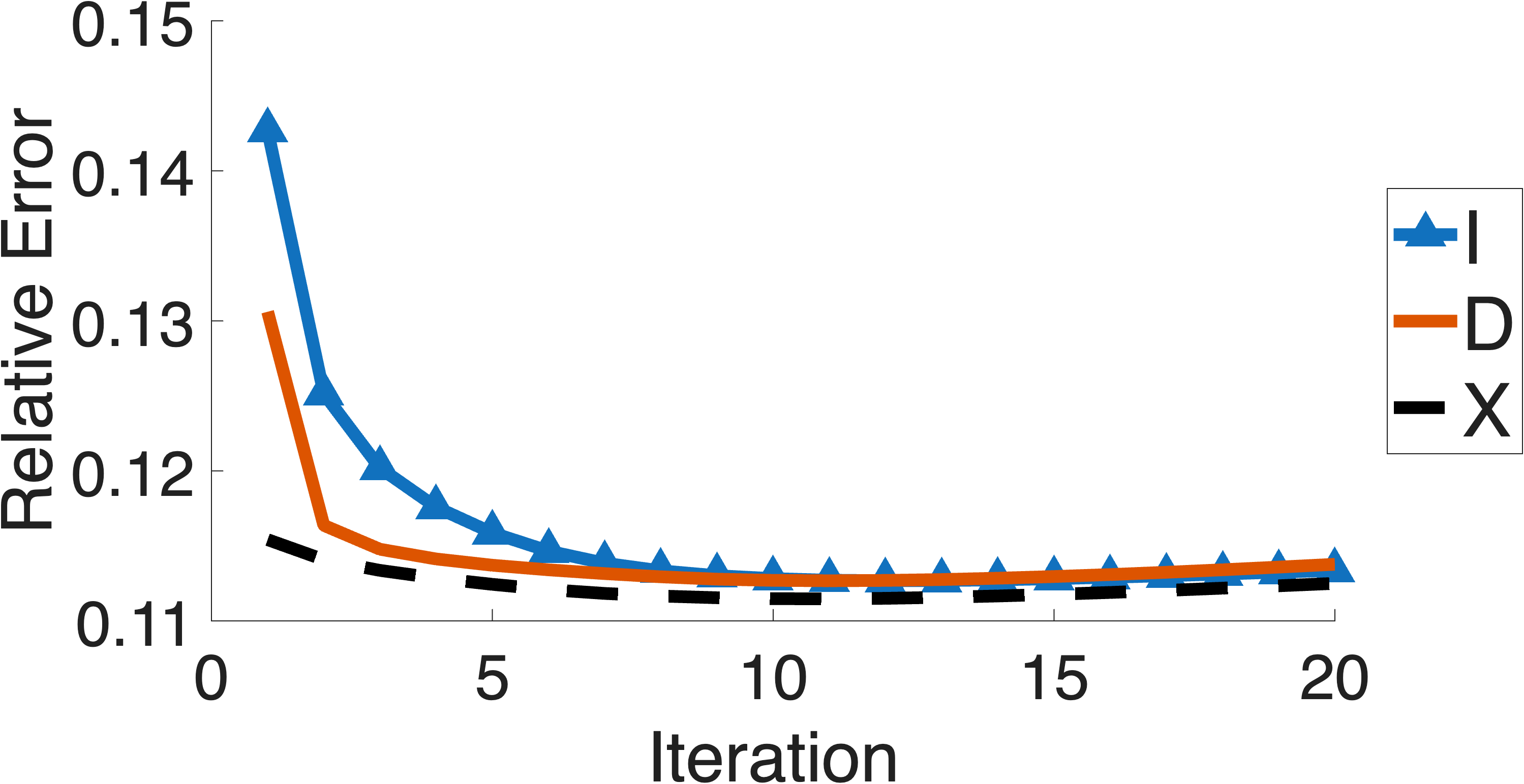}
\caption{DWT, $\sigma=3.5$}
\end{subfigure}
\vspace{0.8em}
\begin{subfigure}[t]{0.45\textwidth}
\centering
\includegraphics[width=\linewidth]{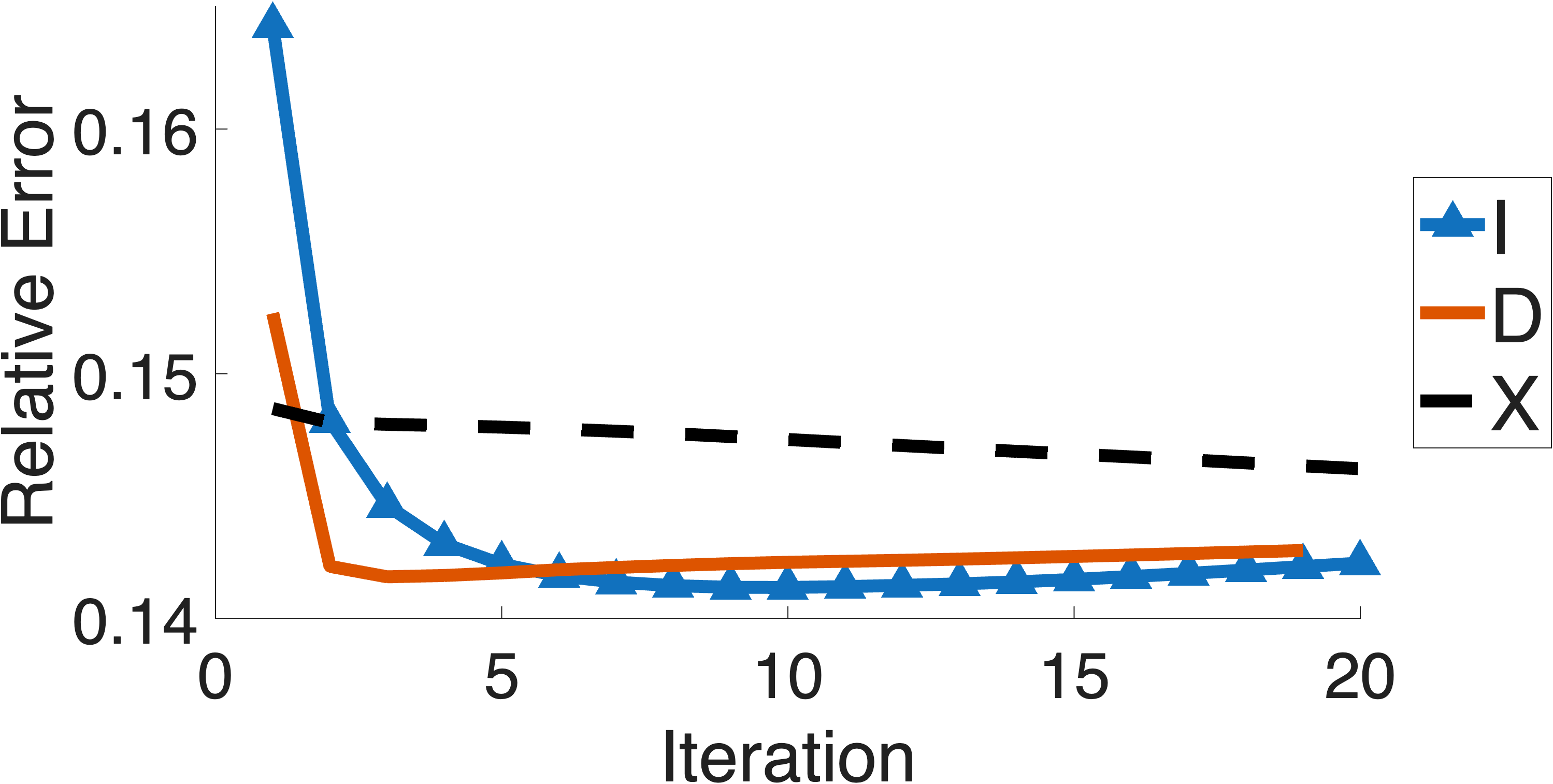}
\caption{HWT, $\sigma=4.5$}
\end{subfigure}\hfill
\begin{subfigure}[t]{0.45\textwidth}
\centering
\includegraphics[width=\linewidth]{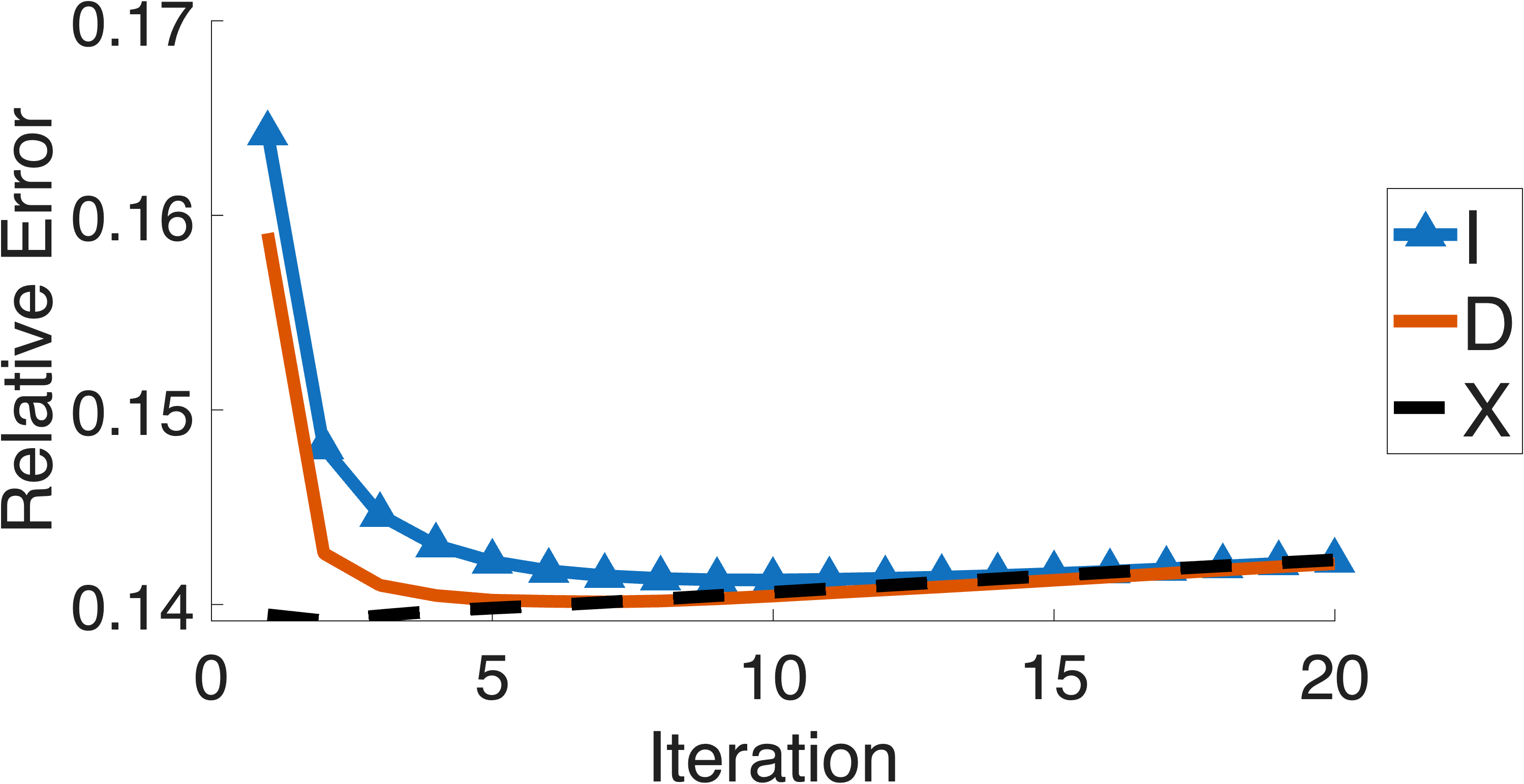}
\caption{DWT, $\sigma=4.5$}
\end{subfigure}
\caption{RRE convergence curves for IRLS (two-level).}
\label{fig:RRENirntv}
\end{figure}

\begin{figure}[tbp]
\centering
\begin{subfigure}[t]{0.45\textwidth}
\centering
\includegraphics[width=\linewidth]{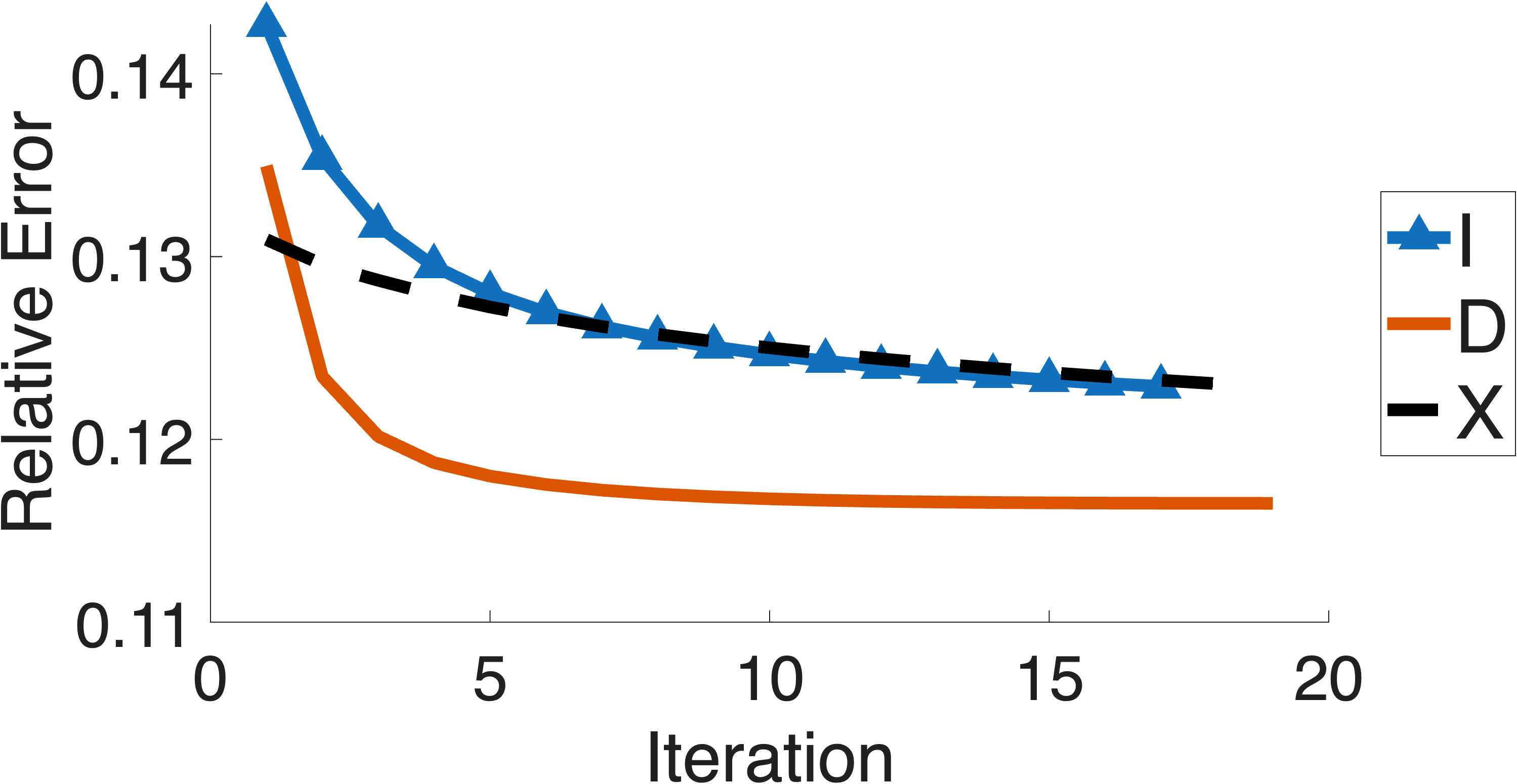}
\caption{HWT, $\sigma=3.5$}
\end{subfigure}\hfill
\begin{subfigure}[t]{0.45\textwidth}
\centering
\includegraphics[width=\linewidth]{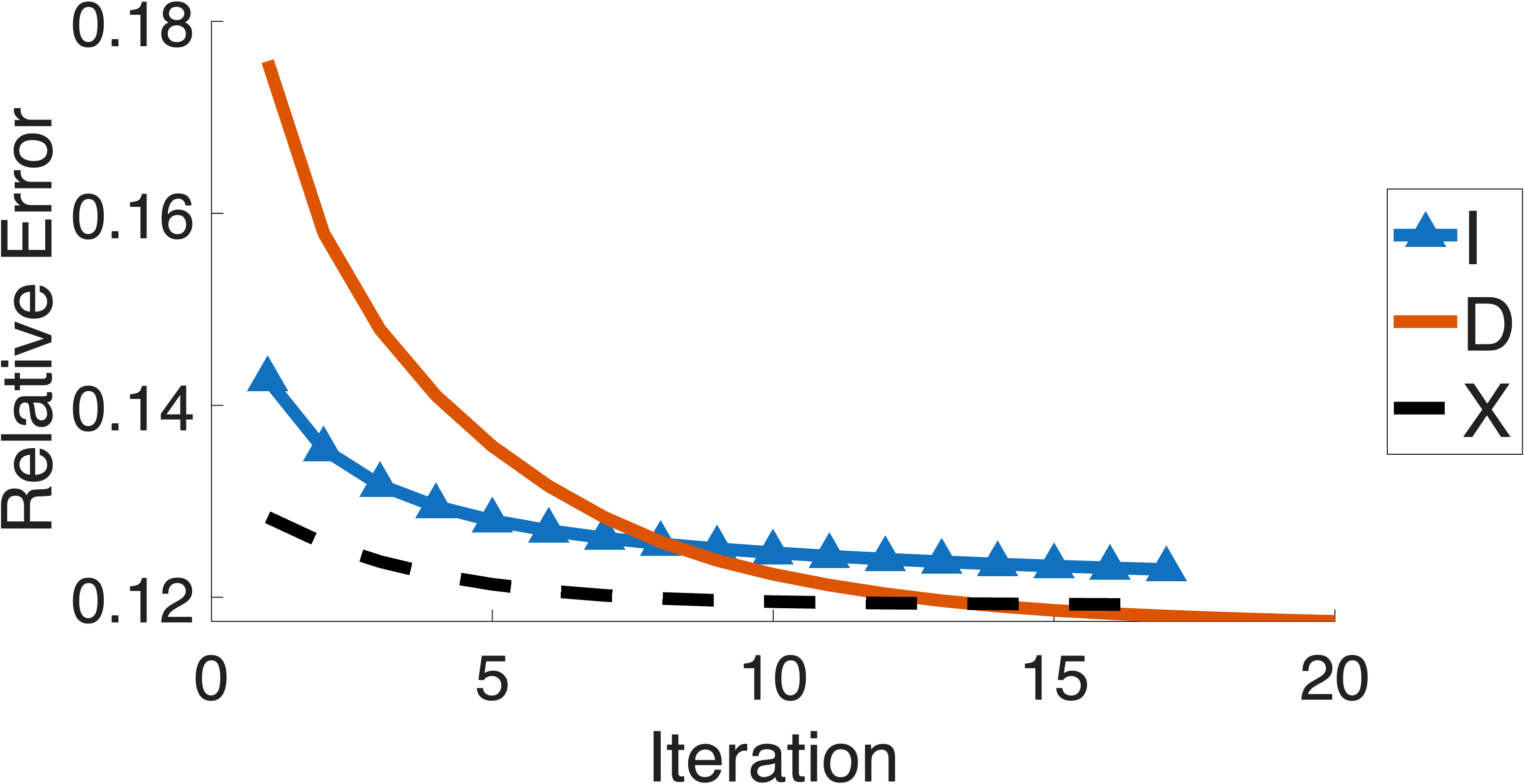}
\caption{DWT, $\sigma=3.5$}
\end{subfigure}
\vspace{0.8em}
\begin{subfigure}[t]{0.45\textwidth}
\centering
\includegraphics[width=\linewidth]{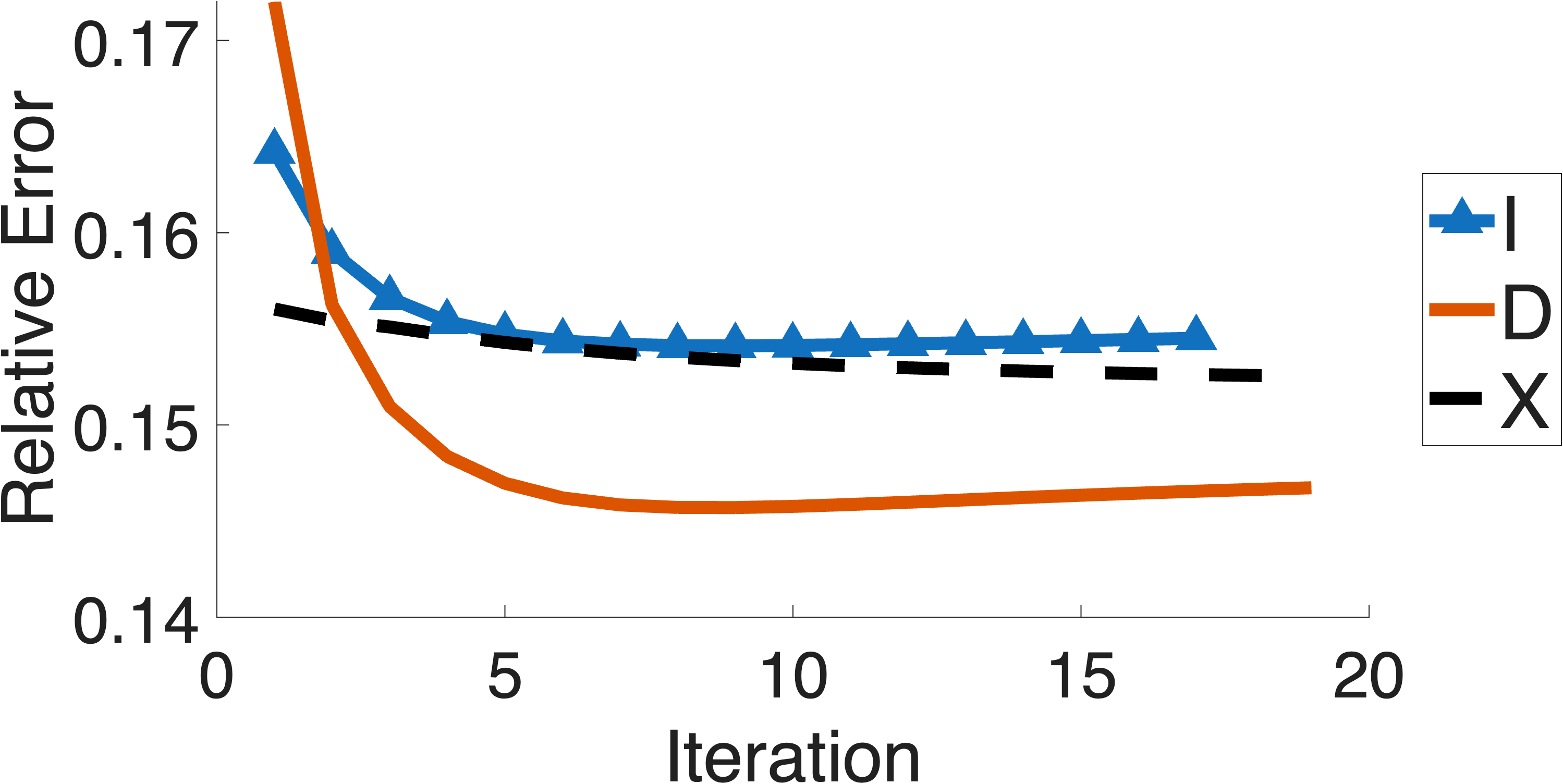}
\caption{HWT, $\sigma=4.5$}
\end{subfigure}\hfill
\begin{subfigure}[t]{0.45\textwidth}
\centering
\includegraphics[width=\linewidth]{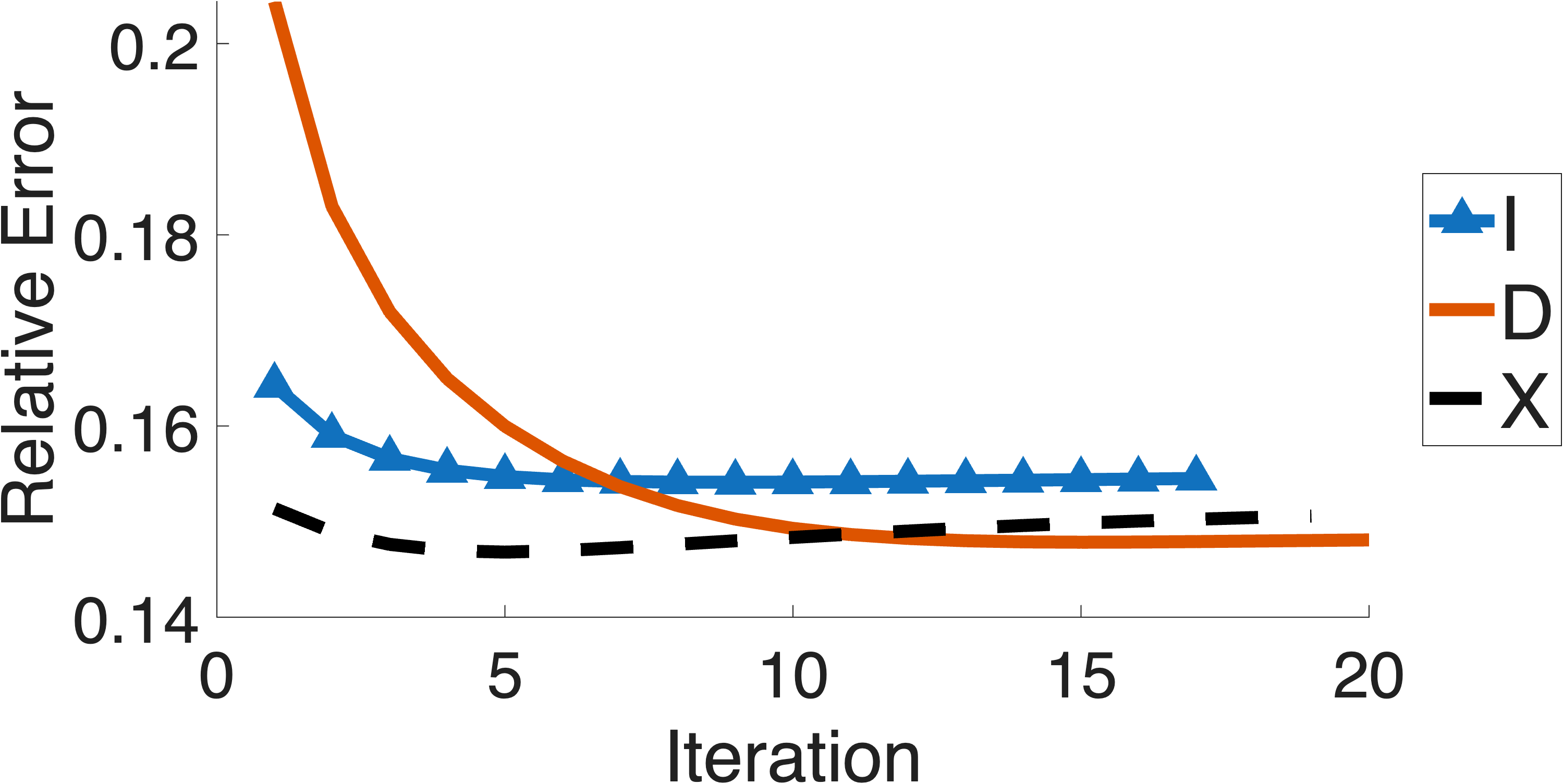}
\caption{DWT, $\sigma=4.5$}
\end{subfigure}
\caption{RRE convergence curves for MM (two-level).}
\label{fig:RREMM}
\end{figure}

\begin{figure}[tbp]
\centering
\begin{subfigure}[t]{0.45\textwidth}
\centering
\includegraphics[width=\linewidth]{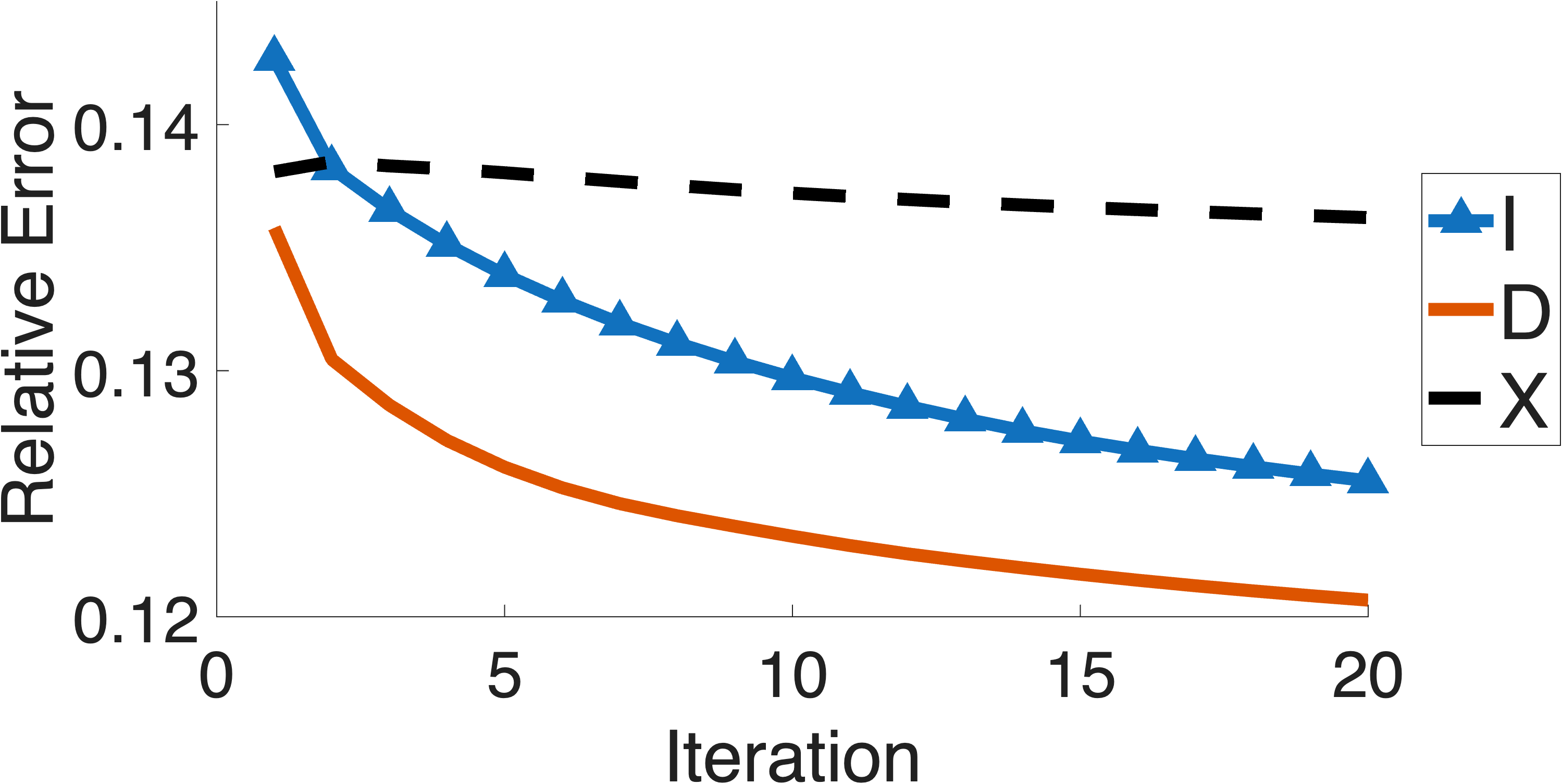}
\caption{HWT, $\sigma=3.5$}
\end{subfigure}\hfill
\begin{subfigure}[t]{0.45\textwidth}
\centering
\includegraphics[width=\linewidth]{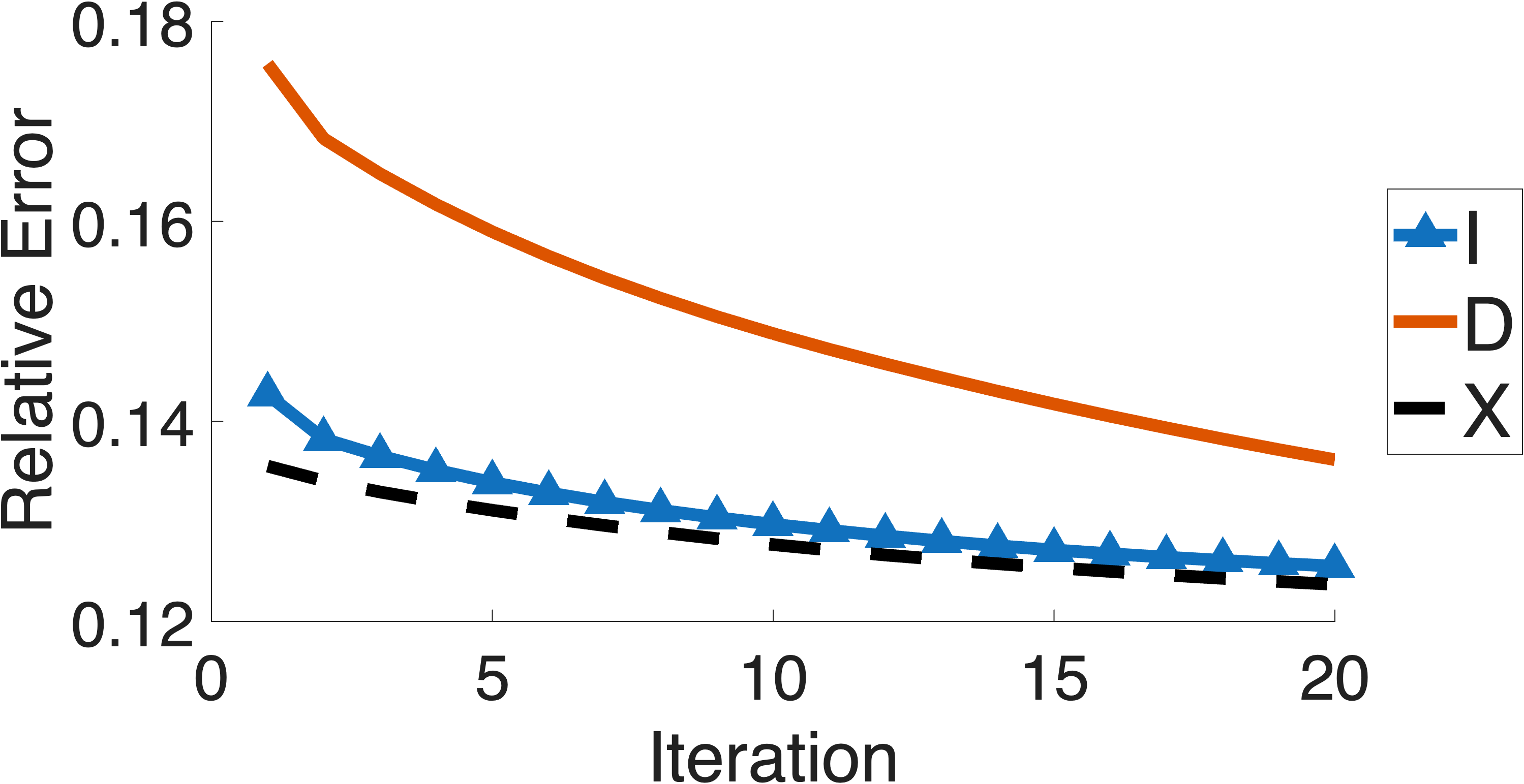}
\caption{DWT, $\sigma=3.5$}
\end{subfigure}
\vspace{0.8em}
\begin{subfigure}[t]{0.45\textwidth}
\centering
\includegraphics[width=\linewidth]{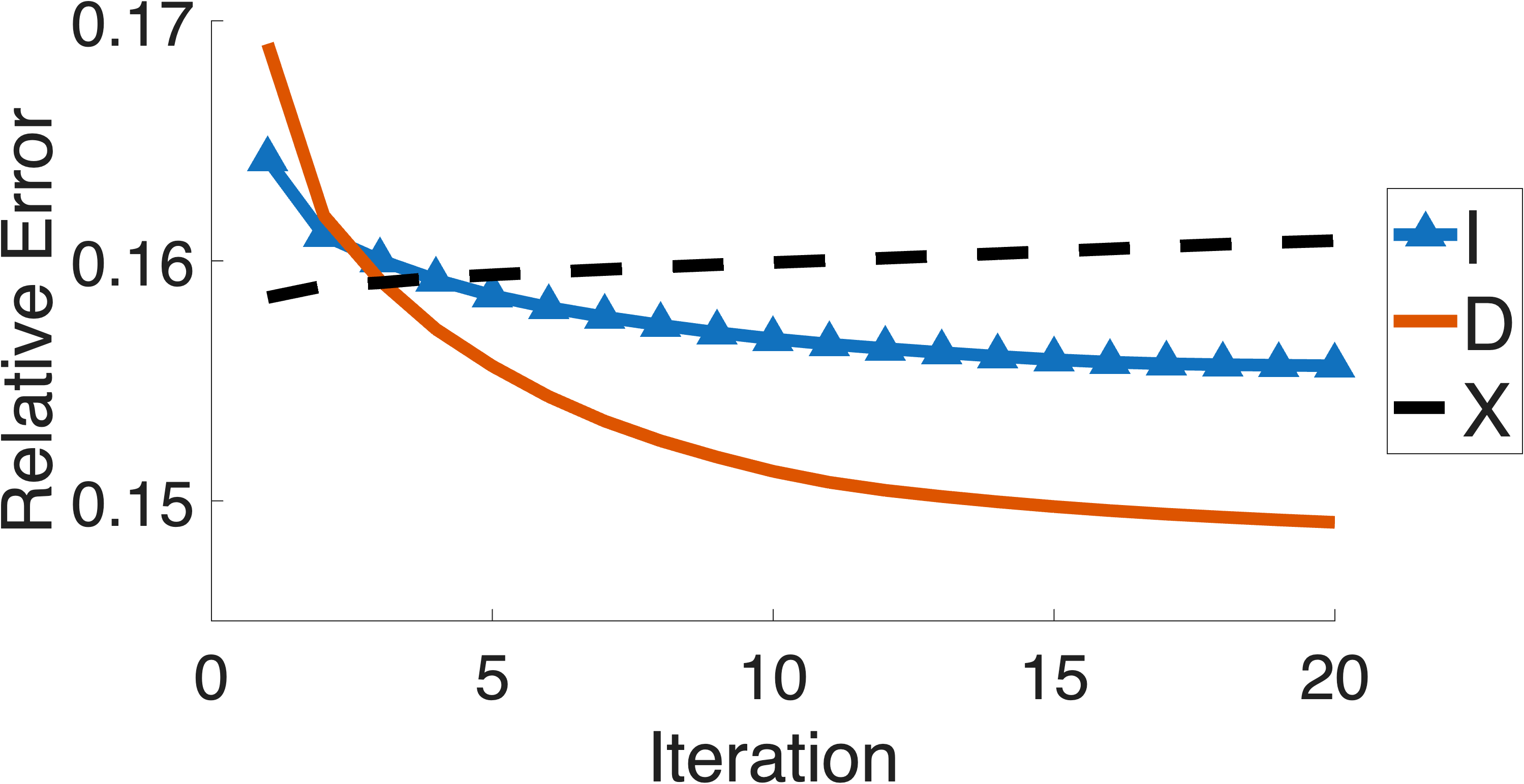}
\caption{HWT, $\sigma=4.5$}
\end{subfigure}\hfill
\begin{subfigure}[t]{0.45\textwidth}
\centering
\includegraphics[width=\linewidth]{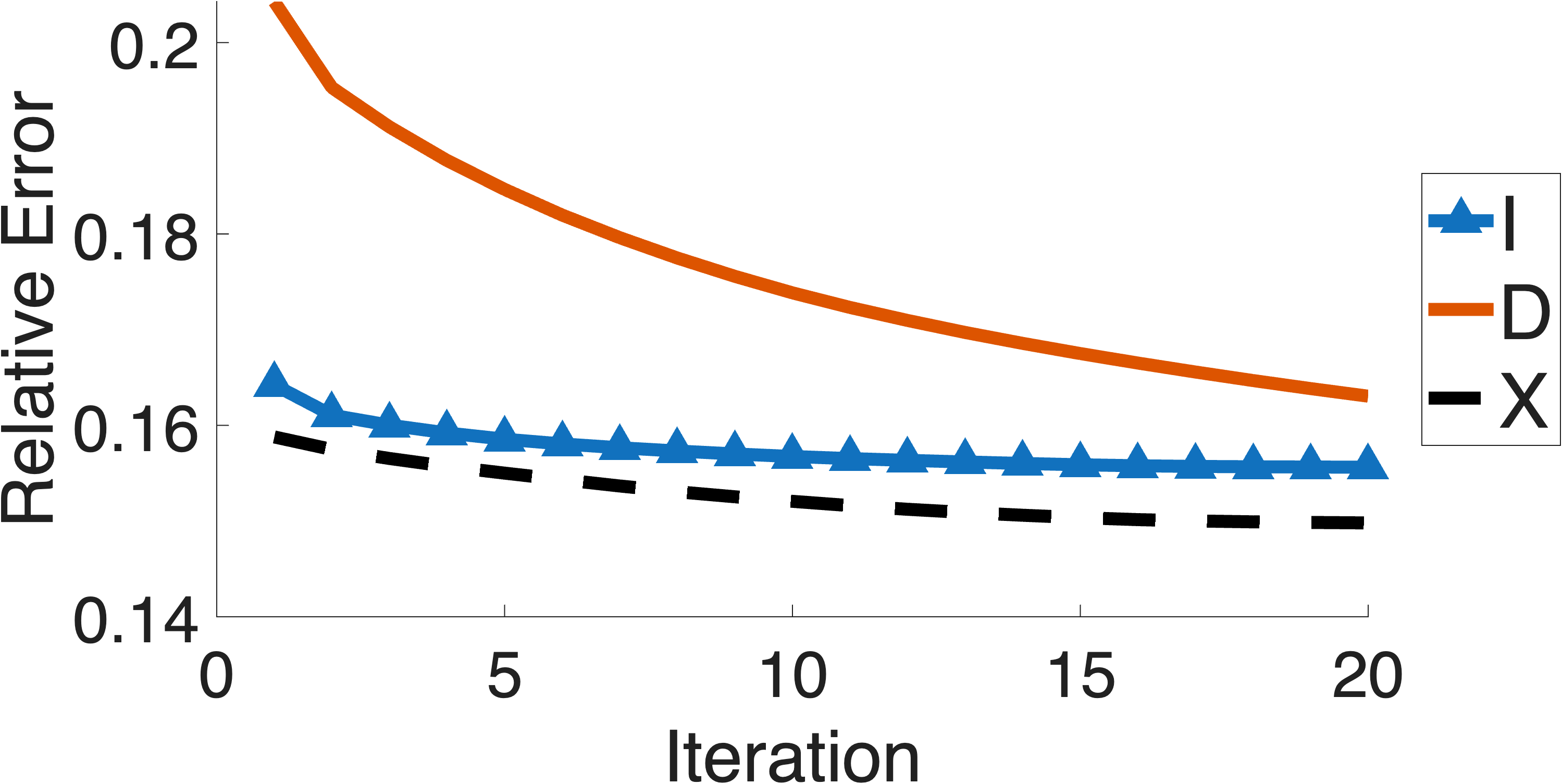}
\caption{DWT, $\sigma=4.5$}
\end{subfigure}
\caption{RRE convergence curves for SB (two-level).}
\label{fig:RRENirnSB}
\end{figure}

Figures~\ref{fig:SolutionTV}-\ref{fig:SolutionSB} display the reconstructed images produced with IRLS, MM, and SB under all five combinations (I, HWT-D, HWT-X, DWT-D, DWT-X), and Figures~\ref{fig:RRENirntv}-\ref{fig:RRENirnSB} show the corresponding RRE curves.

Across all three solvers, the multilevel approaches (D and X) generally outperform Approach~I in reconstruction quality and converge in fewer fine-level outer iterations. Approach~I requires more than 10 outer iterations before converging, whereas the multilevel variants stabilize within 3-5. This reduction indicates that the coarse-level solution or auxiliary quantities provide a fine-level starting point already close to convergence, so fewer GCV-driven refinements are needed at full resolution. Table~\ref{tab:time_all_solvers} shows that IRLS with multilevel achieves speedups exceeding $15\times$ over Approach~I for both blur levels, while for MM and SB the single-level times are already small, so the relative computational savings are more modest and are not uniform across all cases. Critically, the four multilevel combinations have broadly comparable runtimes for each solver, confirming that the choice among them is driven by reconstruction quality rather than computational cost.

Regarding the central finding, under Approach~D, HWT generally yields lower RRE than DWT across the three solvers, while under Approach~X, DWT produces lower RRE in most cases. The same qualitative pattern is observed across all three solvers and both blur levels: HWT tends to outperform DWT under Approach~D, whereas DWT tends to outperform HWT under Approach~X.

\begin{table}[htbp]
\caption{Reconstruction time (seconds) for the two-level method.}
\label{tab:time_all_solvers}
\centering
\begin{tabular}{llccccc}
\toprule
Solver & $\sigma$ & I & HWT-D & HWT-X & DWT-D & DWT-X \\
\midrule
IRLS & 3.5 & 243 & 8 & 7 & 8 & 8 \\
IRLS & 4.5 & 240 & 8 & 7 & 9 & 8 \\
\midrule
MM & 3.5 & 2 & 3 & 2 & 3 & 3 \\
MM & 4.5 & 2 & 3 & 3 & 5 & 4 \\
\midrule
SB & 3.5 & 3 & 4 & 4 & 4 & 4 \\
SB & 4.5 & 4 & 4 & 4 & 7 & 5 \\
\bottomrule
\end{tabular}
\end{table}

\subsection{Wavelet Comparison}

\begin{figure}[tbp]
\centering
\begin{subfigure}[t]{0.24\textwidth}
\centering
\includegraphics[width=\linewidth]{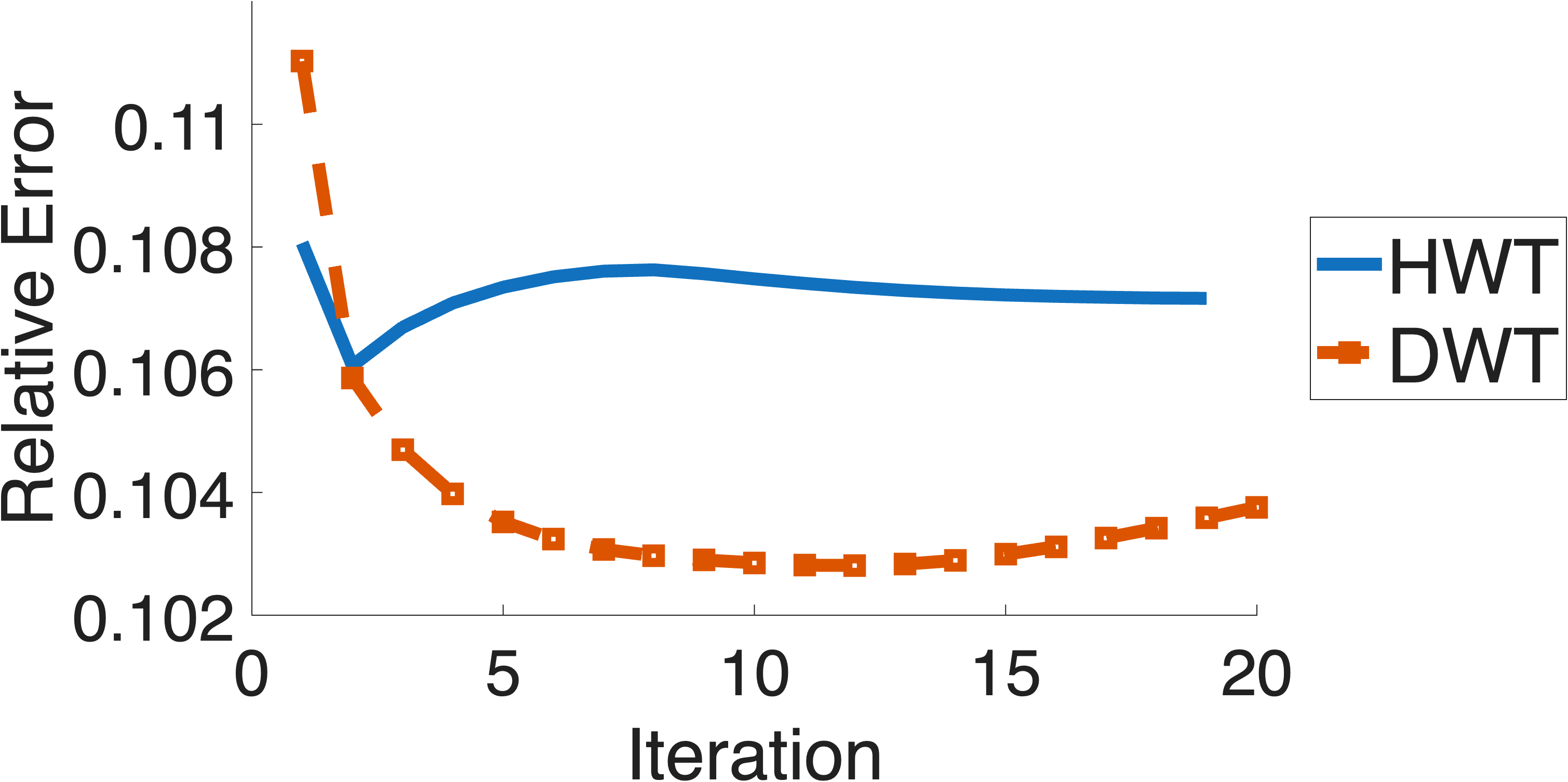}
\caption{IRLS, D, $\sigma=3.5$}
\end{subfigure}\hfill
\begin{subfigure}[t]{0.24\textwidth}
\centering
\includegraphics[width=\linewidth]{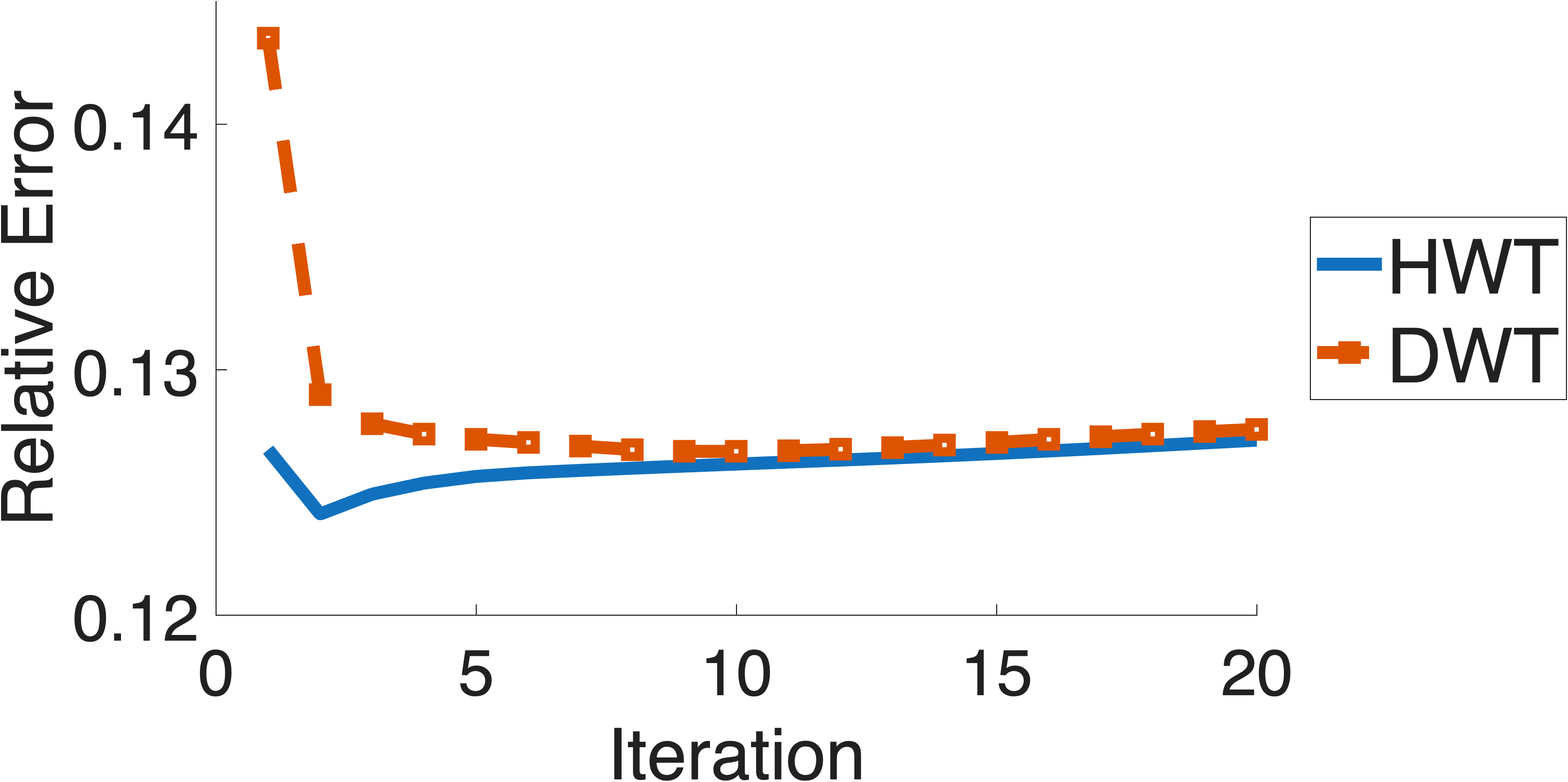}
\caption{IRLS, D, $\sigma=4.5$}
\end{subfigure}\hfill
\begin{subfigure}[t]{0.24\textwidth}
\centering
\includegraphics[width=\linewidth]{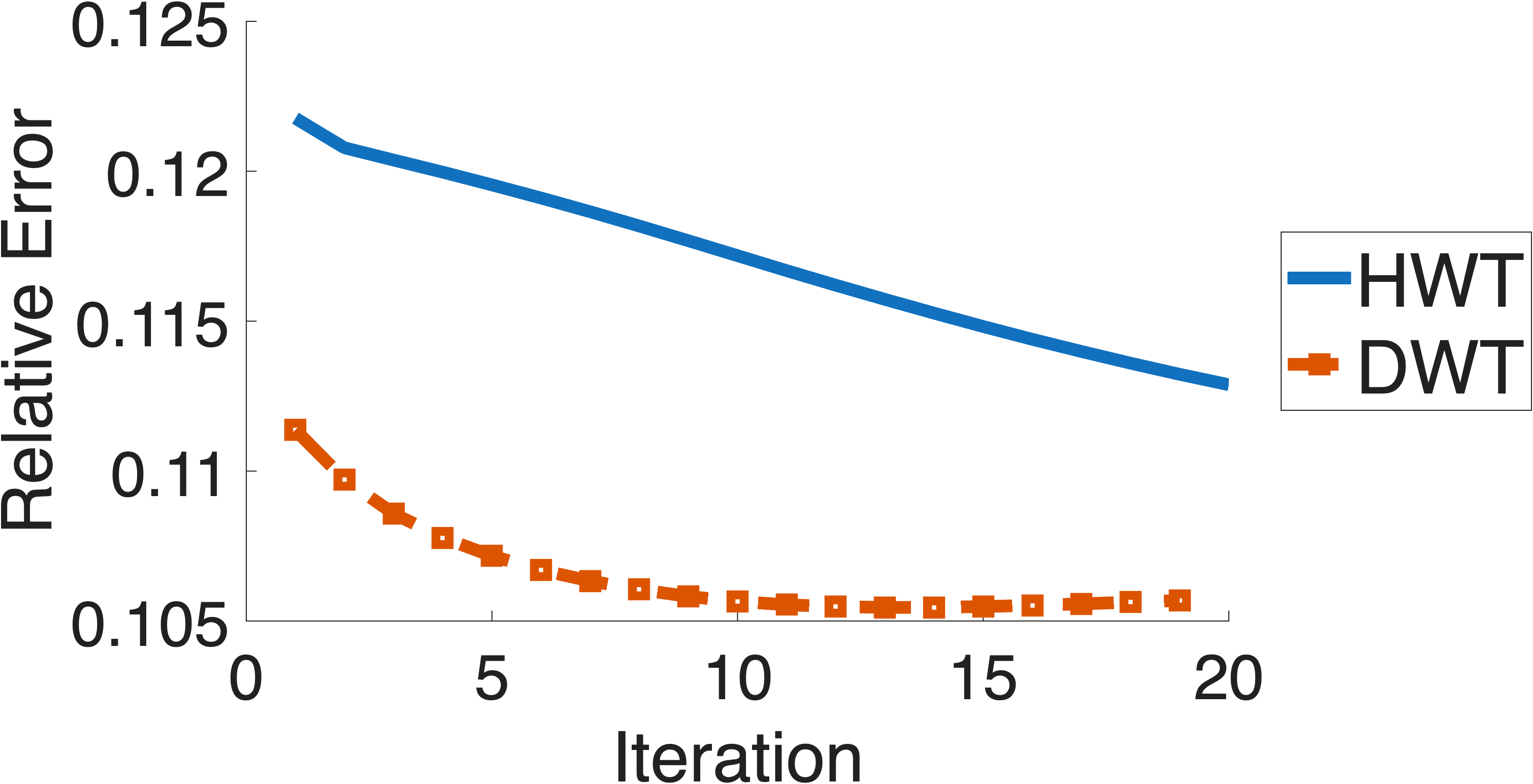}
\caption{IRLS, X, $\sigma=3.5$}
\end{subfigure}
\vspace{0.6em}
\begin{subfigure}[t]{0.24\textwidth}
\centering
\includegraphics[width=\linewidth]{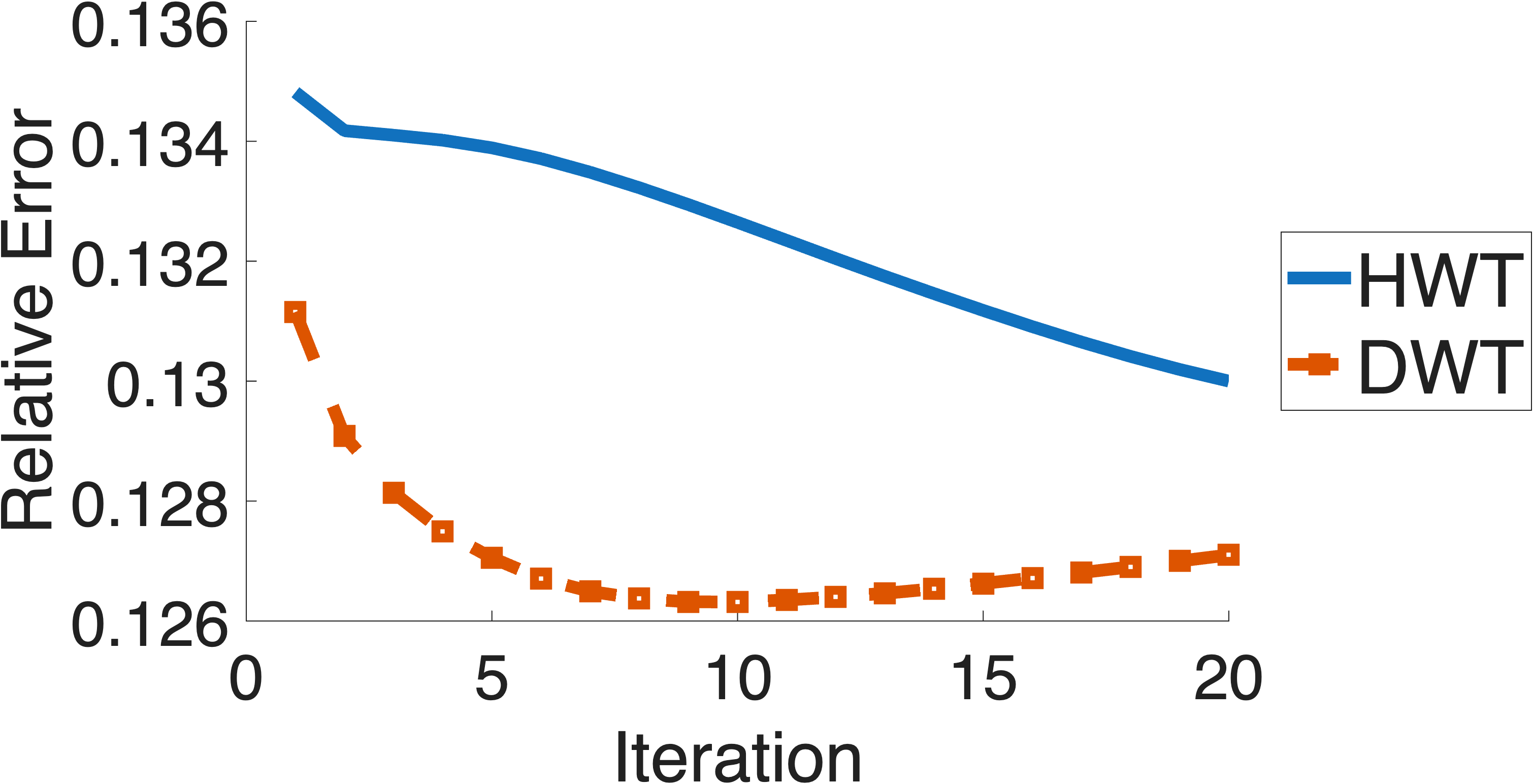}
\caption{IRLS, X, $\sigma=4.5$}
\end{subfigure}\hfill
\begin{subfigure}[t]{0.24\textwidth}
\centering
\includegraphics[width=\linewidth]{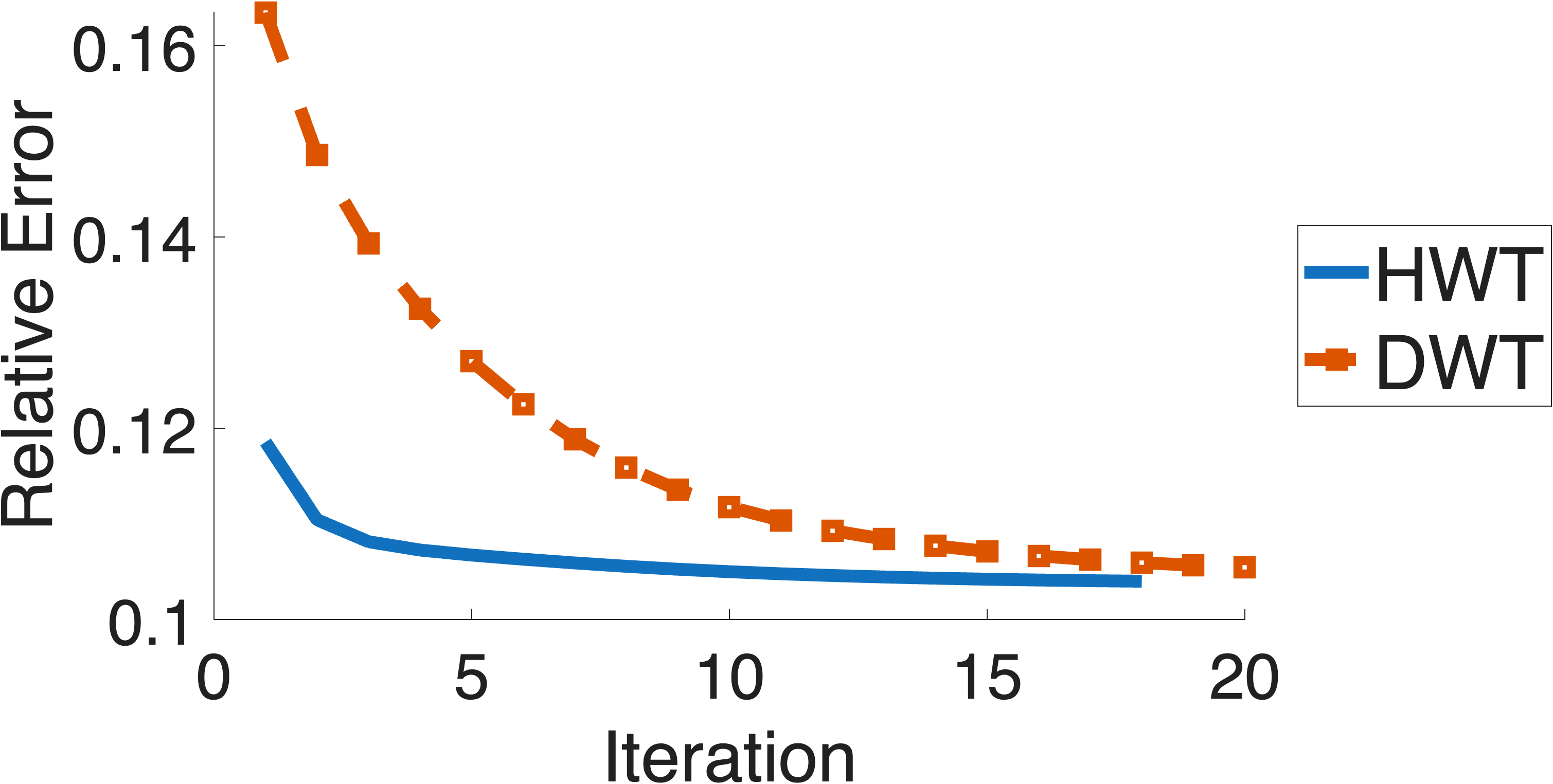}
\caption{MM, D, $\sigma=3.5$}
\end{subfigure}\hfill
\begin{subfigure}[t]{0.24\textwidth}
\centering
\includegraphics[width=\linewidth]{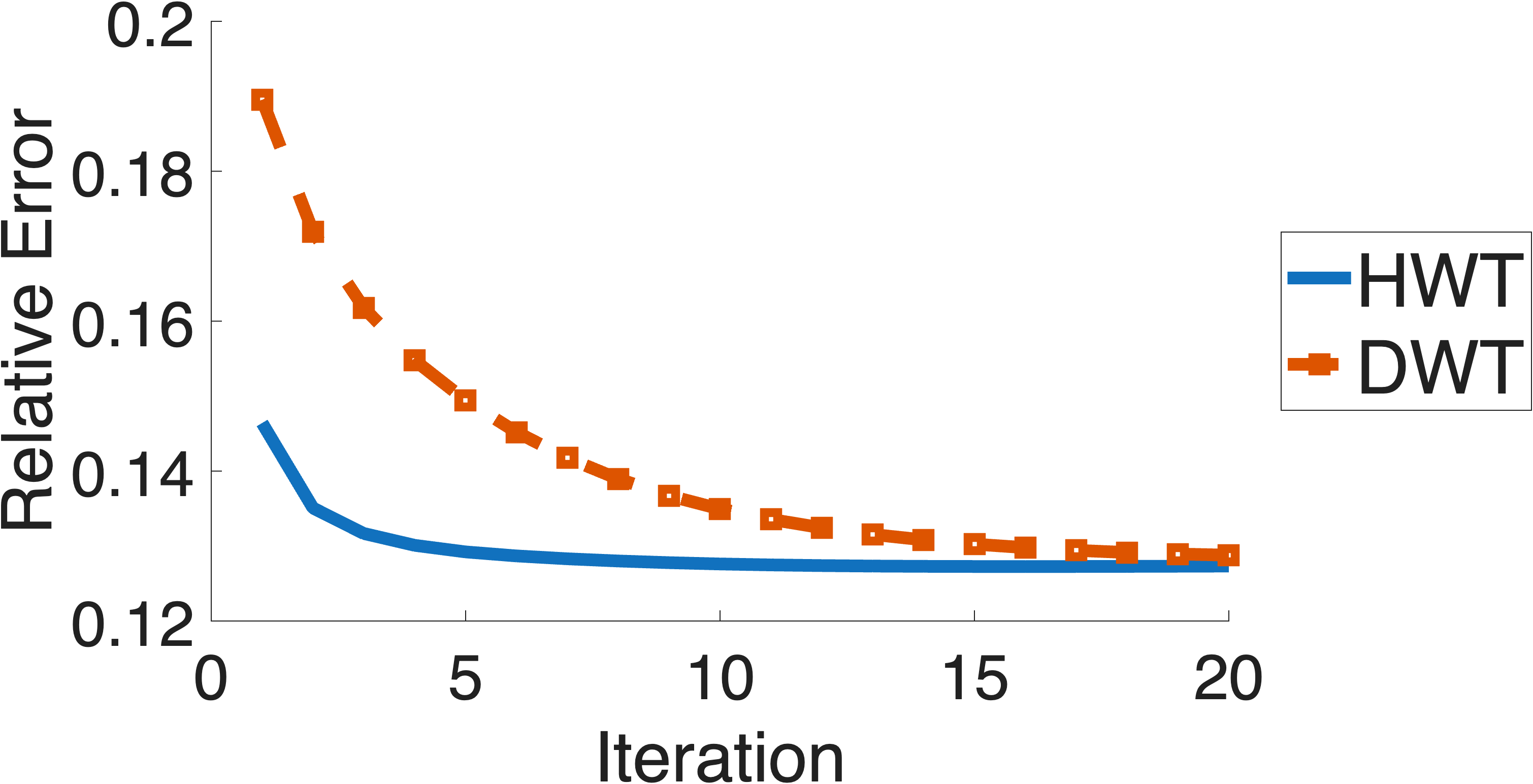}
\caption{MM, D, $\sigma=4.5$}
\end{subfigure}
\vspace{0.6em}
\begin{subfigure}[t]{0.24\textwidth}
\centering
\includegraphics[width=\linewidth]{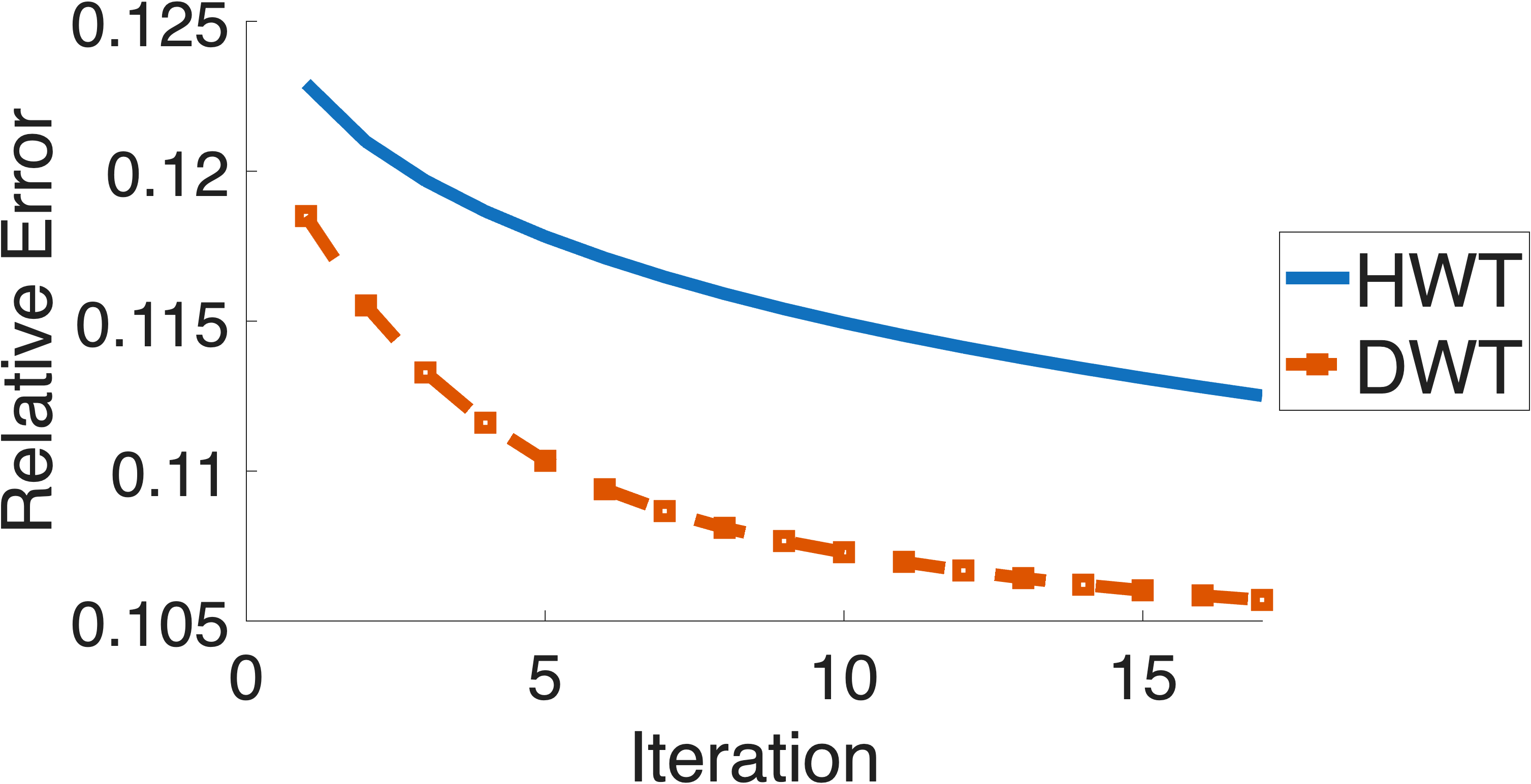}
\caption{MM, X, $\sigma=3.5$}
\end{subfigure}\hfill
\begin{subfigure}[t]{0.24\textwidth}
\centering
\includegraphics[width=\linewidth]{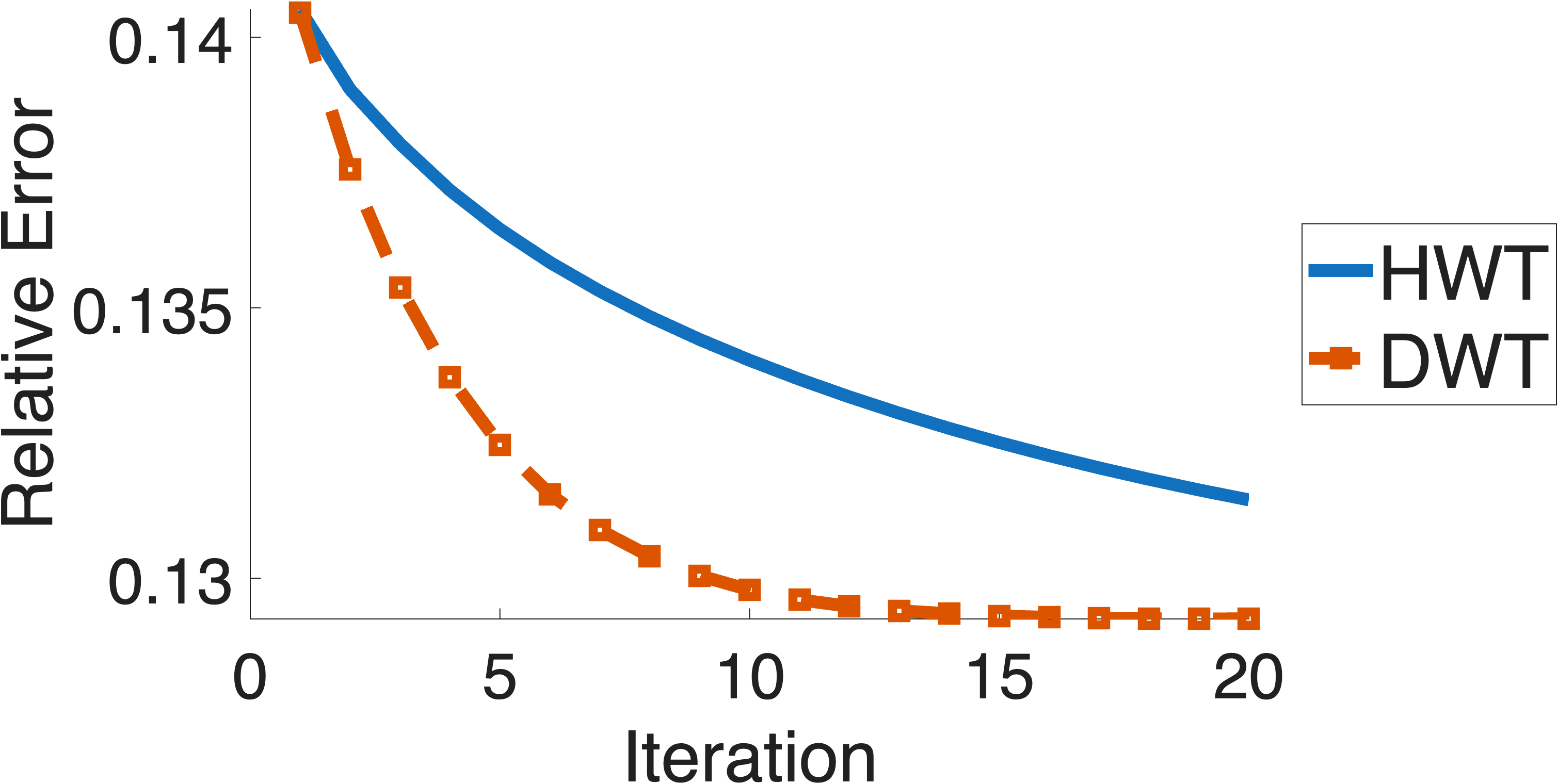}
\caption{MM, X, $\sigma=4.5$}
\end{subfigure}\hfill
\begin{subfigure}[t]{0.24\textwidth}
\centering
\includegraphics[width=\linewidth]{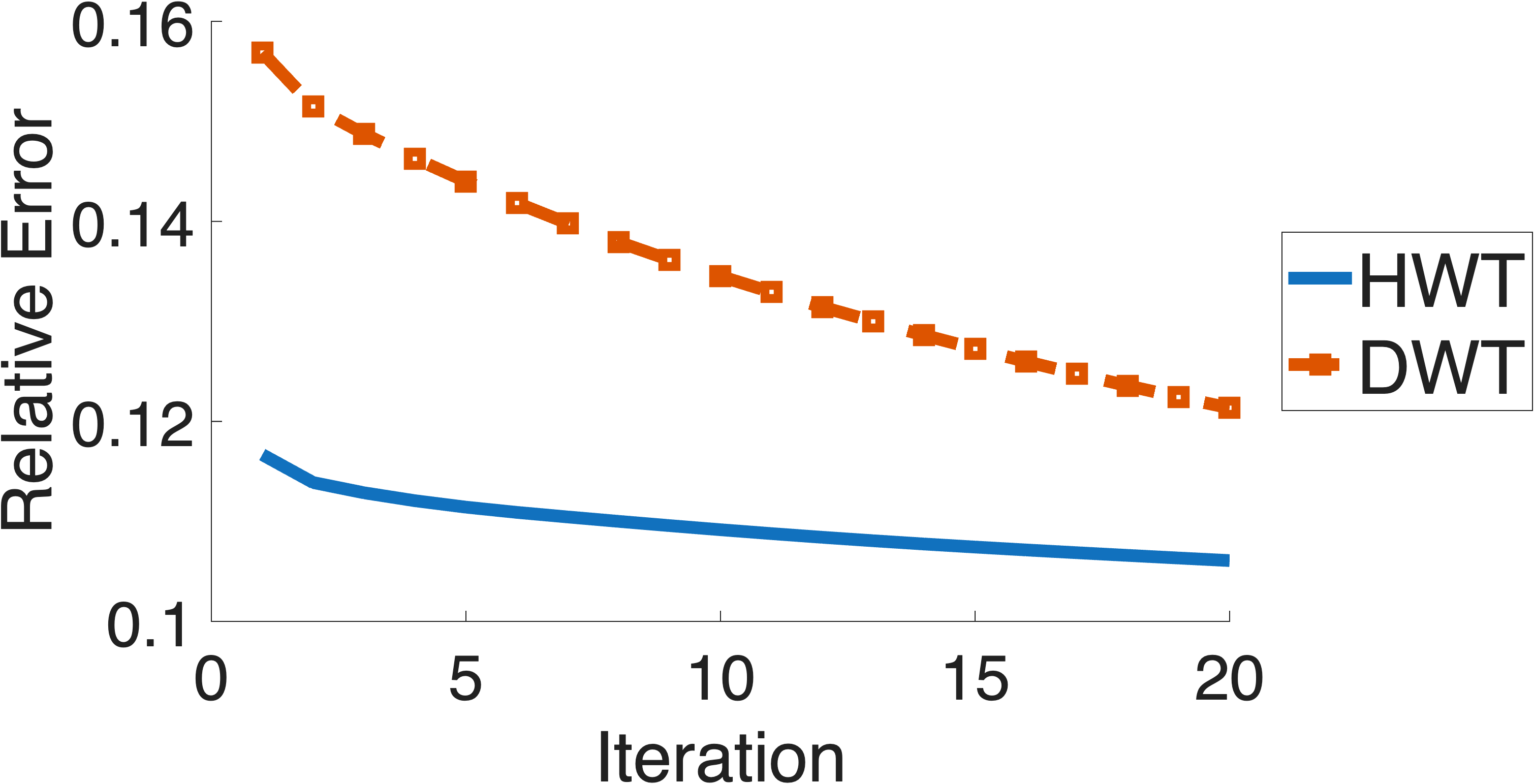}
\caption{SB, D, $\sigma=3.5$}
\end{subfigure}
\vspace{0.6em}
\begin{subfigure}[t]{0.24\textwidth}
\centering
\includegraphics[width=\linewidth]{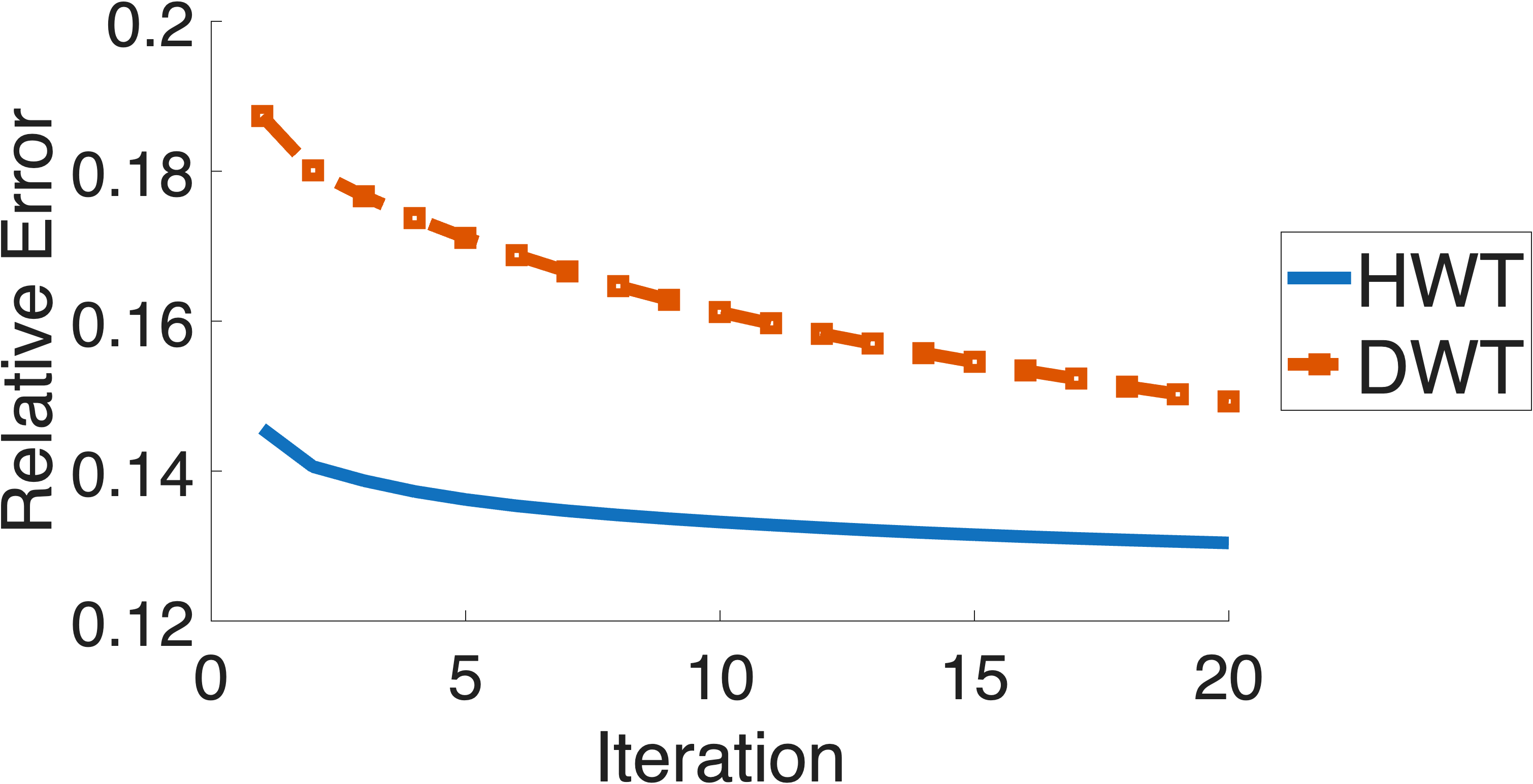}
\caption{SB, D, $\sigma=4.5$}
\end{subfigure}\hfill
\begin{subfigure}[t]{0.24\textwidth}
\centering
\includegraphics[width=\linewidth]{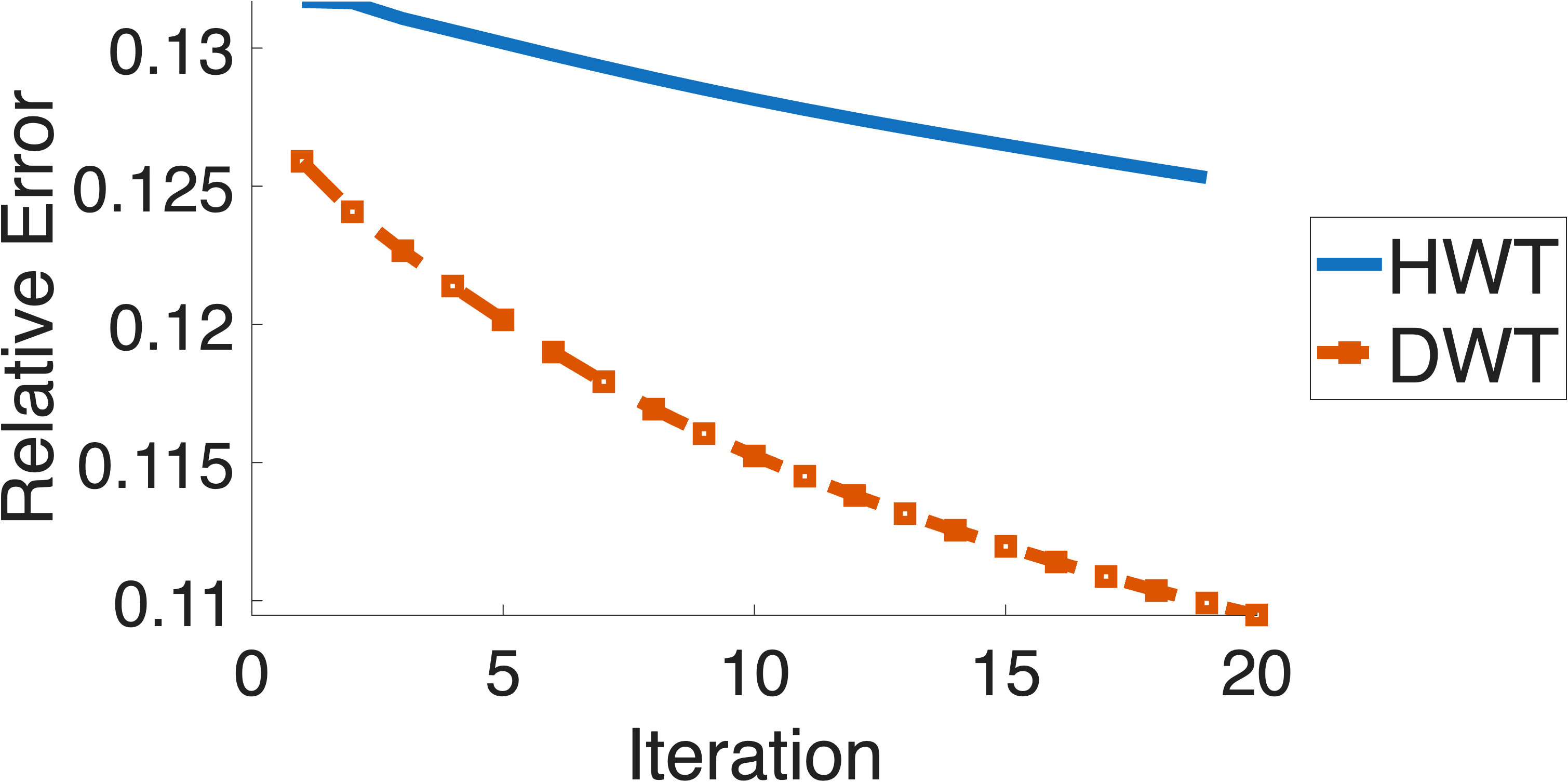}
\caption{SB, X, $\sigma=3.5$}
\end{subfigure}\hfill
\begin{subfigure}[t]{0.24\textwidth}
\centering
\includegraphics[width=\linewidth]{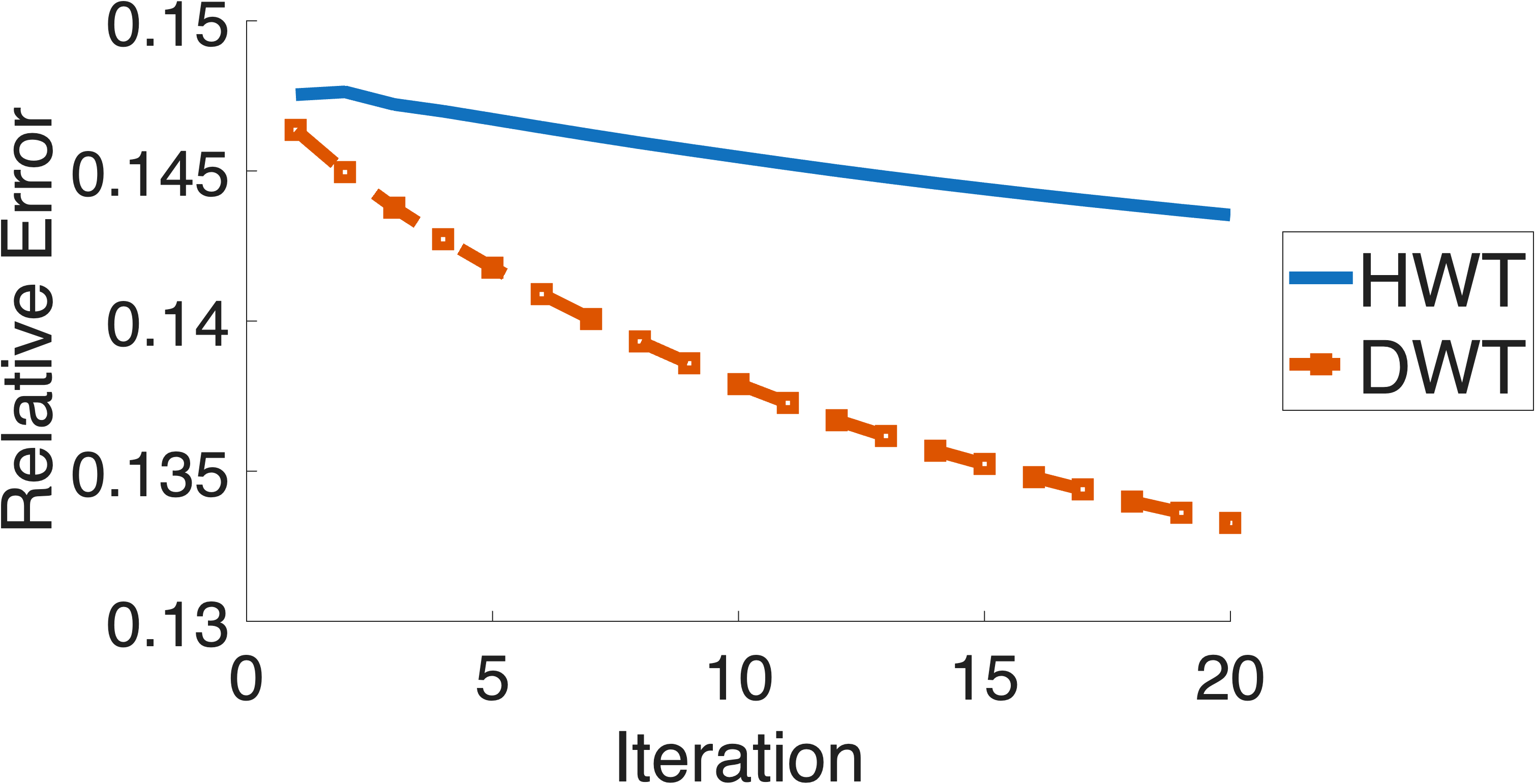}
\caption{SB, X, $\sigma=4.5$}
\end{subfigure}
\caption{RRE convergence curves comparing HWT and DWT (two-level).}
\label{fig:RREHWTvsDWT}
\end{figure}

Figure~\ref{fig:RREHWTvsDWT} shows RRE curves comparing HWT and DWT under each approach and solver. A general interaction pattern is observed: HWT-D converges to lower RRE than DWT-D in most cases, while DWT-X converges to lower RRE than HWT-X. This asymmetry suggests that the wavelet family best suited to the multilevel cycle depends on the type of information transferred between levels. Within the experiments considered here, Haar generally performs better when solver-specific auxiliary quantities are transferred, whereas Daubechies generally performs better when only the solution is transferred. In terms of runtime, HWT is marginally faster than DWT (at most $2\times$), reflecting its simpler two-tap filter.

\subsection{Three-Level Results} \label{sec:threelevel}

\begin{figure}[tbp]
\centering
\begin{subfigure}[t]{0.23\textwidth}
\centering
\includegraphics[width=\linewidth]{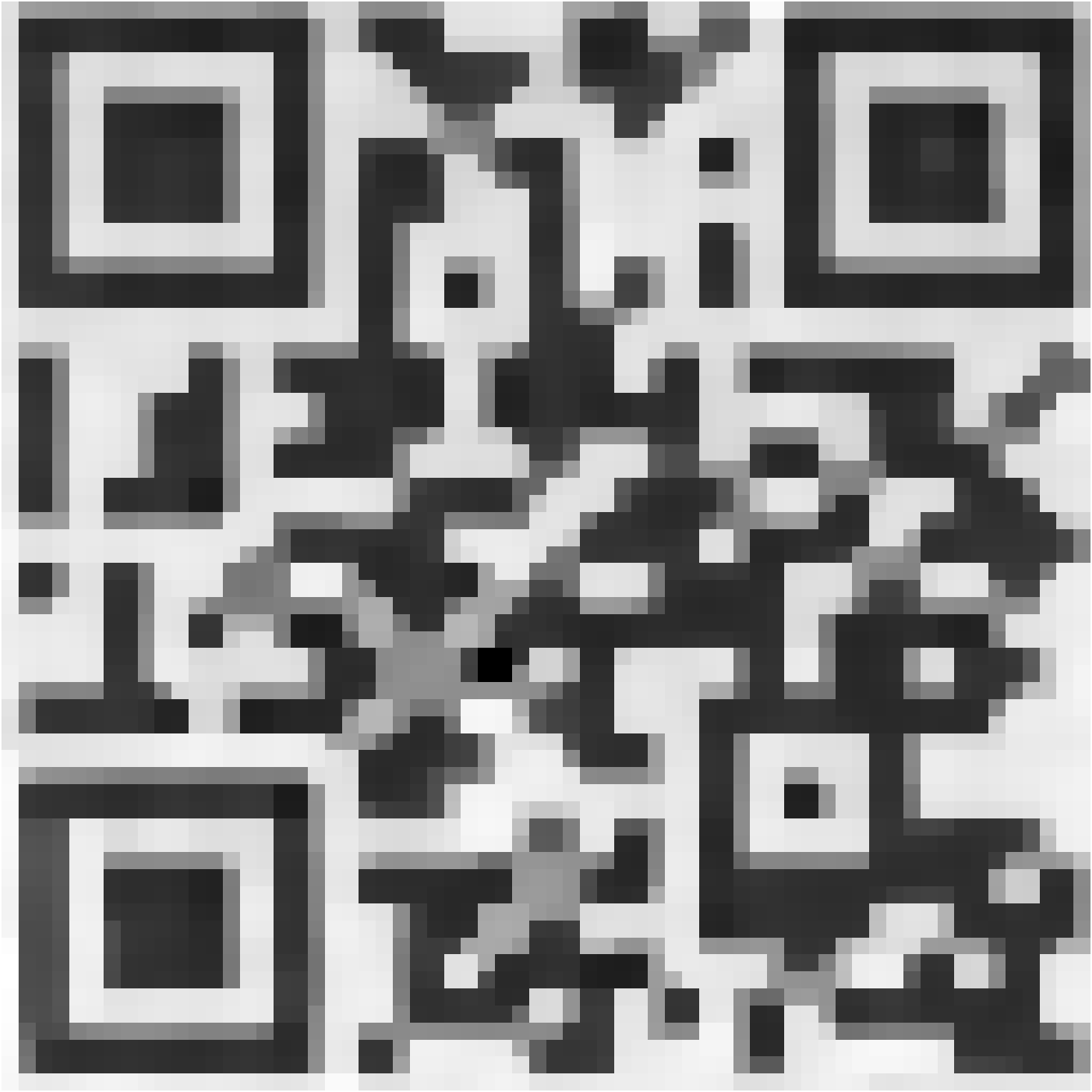}
\caption{HWT-D, $\sigma=3.5$}
\end{subfigure}\hfill
\begin{subfigure}[t]{0.23\textwidth}
\centering
\includegraphics[width=\linewidth]{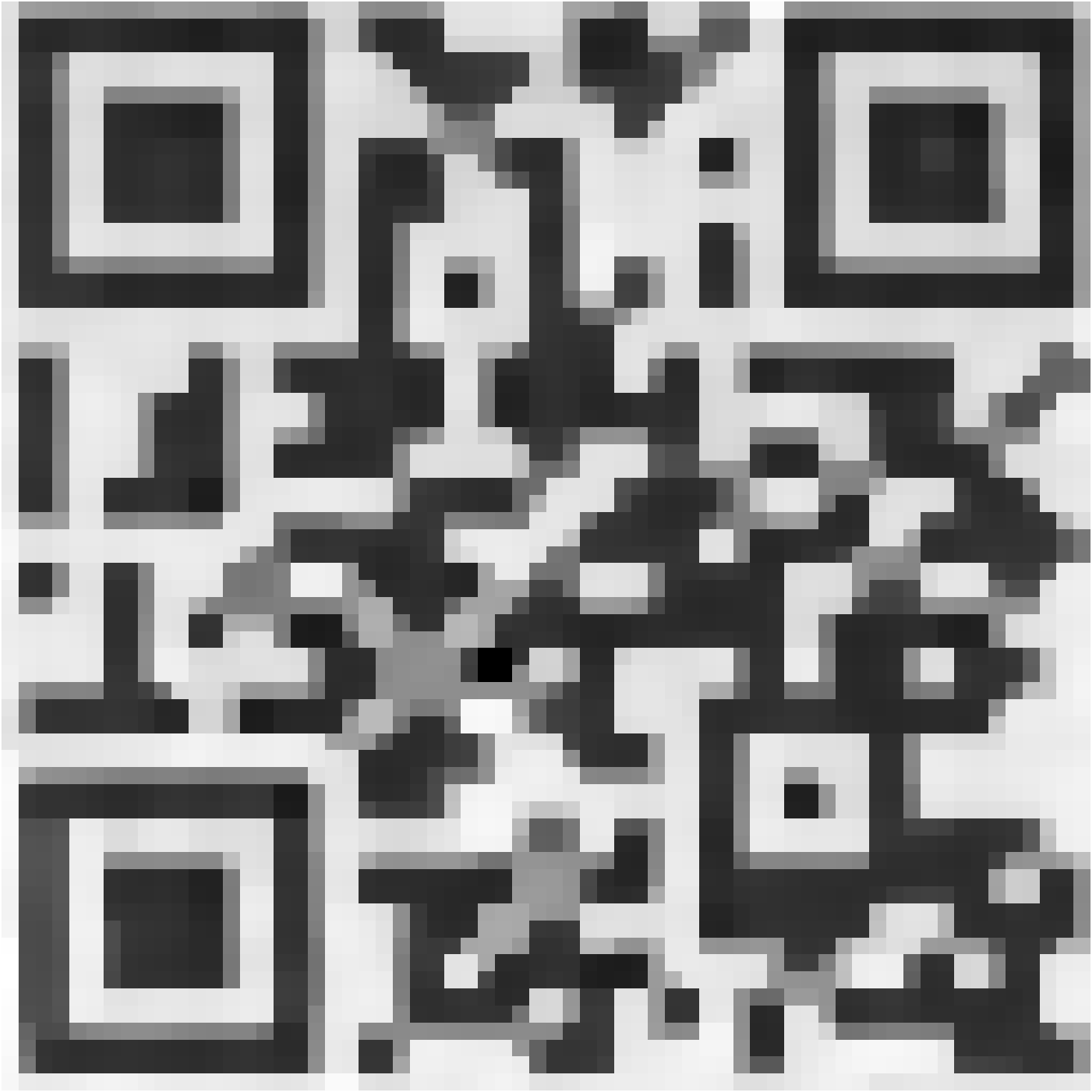}
\caption{HWT-X, $\sigma=3.5$}
\end{subfigure}\hfill
\begin{subfigure}[t]{0.23\textwidth}
\centering
\includegraphics[width=\linewidth]{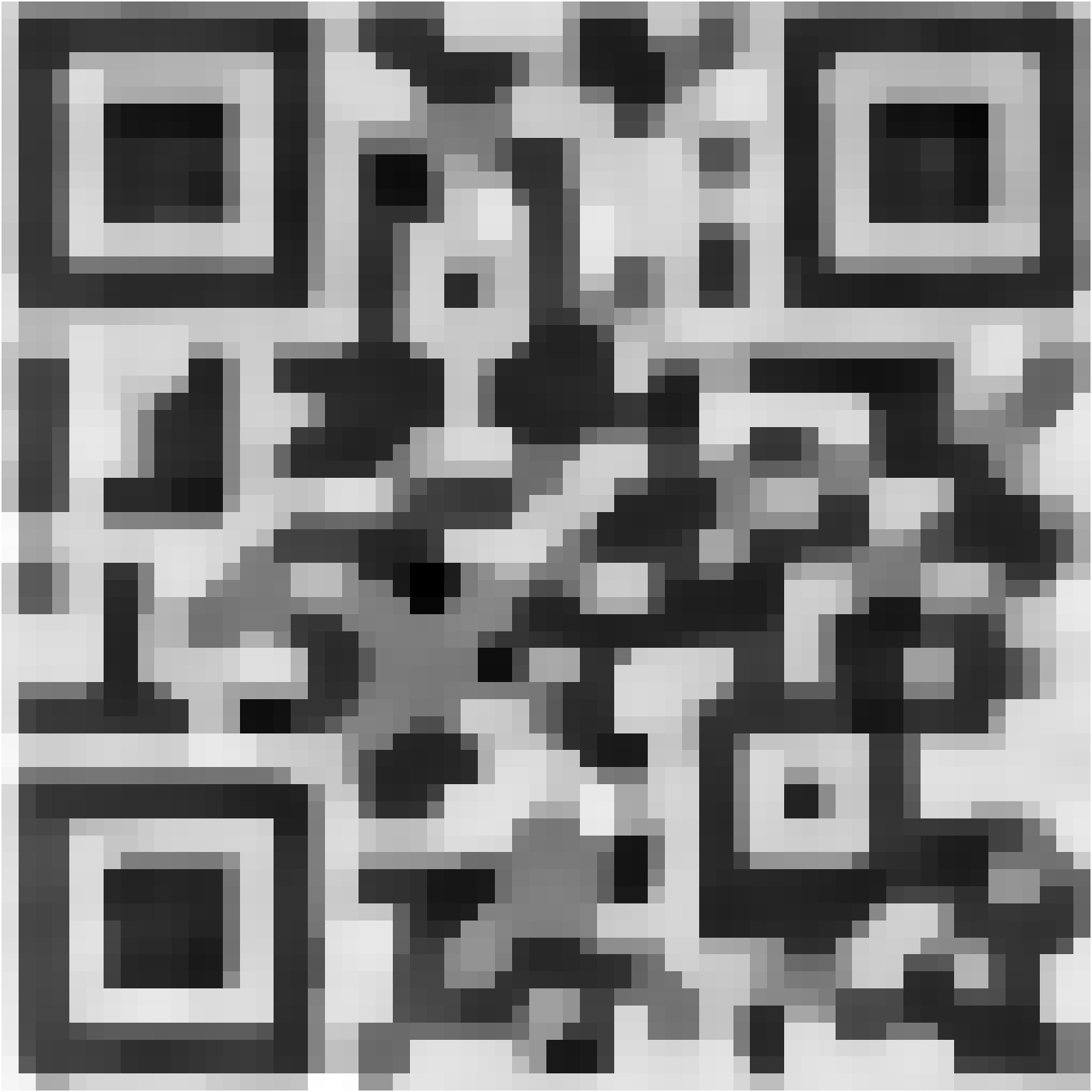}
\caption{HWT-D, $\sigma=4.5$}
\end{subfigure}\hfill
\begin{subfigure}[t]{0.23\textwidth}
\centering
\includegraphics[width=\linewidth]{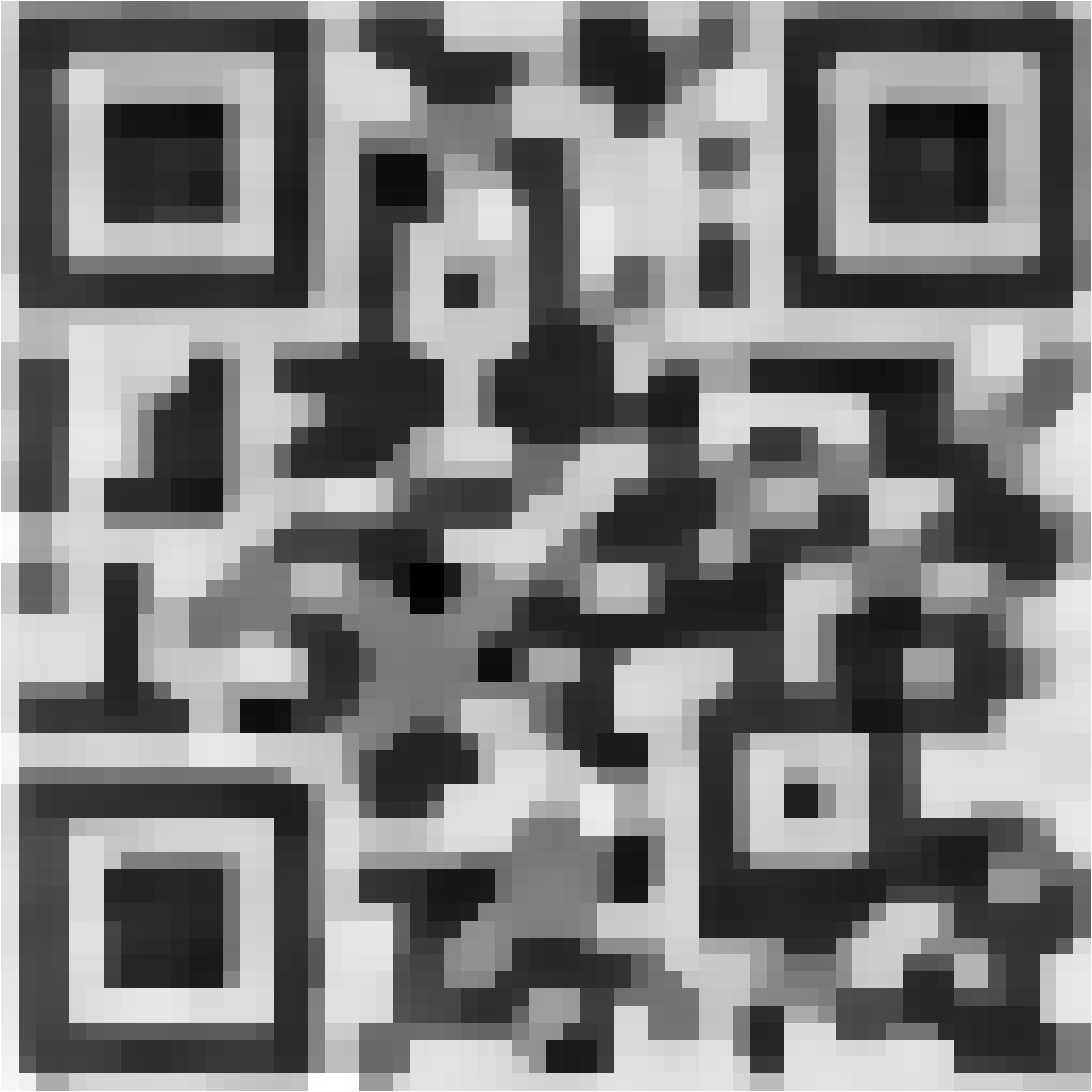}
\caption{HWT-X, $\sigma=4.5$}
\end{subfigure}
\vspace{0.8em}
\begin{subfigure}[t]{0.23\textwidth}
\centering
\includegraphics[width=\linewidth]{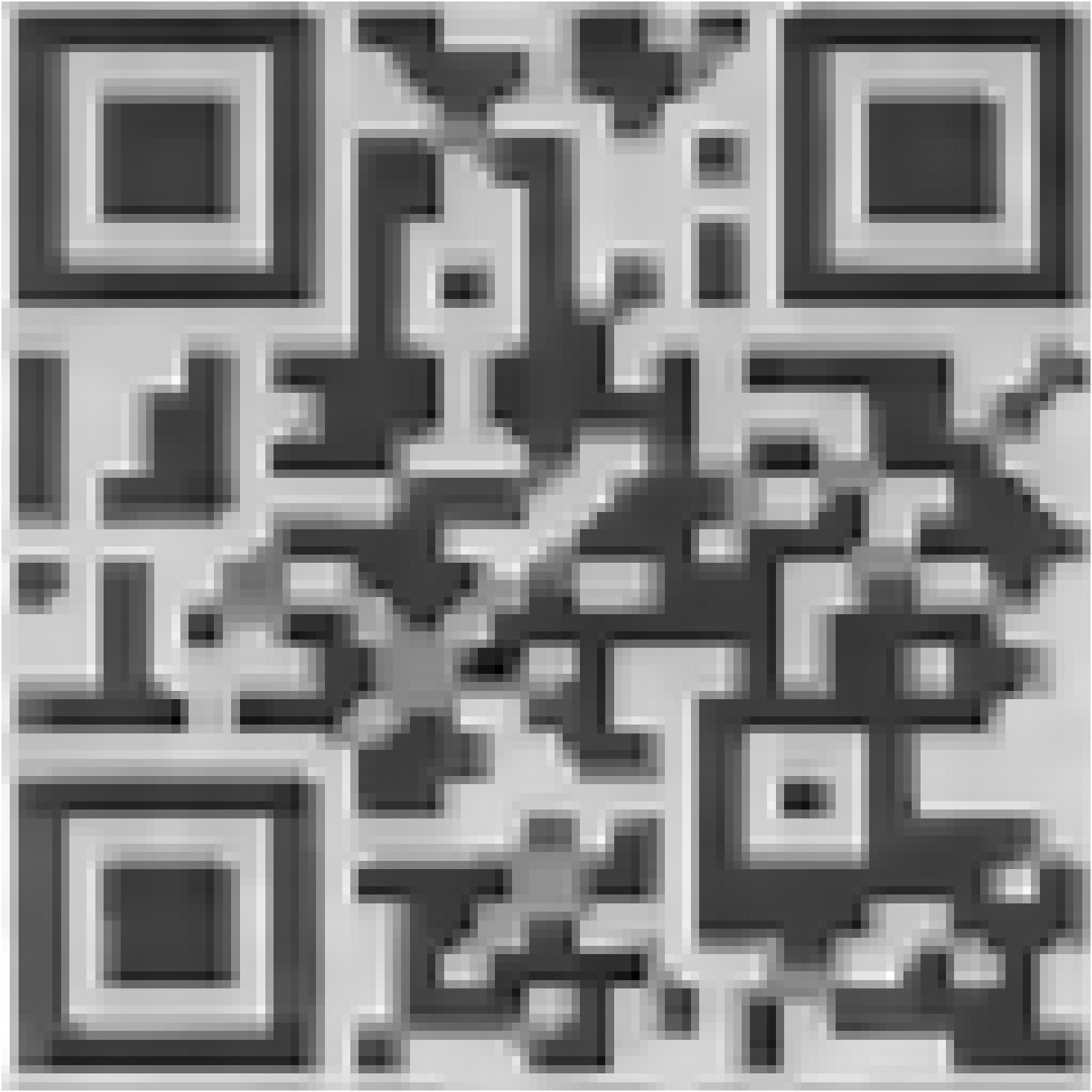}
\caption{DWT-D, $\sigma=3.5$}
\end{subfigure}\hfill
\begin{subfigure}[t]{0.23\textwidth}
\centering
\includegraphics[width=\linewidth]{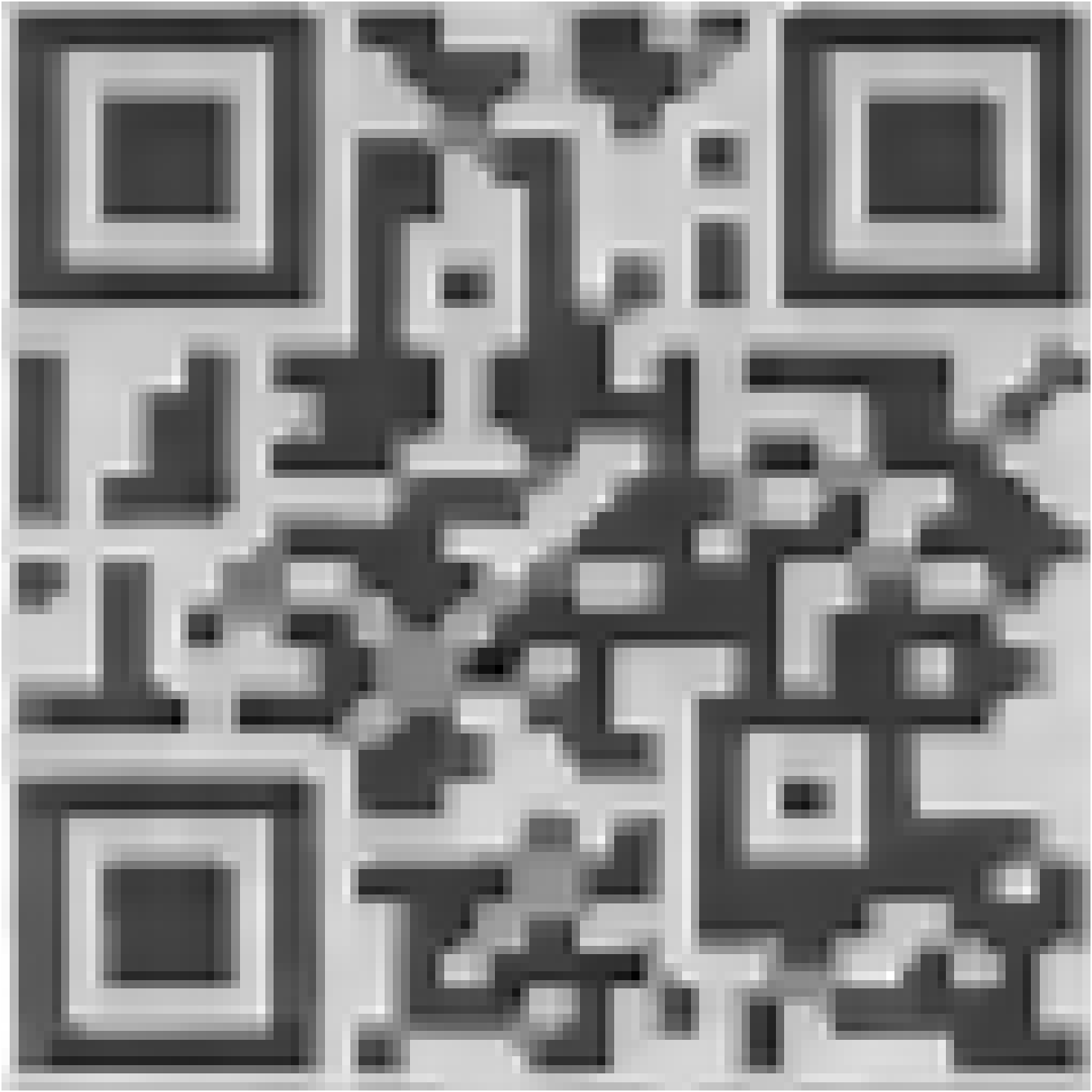}
\caption{DWT-X, $\sigma=3.5$}
\end{subfigure}\hfill
\begin{subfigure}[t]{0.23\textwidth}
\centering
\includegraphics[width=\linewidth]{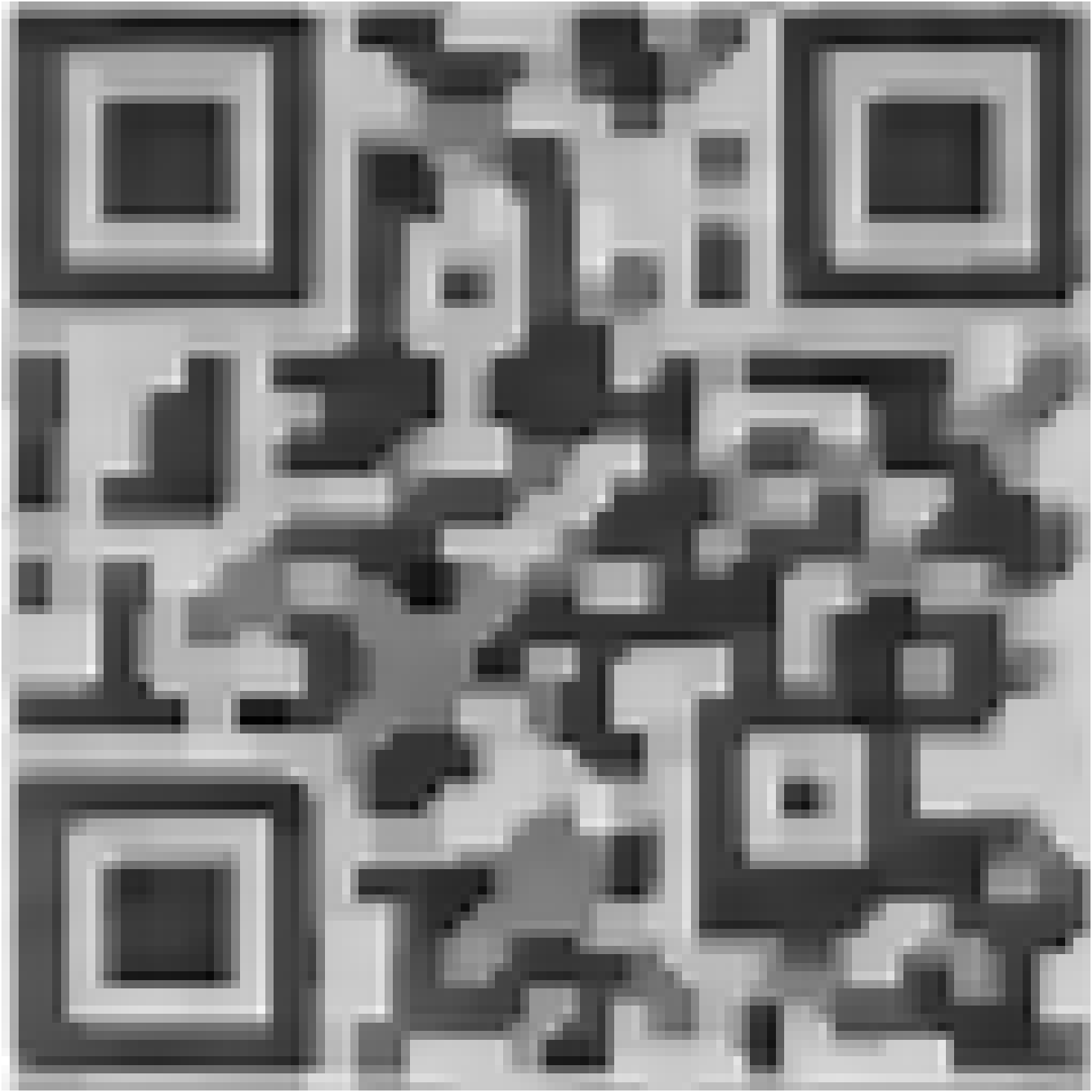}
\caption{DWT-D, $\sigma=4.5$}
\end{subfigure}\hfill
\begin{subfigure}[t]{0.23\textwidth}
\centering
\includegraphics[width=\linewidth]{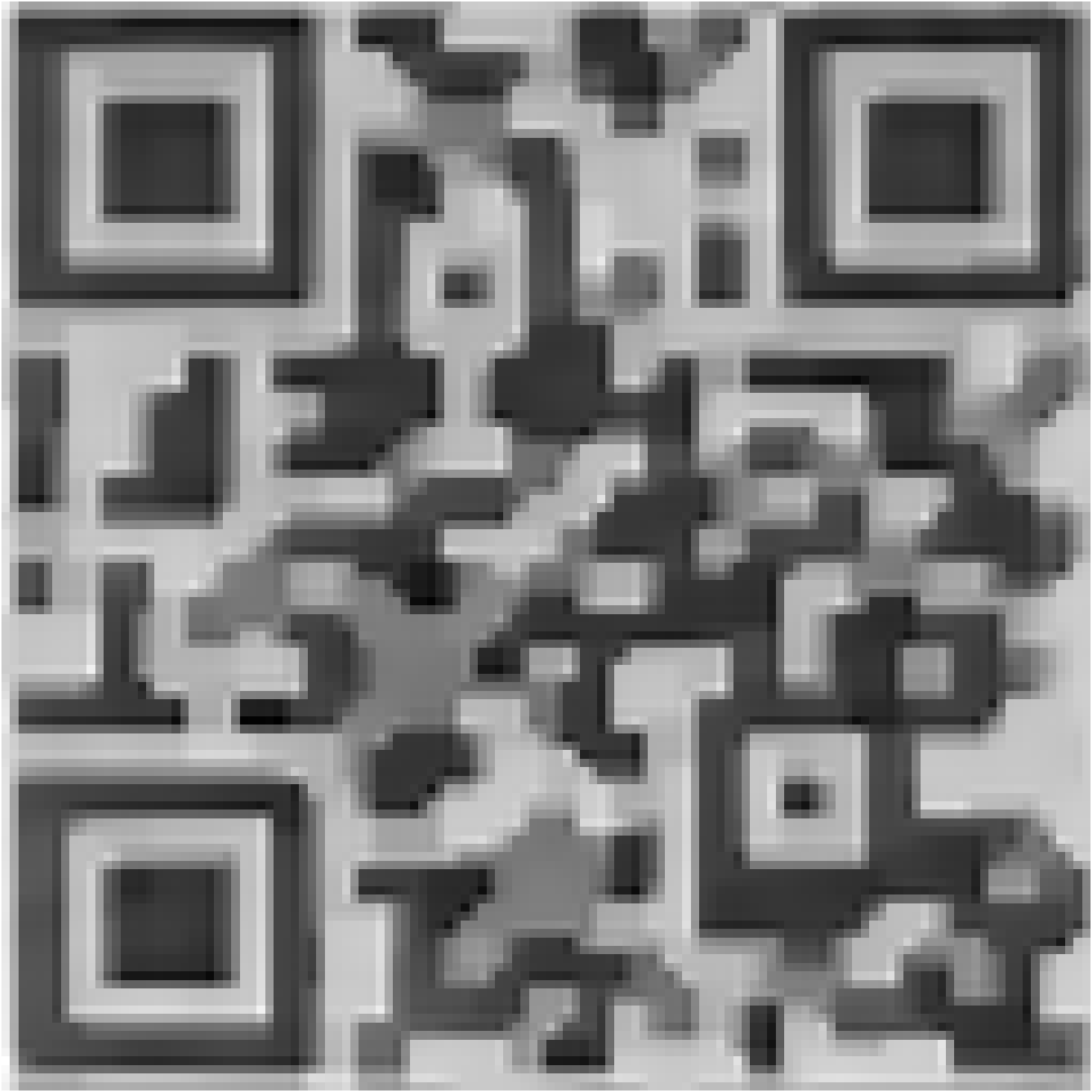}
\caption{DWT-X, $\sigma=4.5$}
\end{subfigure}
\caption{Three-level IRLS reconstructions. Top row: HWT; bottom row: DWT.}
\label{fig:IRLS3levelQ}
\end{figure}

\begin{figure}[tbp]
\centering
\begin{subfigure}[t]{0.23\textwidth}
\centering
\includegraphics[width=\linewidth]{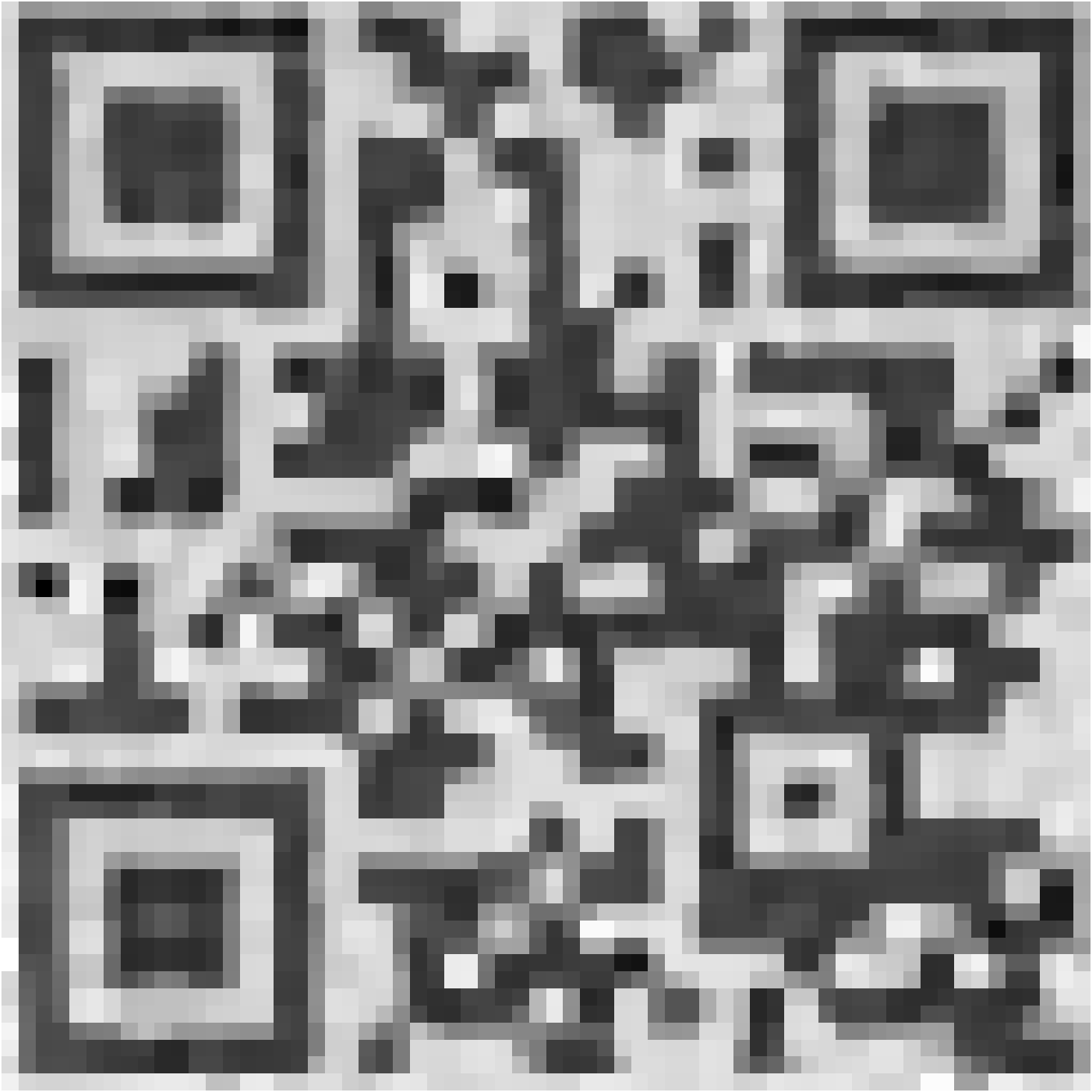}
\caption{HWT-D, $\sigma=3.5$}
\end{subfigure}\hfill
\begin{subfigure}[t]{0.23\textwidth}
\centering
\includegraphics[width=\linewidth]{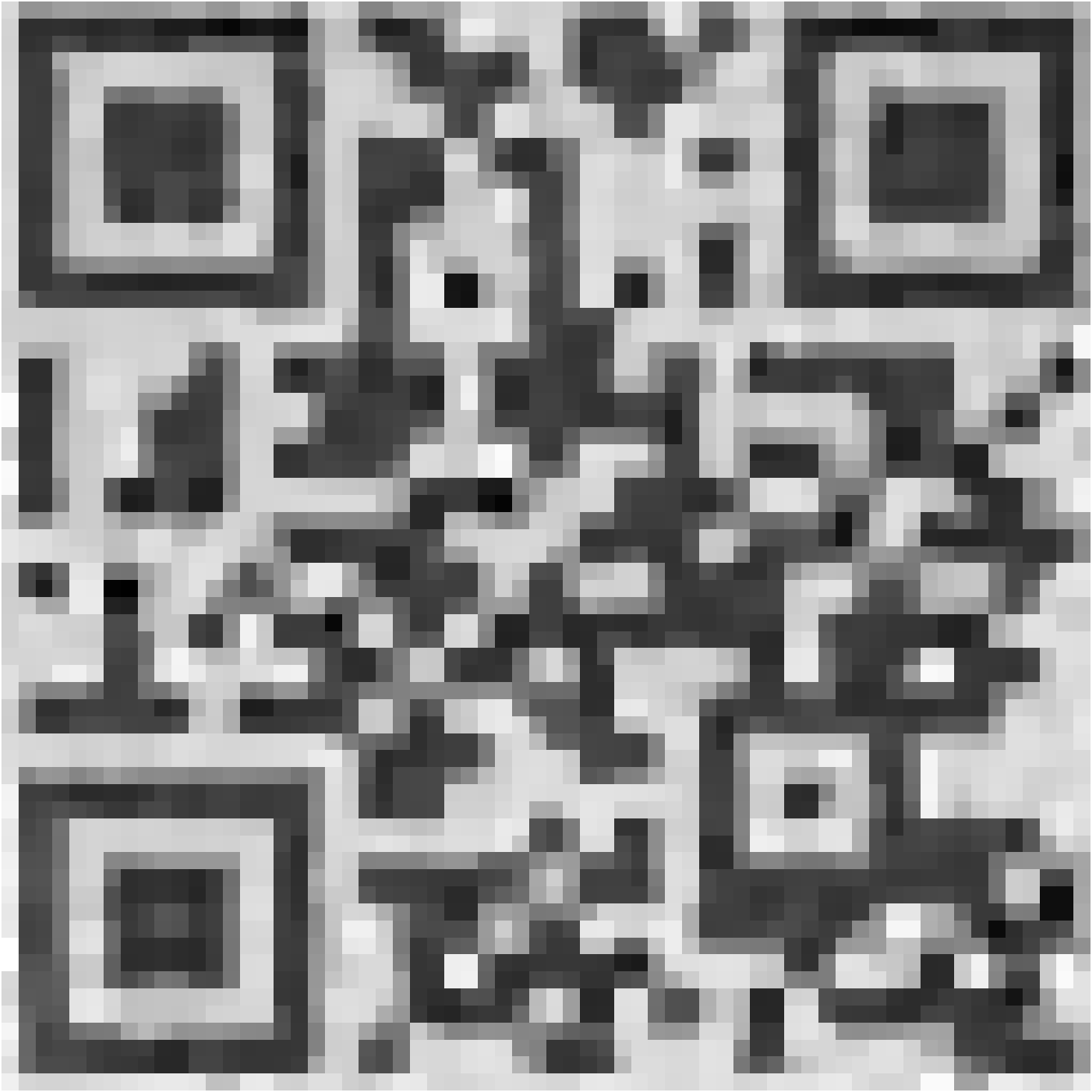}
\caption{HWT-X, $\sigma=3.5$}
\end{subfigure}\hfill
\begin{subfigure}[t]{0.23\textwidth}
\centering
\includegraphics[width=\linewidth]{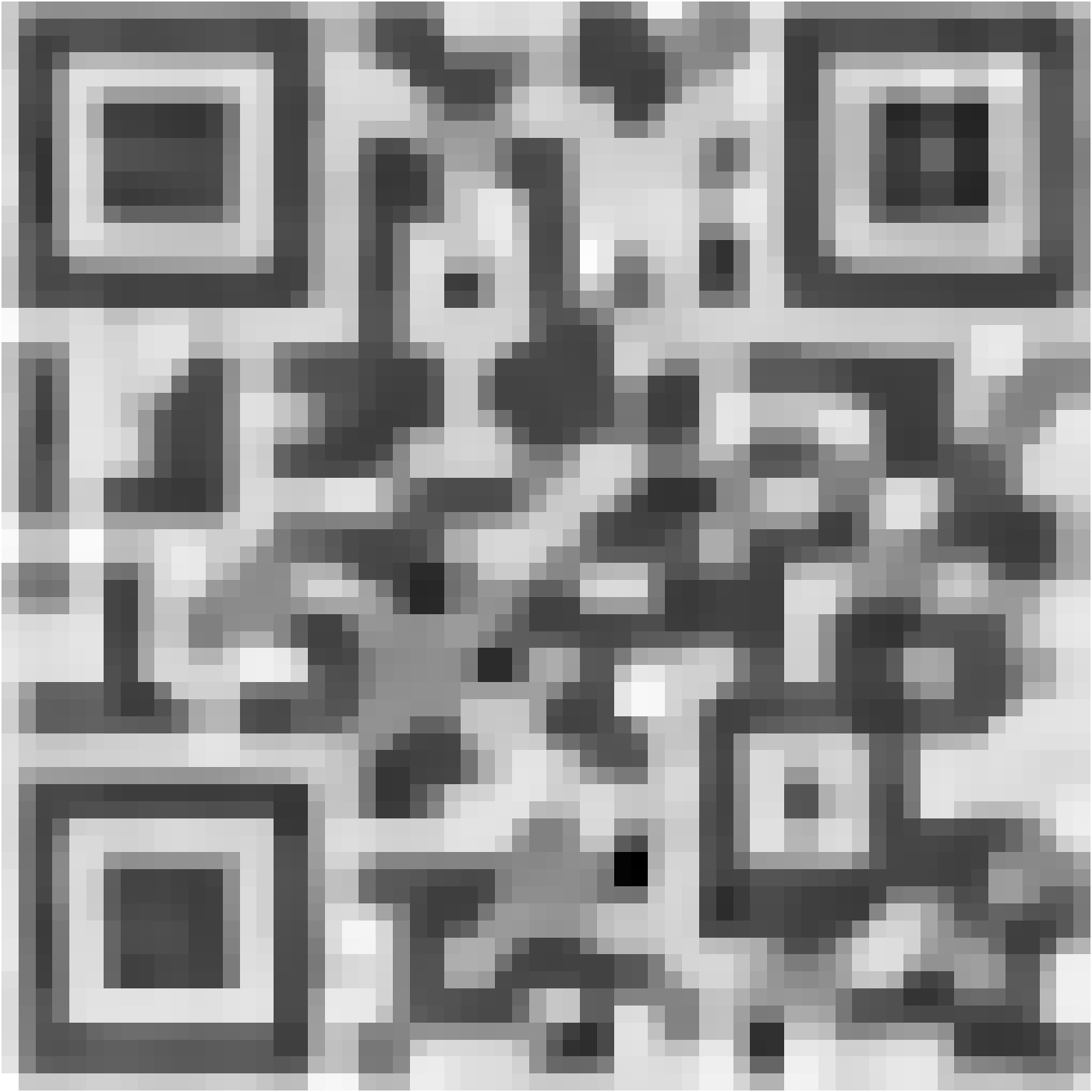}
\caption{HWT-D, $\sigma=4.5$}
\end{subfigure}\hfill
\begin{subfigure}[t]{0.23\textwidth}
\centering
\includegraphics[width=\linewidth]{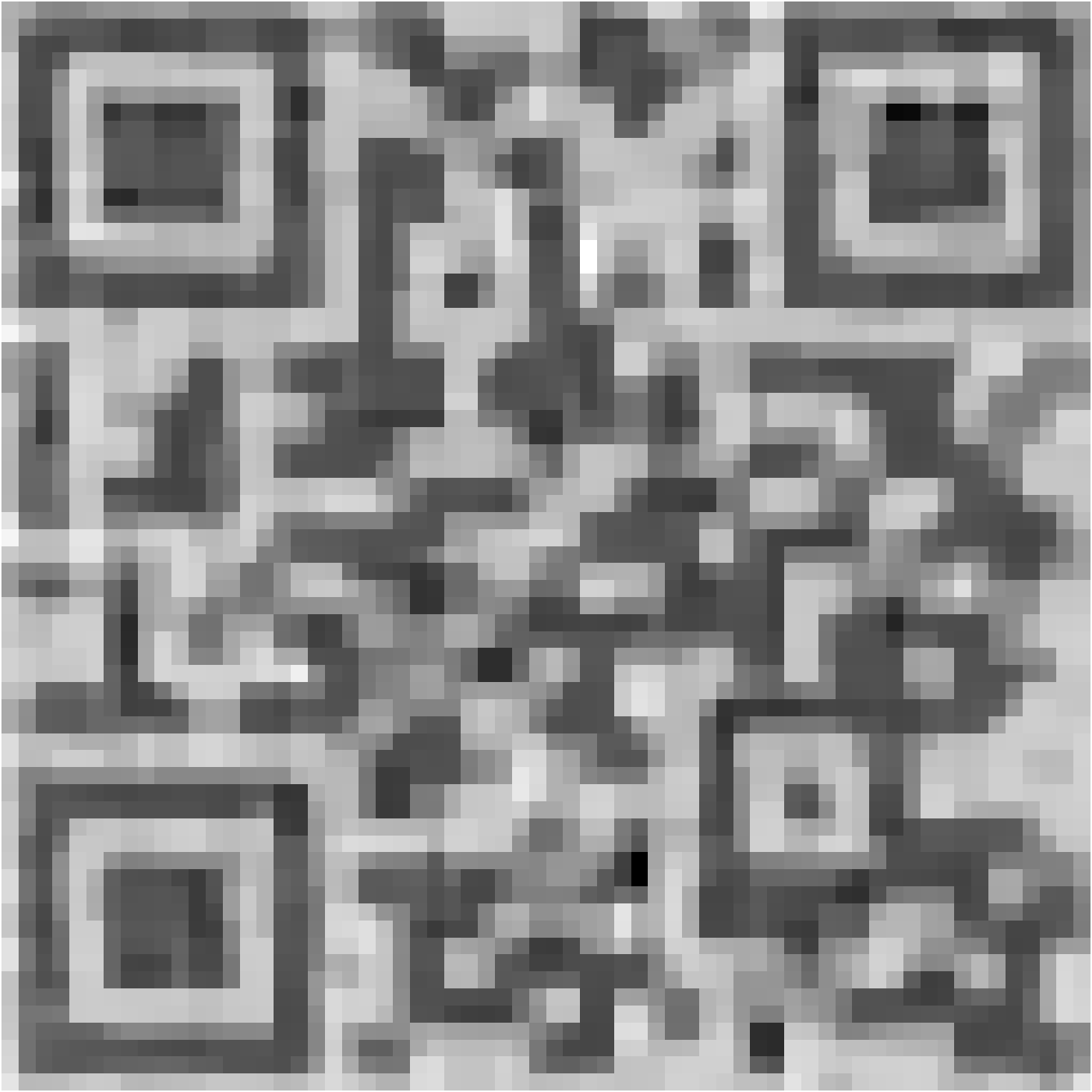}
\caption{HWT-X, $\sigma=4.5$}
\end{subfigure}
\vspace{0.8em}
\begin{subfigure}[t]{0.23\textwidth}
\centering
\includegraphics[width=\linewidth]{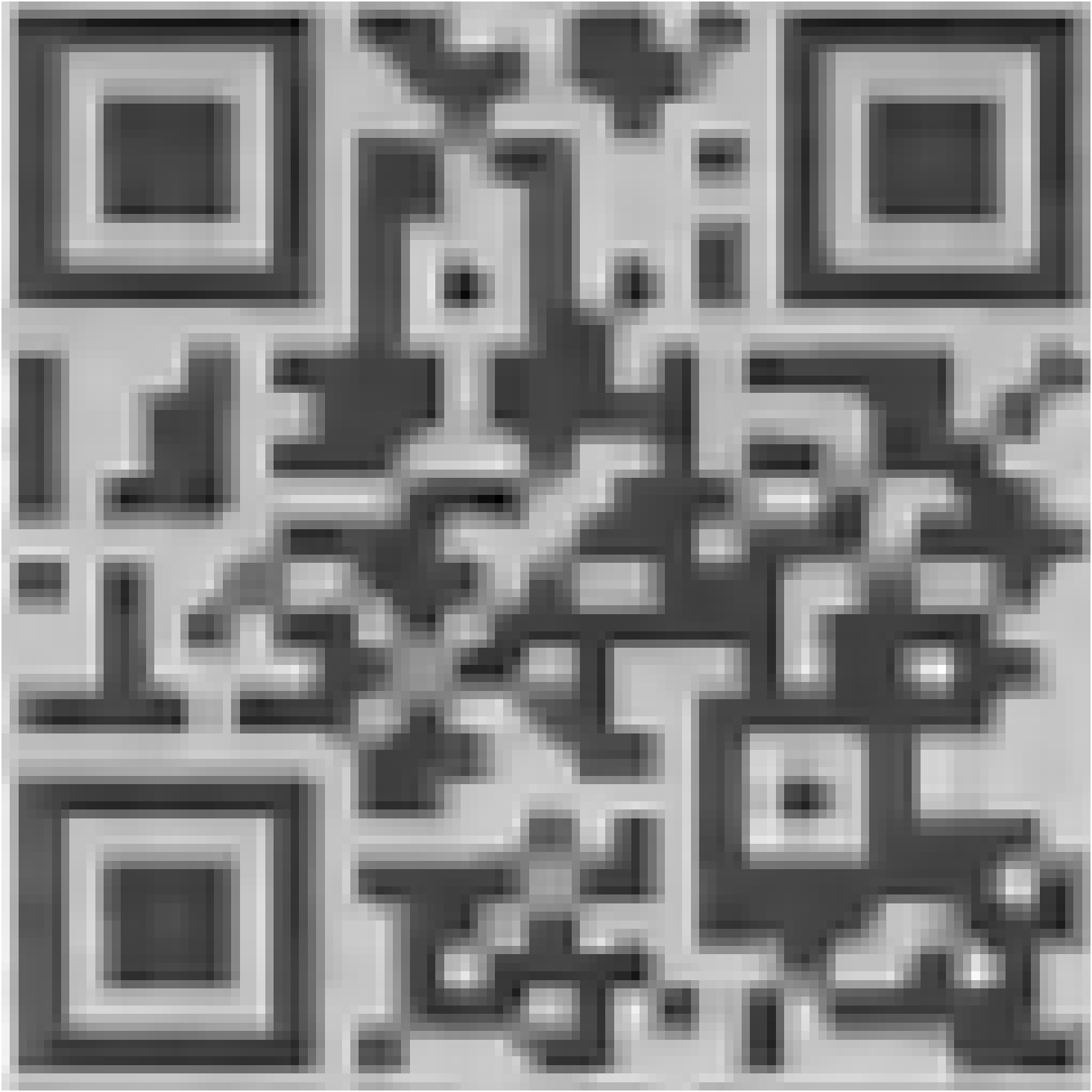}
\caption{DWT-D, $\sigma=3.5$}
\end{subfigure}\hfill
\begin{subfigure}[t]{0.23\textwidth}
\centering
\includegraphics[width=\linewidth]{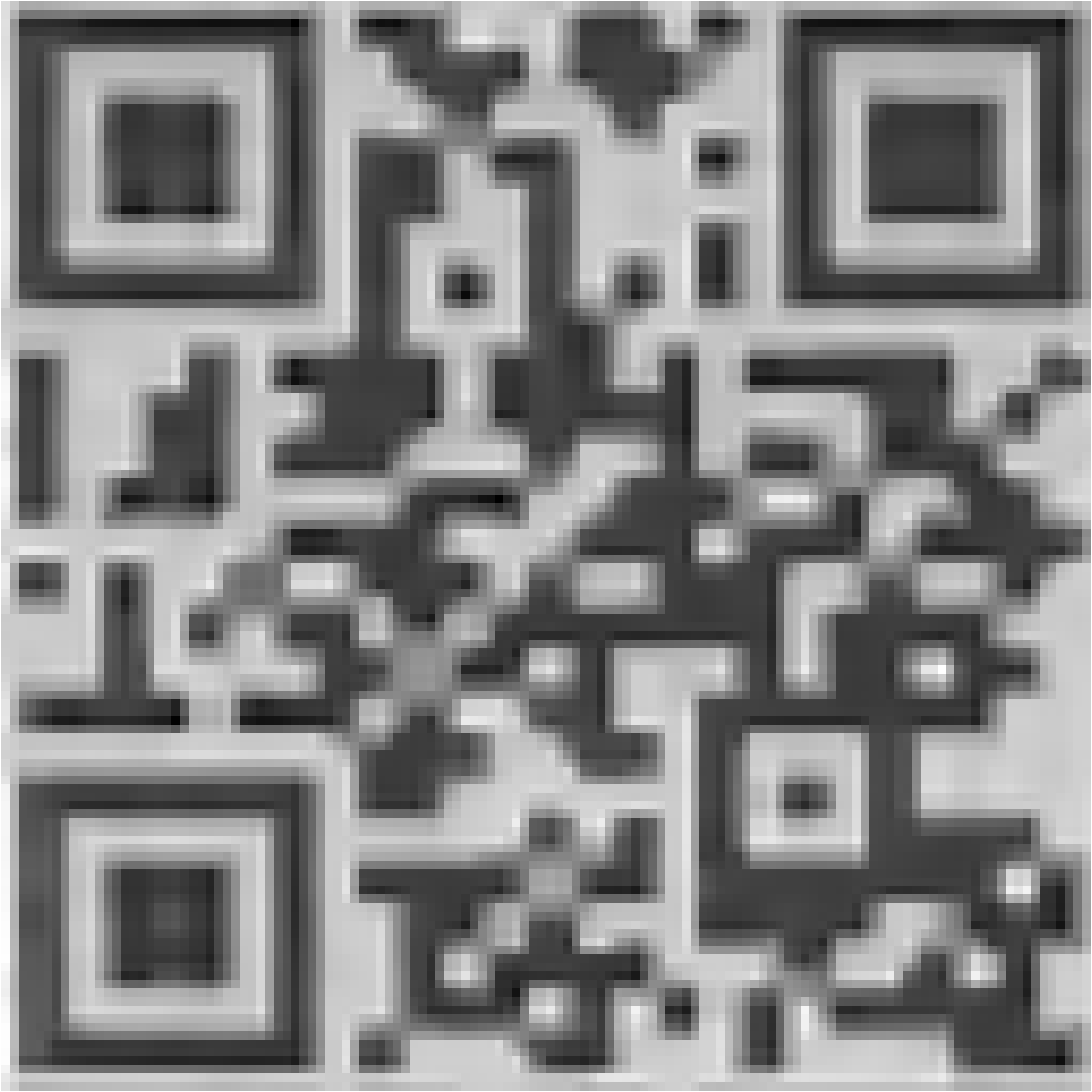}
\caption{DWT-X, $\sigma=3.5$}
\end{subfigure}\hfill
\begin{subfigure}[t]{0.23\textwidth}
\centering
\includegraphics[width=\linewidth]{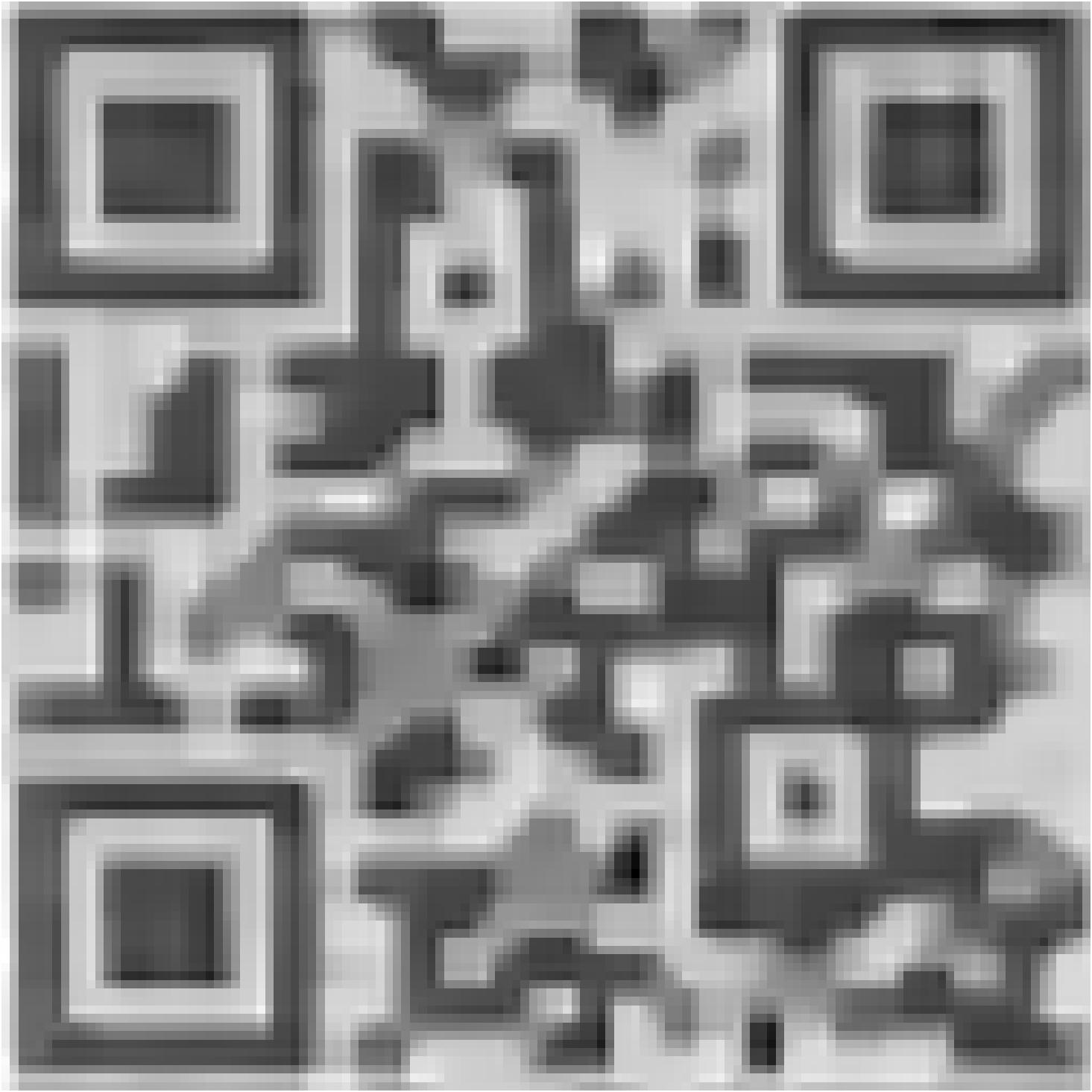}
\caption{DWT-D, $\sigma=4.5$}
\end{subfigure}\hfill
\begin{subfigure}[t]{0.23\textwidth}
\centering
\includegraphics[width=\linewidth]{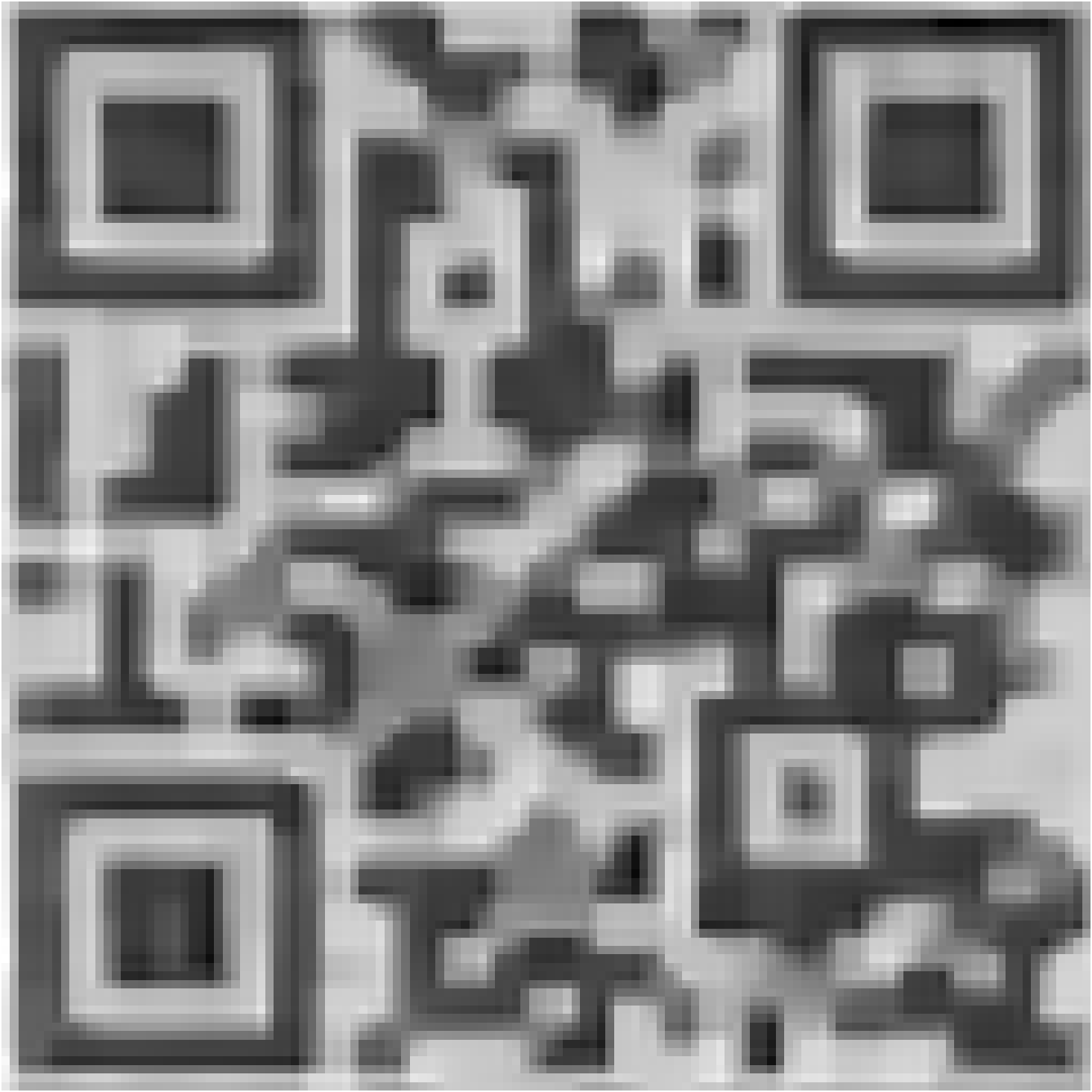}
\caption{DWT-X, $\sigma=4.5$}
\end{subfigure}
\caption{Three-level MM reconstructions. Top row: HWT; bottom row: DWT.}
\label{fig:MM3levelQ}
\end{figure}

\begin{figure}[tbp]
\centering
\begin{subfigure}[t]{0.23\textwidth}
\centering
\includegraphics[width=\linewidth]{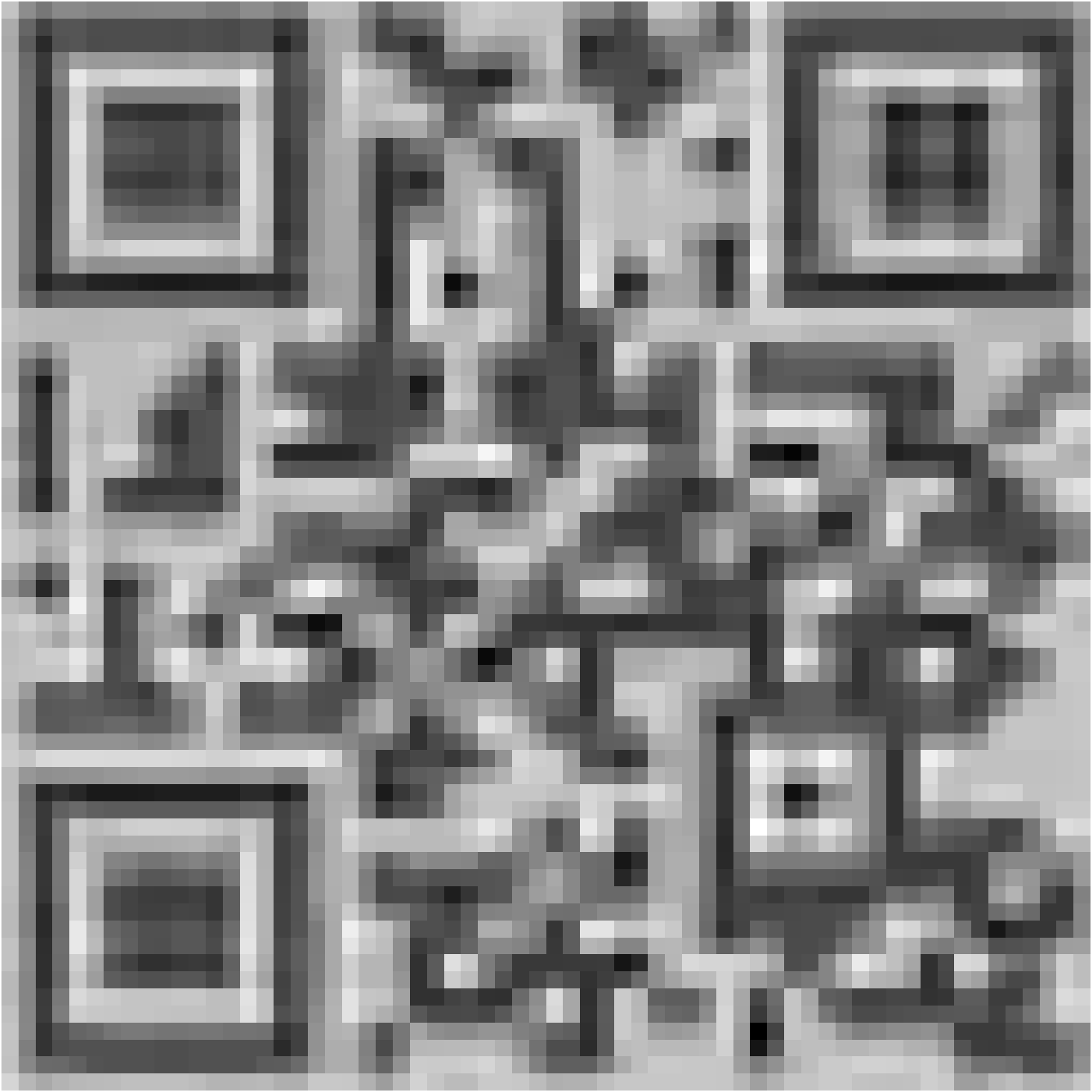}
\caption{HWT-D, $\sigma=3.5$}
\end{subfigure}\hfill
\begin{subfigure}[t]{0.23\textwidth}
\centering
\includegraphics[width=\linewidth]{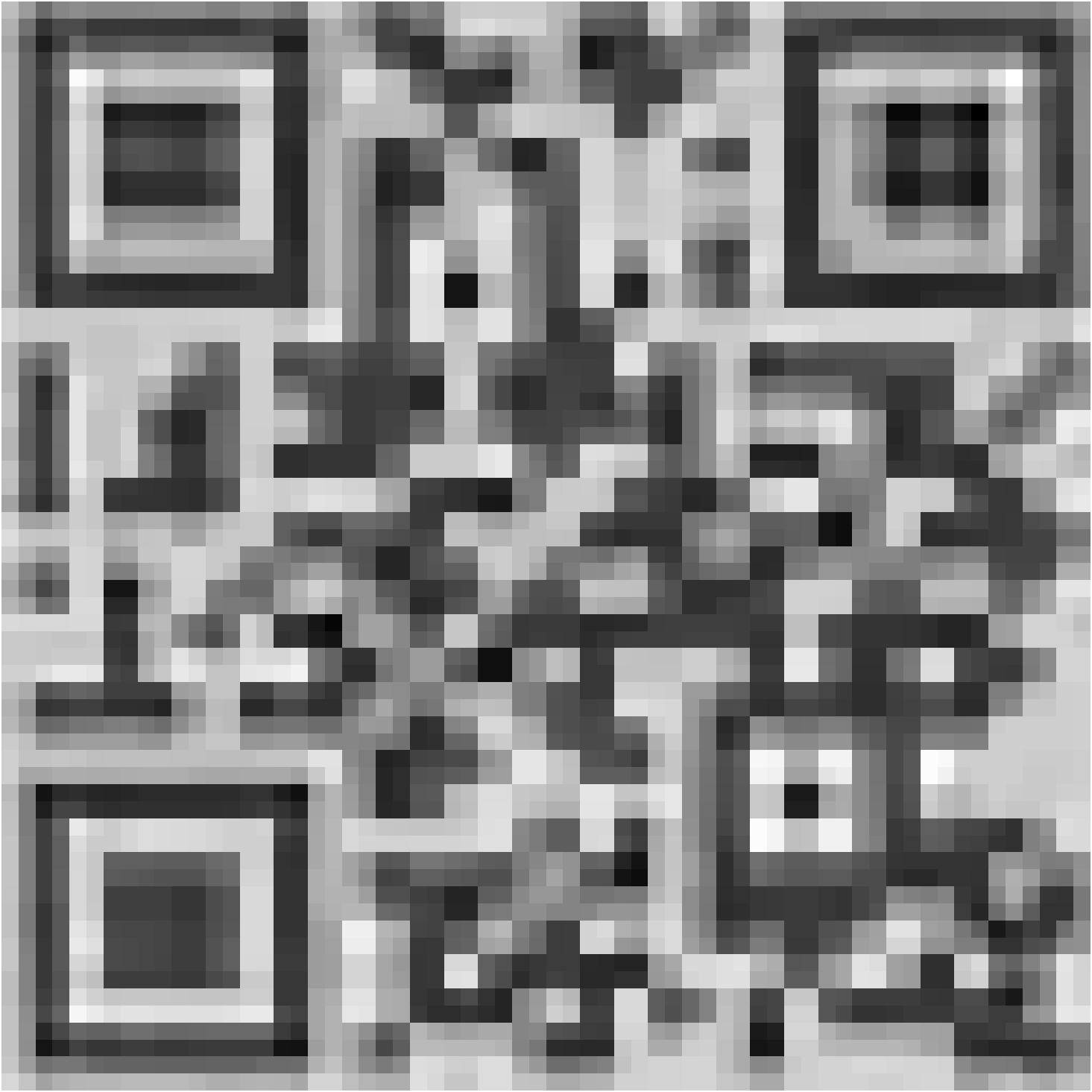}
\caption{HWT-X, $\sigma=3.5$}
\end{subfigure}\hfill
\begin{subfigure}[t]{0.23\textwidth}
\centering
\includegraphics[width=\linewidth]{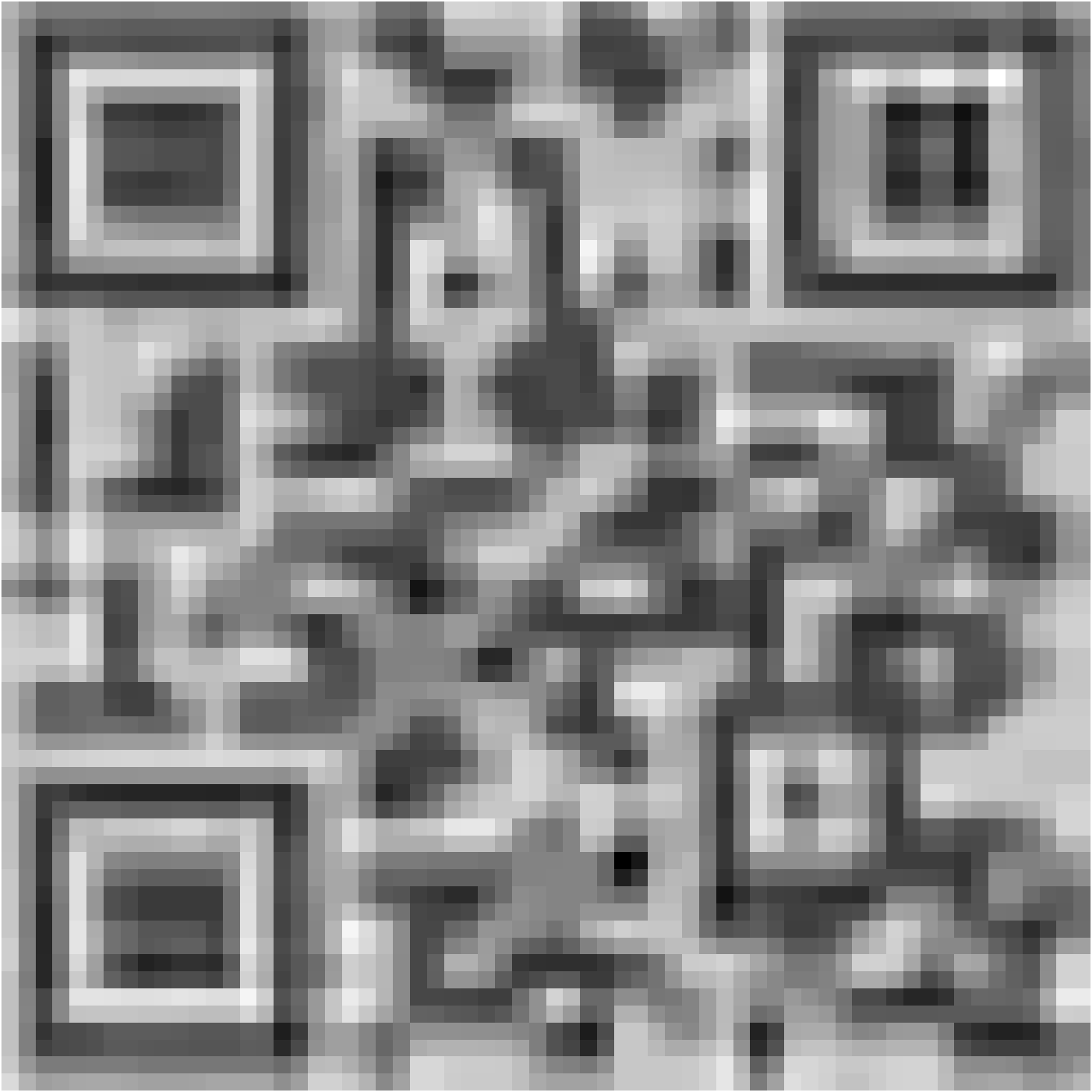}
\caption{HWT-D, $\sigma=4.5$}
\end{subfigure}\hfill
\begin{subfigure}[t]{0.23\textwidth}
\centering
\includegraphics[width=\linewidth]{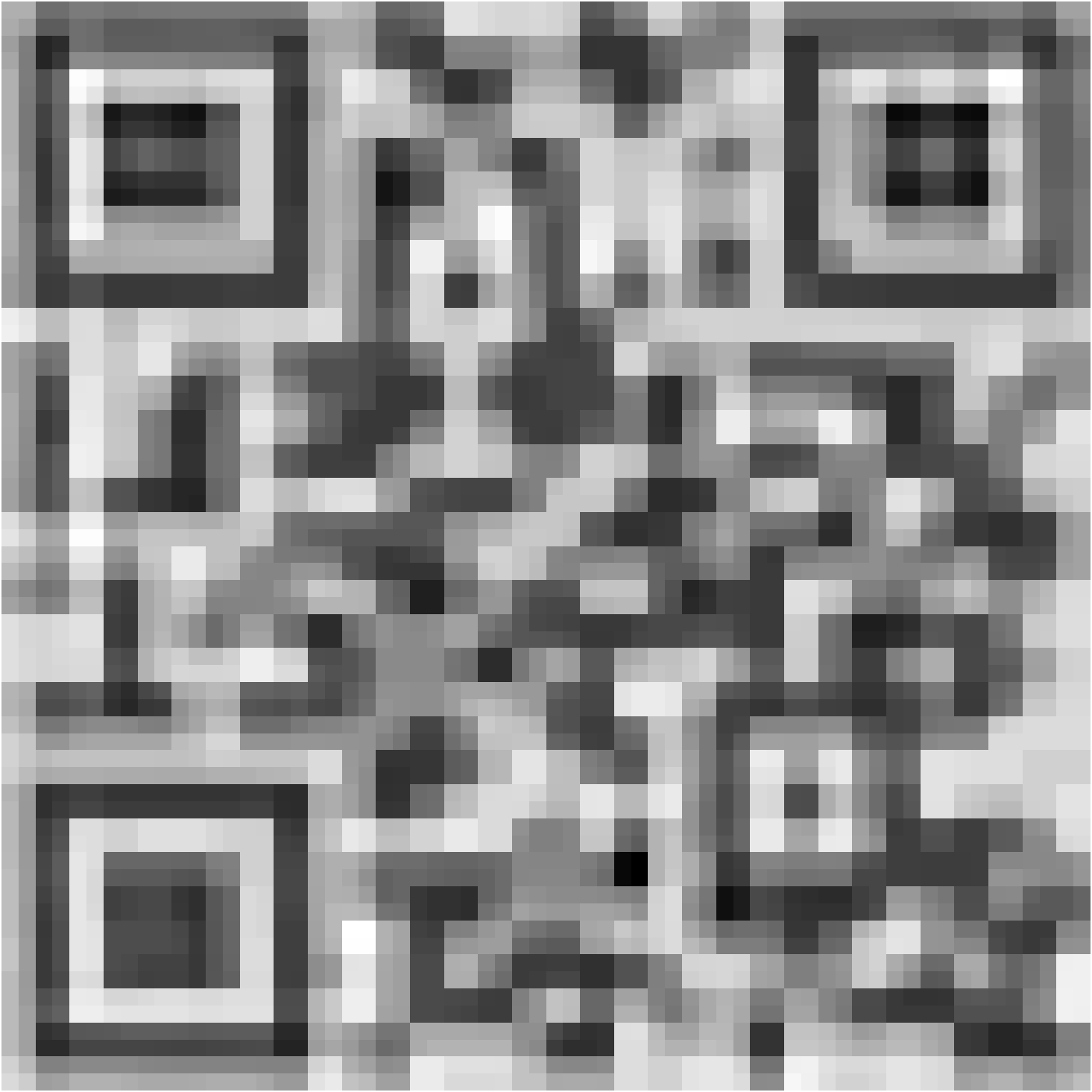}
\caption{HWT-X, $\sigma=4.5$}
\end{subfigure}
\vspace{0.8em}
\begin{subfigure}[t]{0.23\textwidth}
\centering
\includegraphics[width=\linewidth]{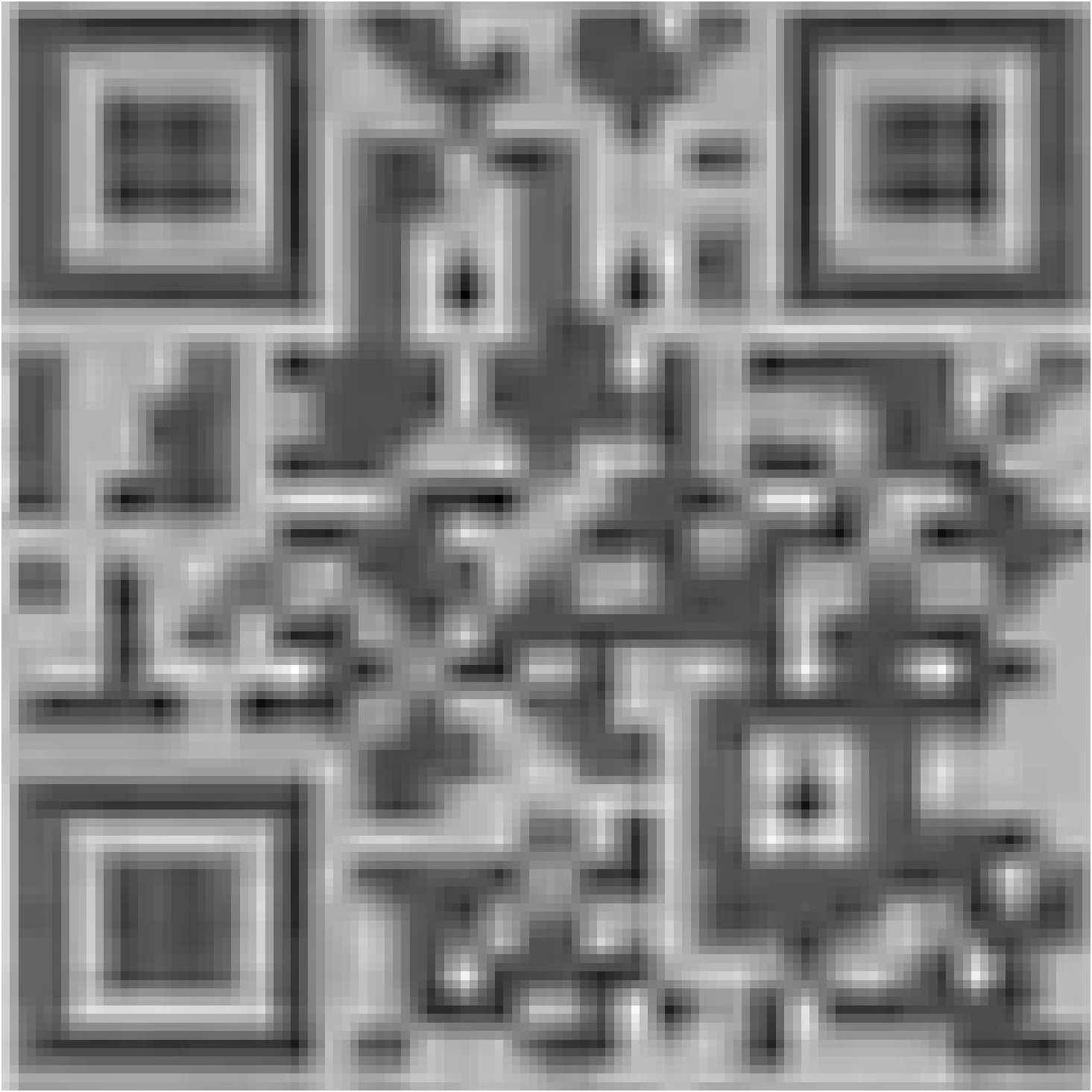}
\caption{DWT-D, $\sigma=3.5$}
\end{subfigure}\hfill
\begin{subfigure}[t]{0.23\textwidth}
\centering
\includegraphics[width=\linewidth]{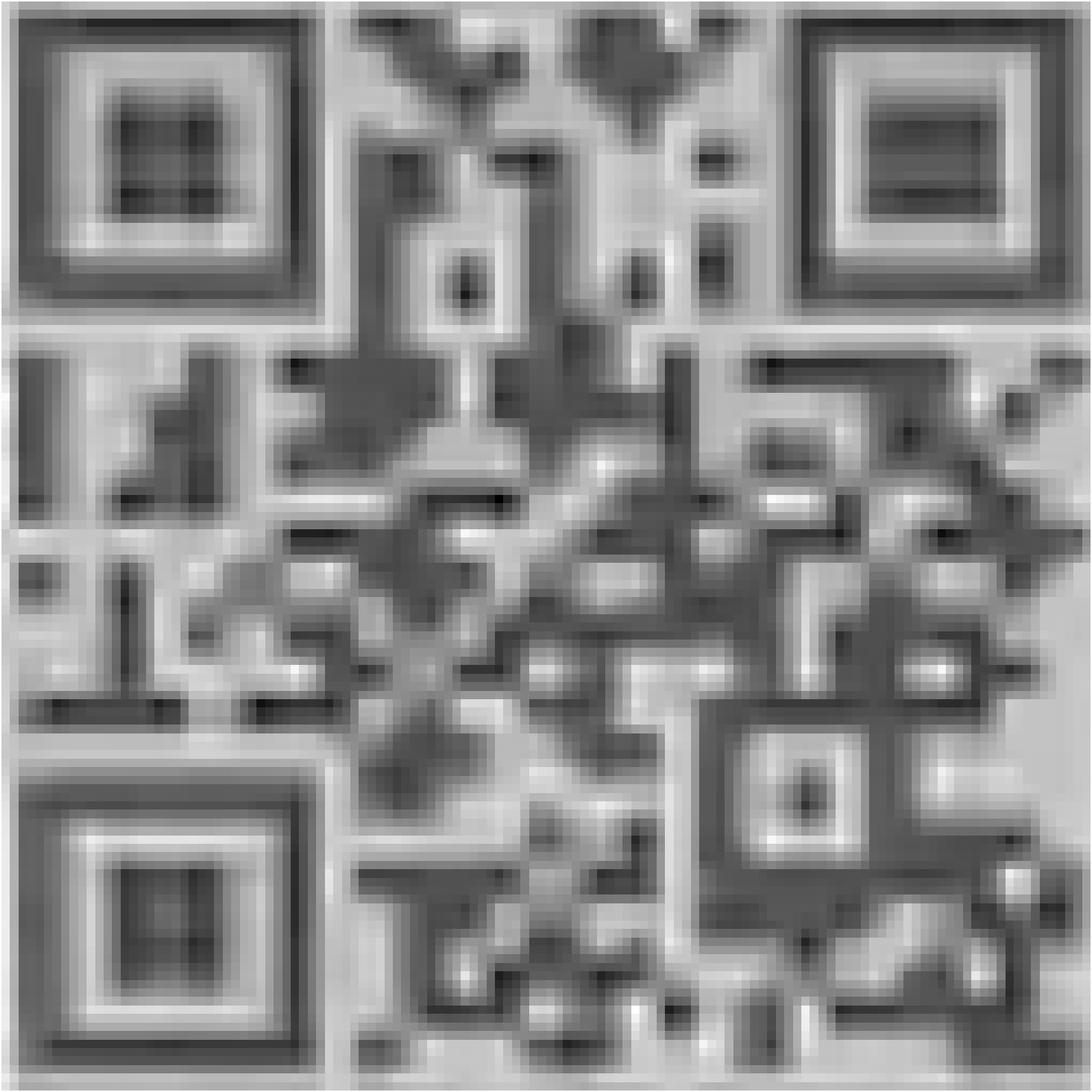}
\caption{DWT-X, $\sigma=3.5$}
\end{subfigure}\hfill
\begin{subfigure}[t]{0.23\textwidth}
\centering
\includegraphics[width=\linewidth]{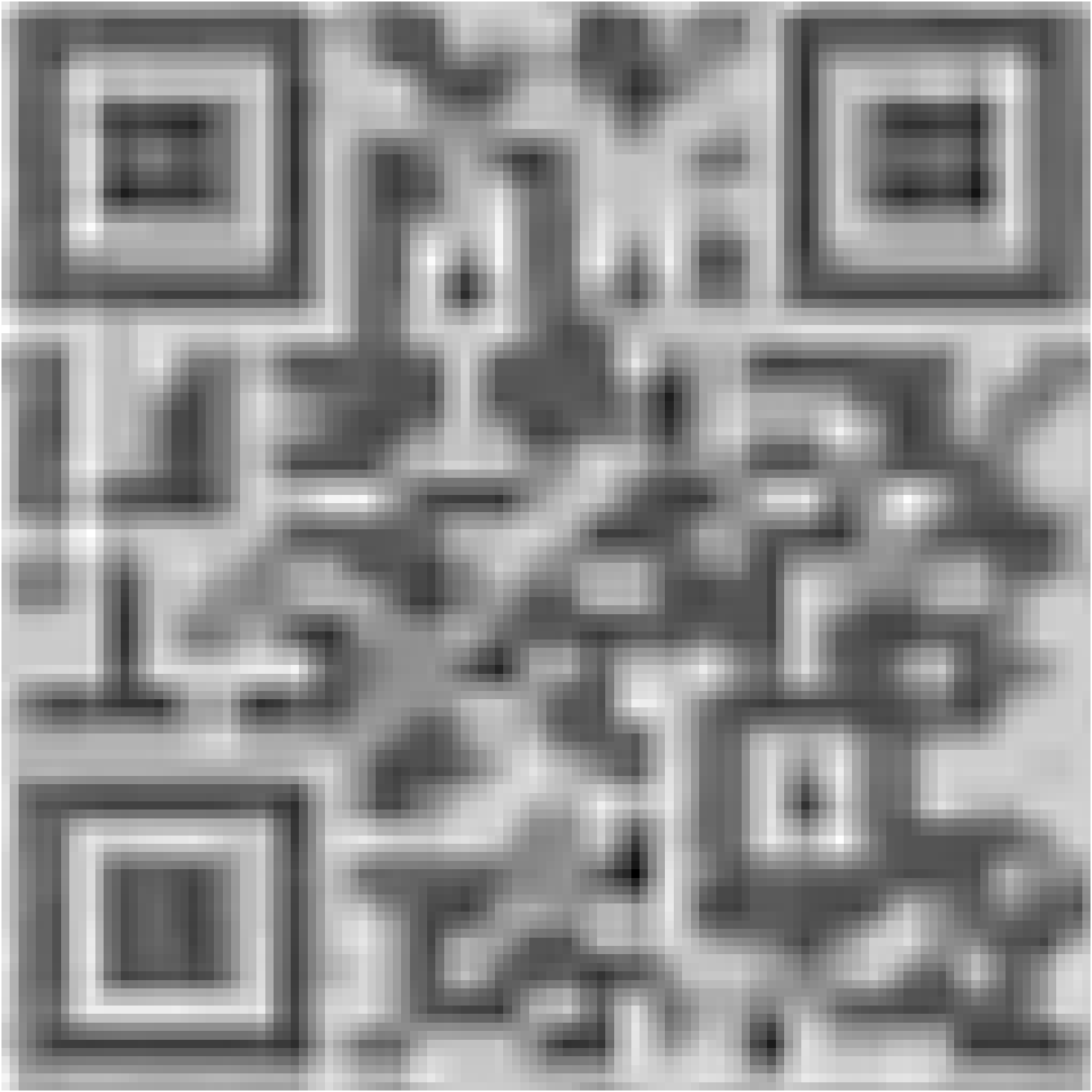}
\caption{DWT-D, $\sigma=4.5$}
\end{subfigure}\hfill
\begin{subfigure}[t]{0.23\textwidth}
\centering
\includegraphics[width=\linewidth]{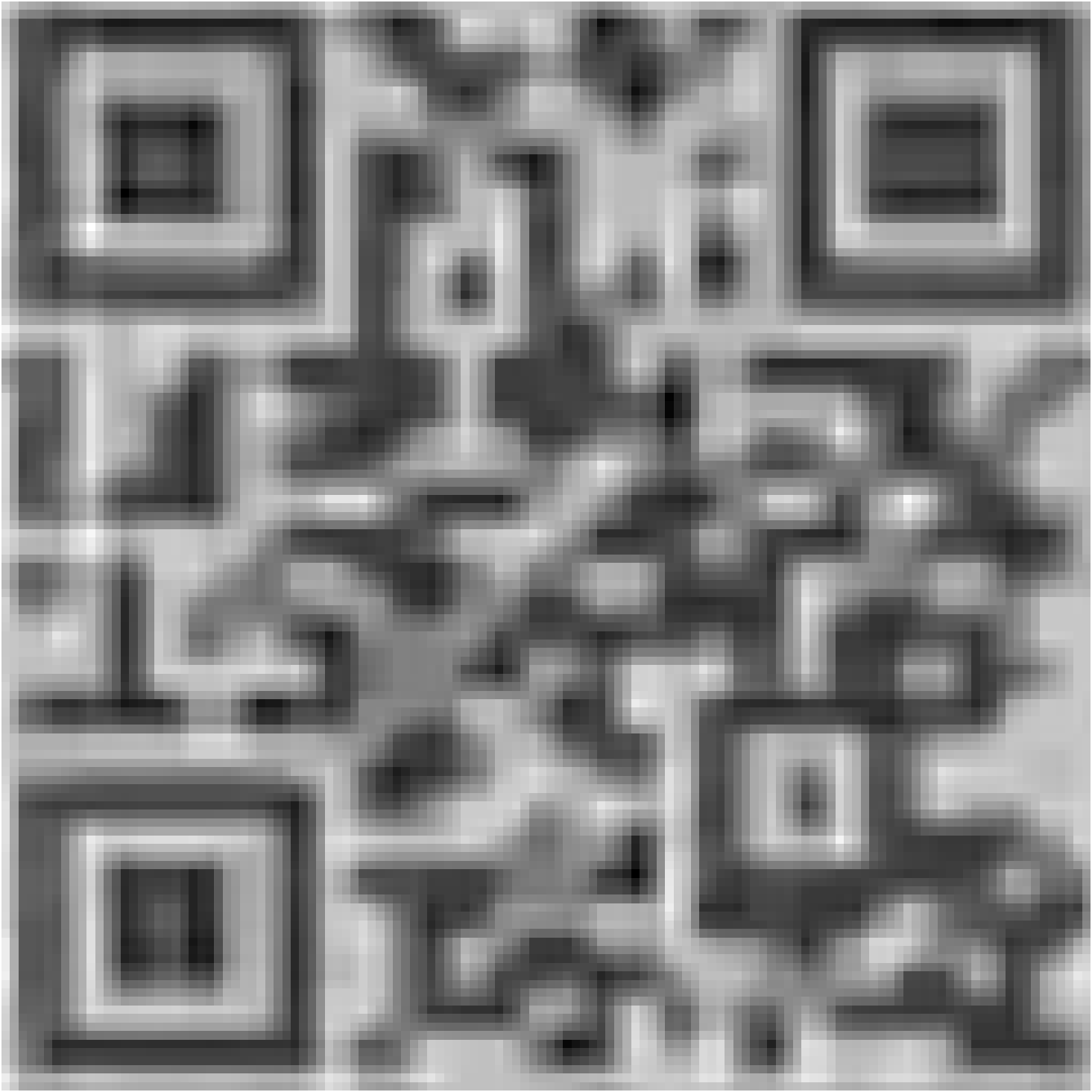}
\caption{DWT-X, $\sigma=4.5$}
\end{subfigure}
\caption{Three-level SB reconstructions. Top row: HWT; bottom row: DWT.}
\label{fig:SB3levelQ}
\end{figure}

Figures~\ref{fig:IRLS3levelQ}-\ref{fig:SB3levelQ} display the finest-level reconstructions for IRLS, MM, and SB, and Tables~\ref{tab:ssim_summary}-\ref{tab:time_summary} summarize quality metrics and runtimes. The three-level results are broadly consistent with the two-level findings. For the Voronoi image, and for the QR image as measured by SSIM (Table~\ref{tab:ssim_summary}), the qualitative interaction between wavelet family and information transfer strategy remains evident across both blur levels, although the quantitative differences become smaller for IRLS. For the QR image, however, RRE (Table~\ref{tab:rre_summary}) shows the opposite pattern under Approach~D: DWT-D attains lower RRE than HWT-D for every solver and blur level, even though HWT-D remains preferred under SSIM. A plausible explanation is that the QR image is binary and high-frequency throughout, so global $\ell_2$ pixel error is sensitive to amplitude and contrast errors that Daubechies' smoother coarse approximation may reduce, whereas SSIM's local structural comparison continues to reward the sharper edge localization that Haar provides. We do not have a definitive explanation for this metric-dependent reversal and flag it as a limitation of the study. The four multilevel combinations remain broadly comparable in runtime within each solver, and runtimes remain modest (2-15 seconds) for the $128\times128$ experiments considered here.

\begin{table}[htbp]
\caption{SSIM at the finest level (three-level method, QR image).}
\label{tab:ssim_summary}
\centering
\begin{tabular}{llcccc}
\toprule
Solver & $\sigma$ & HWT-D & HWT-X & DWT-D & DWT-X \\
\midrule
IRLS & 3.5 & 0.73 & 0.70 & 0.65 & 0.66 \\
IRLS & 4.5 & 0.56 & 0.54 & 0.54 & 0.54 \\
\midrule
MM & 3.5 & 0.54 & 0.55 & 0.53 & 0.57 \\
MM & 4.5 & 0.49 & 0.47 & 0.49 & 0.50 \\
\midrule
SB & 3.5 & 0.46 & 0.48 & 0.45 & 0.49 \\
SB & 4.5 & 0.42 & 0.43 & 0.41 & 0.43 \\
\bottomrule
\end{tabular}
\end{table}

\begin{table}[htbp]
\caption{RRE at the finest level (three-level method, QR image).}
\label{tab:rre_summary}
\centering
\begin{tabular}{llcccc}
\toprule
Solver & $\sigma$ & HWT-D & HWT-X & DWT-D & DWT-X \\
\midrule
IRLS & 3.5 & 0.26 & 0.28 & 0.25 & 0.25 \\
IRLS & 4.5 & 0.31 & 0.33 & 0.30 & 0.30 \\
\midrule
MM & 3.5 & 0.32 & 0.34 & 0.27 & 0.26 \\
MM & 4.5 & 0.34 & 0.37 & 0.31 & 0.31 \\
\midrule
SB & 3.5 & 0.37 & 0.38 & 0.31 & 0.31 \\
SB & 4.5 & 0.38 & 0.39 & 0.34 & 0.33 \\
\bottomrule
\end{tabular}
\end{table}

\begin{table}[htbp]
\caption{Reconstruction time in seconds (three-level method, QR image).}
\label{tab:time_summary}
\centering
\begin{tabular}{llcccc}
\toprule
Solver & $\sigma$ & HWT-D & HWT-X & DWT-D & DWT-X \\
\midrule
IRLS & 3.5 & 7 & 6 & 8 & 7 \\
IRLS & 4.5 & 7 & 7 & 12 & 10 \\
\midrule
MM & 3.5 & 3 & 2 & 4 & 4 \\
MM & 4.5 & 5 & 6 & 8 & 7 \\
\midrule
SB & 3.5 & 2 & 2 & 2 & 2 \\
SB & 4.5 & 3 & 2 & 4 & 4 \\
\bottomrule
\end{tabular}
\end{table}
\section{Conclusions}\label{sec:conclusion}

We have introduced wavelet-based multilevel V-cycle methods for $\ell_1$-regularized image deblurring that embed three solvers (IRLS, MM, and SB) with automatic GCV parameter selection and two inter-level information transfer strategies. The numerical experiments support three main conclusions. First, multilevel variants provide substantial speedups over their single-level counterparts, most dramatically for IRLS, where multilevel initialization reduces reconstruction time by more than $15\times$; for MM and SB, where single-level runtimes are already small, the relative speedups are more modest and occasionally slightly negative, but wavelet-based coarsening remains an effective strategy for this class of problems. Second, since the four multilevel combinations (HWT-D, HWT-X, DWT-D, and DWT-X) have broadly comparable runtimes within each solver, the practical choice among them is determined primarily by reconstruction quality rather than computational cost. Third, and most importantly, the wavelet family and the information transfer strategy interact in a way that is evident across most of the experiments considered. Approach~D (auxiliary transfer) generally performs best with the Haar wavelet, whereas Approach~X (solution transfer only) generally performs best with the Daubechies wavelet. A plausible explanation is that Haar's piecewise-constant coarse representation localizes edges more sharply, yielding informative auxiliary warm-starts, whereas the higher approximation order of the Daubechies wavelet produces a more accurate interpolated solution when only the coarse image is transferred. This interaction holds consistently for the Voronoi test image and for SSIM on the QR test image; for the three-level QR image, RRE reverses the pattern under Approach~D (Table~\ref{tab:rre_summary}), suggesting that the interaction may also depend on image content and evaluation metric. This interaction is supported by the numerical experiments considered here, while its theoretical explanation and the source of this metric-dependent exception remain open questions. Future work includes investigating other wavelet families, incorporating alternative hybrid projection techniques within the multilevel framework, and applying the proposed methodology to other inverse problems.

\section*{Acknowledgments}
Malena I. Espa\~nol acknowledges support from the Karen EDGE Fellowship.
Misha E. Kilmer acknowledges support from NSF DMS-2410698.

\bibliographystyle{plain}
\bibliography{references}
\end{document}